\documentclass{amsart}
\usepackage{amsthm,amsmath,amsfonts}
\usepackage{mathrsfs}
\usepackage[colorlinks=false,linktocpage=true]{hyperref}
\usepackage{tocvsec2}
\usepackage{hyperref}
\usepackage{framed}
\usepackage{color}
\usepackage[usenames,dvipsnames,svgnames]{xcolor}
\numberwithin{equation}{section}
\numberwithin{table}{section}
\font\tenscrpt=eusm10
\font\sevenscrpt=eusm10 scaled 700
\font\fivescrpt=eusm10 scaled 500
\newfam\eusmfam
\textfont\eusmfam=\tenscrpt
\scriptfont\eusmfam=\sevenscrpt
\scriptscriptfont\eusmfam=\fivescrpt

\newtheorem{thm}{Theorem}[section]
\newtheorem{cor}{Corollary}[section]
\newtheorem{lem}{Lemma}[section]
\newtheorem{prop}{Proposition}[section]

\theoremstyle{definition}
\newtheorem{defn}{Definition}[section]

\newtheorem{rem}{Remark}[section]
\newtheorem{notn}{Notation}[section]
\newcommand{\thmref}[1]{Theorem~\ref{#1}}
\newcommand{\secref}[1]{Section~\ref{#1}}

\newcommand{\lemref}[1]{Lemma~\ref{#1}}
\newcommand{\coref}[1]{Corollary~\ref{#1}}
\newcommand{\propref}[1]{Proposition~\ref{#1}}
\newcommand{\defnref}[1]{Definition~\ref{#1}}
\newcommand{\remref}[1]{Remark~\ref{#1}}
\newcommand{\eqnref}[1]{{\rm (\ref{#1})}}

\newcommand{\ceil}[1]{\lceil{#1}\rceil}

\def\qed{\quad\vcenter{\hrule\hbox{\vrule height.6em\kern.6em\vrule}\hrule}}
\newenvironment{pf}{{\bigskip\textit{\newline Proof.}\quad}}{$\qed$\bigskip\newline}
\newenvironment{pf*}[1]{{\bigskip\textit{\newline#1.}\quad}}{$\qed$\bigskip\newline}
\def\ds{\displaystyle}
\def\ts{\textstyle}

\def\h{\delta_n}

\def\hd{\delta_n^d}
\def\sdh{\delta_n^{{d}/{2}}}

\def\d{\delta}

\def\dd{\delta^d}

\def\sdd{\delta^{{d}/{2}}}

\def\ptxy{{K}^{\text{\tiny{\sc{BM}}}^d}_{t;x,y}}
\def\psxy{{K}^{\text{\tiny{\sc{BM}}}^d}_{s;x,y}}

\def\psoxz{{K}^{\text{\tiny{\sc{BM}}}^d}_{s_{1};x}}
\def\puoxz{{K}^{\text{\tiny{\sc{BM}}}^d}_{u_{1};x}}

\def\ptzs{{K}^{\text{\tiny{\sc{BM}}}}_{t;0,s}}
\def\ptzz{{K}^{\text{\tiny{\sc{BM}}}}_{t;0,0}}

\def\ptzsk{{K}^{\text{\tiny{\sc{BM}}}}_{t;0,\frac{s_{k}}{\sqrt2}}}
\def\ptzuk{{K}^{\text{\tiny{\sc{BM}}}}_{t;0,\frac{u_{k}}{\sqrt2}}}
\def\pszuk{{K}^{\text{\tiny{\sc{BM}}}}_{s;0,\frac{u_{k}}{\sqrt2}}}
\def\pszrk{{K}^{\text{\tiny{\sc{BM}}}}_{s;0,\frac{r_{k}}{\sqrt2}}}
\def\ptzskiimo{{K}^{\text{\tiny{\sc{BM}}}}_{s_{k-i};0,\frac{s_{k-i-1}}{\sqrt2}}}
\def\ptzrkiimo{{K}^{\text{\tiny{\sc{BM}}}}_{r_{k-i};0,\frac{r_{k-i-1}}{\sqrt2}}}
\def\ptzukiimo{{K}^{\text{\tiny{\sc{BM}}}}_{u_{k-i};0,\frac{u_{k-i-1}}{\sqrt2}}}

\def\proprtz{K^{{\text{\tiny{\sc{BM}}}^d}}_{r_{1}+r_{2};0}}
\def\prox{K^{{\text{\tiny{\sc{BM}}}^d}}_{r_{1};x}}
\def\prtx{K^{{\text{\tiny{\sc{BM}}}^d}}_{r_{2};x}}
\def\prix{K^{{\text{\tiny{\sc{BM}}}^d}}_{r_{i};x}}
\def\prixpz{K^{{\text{\tiny{\sc{BM}}}^d}}_{r_{i};x+z}}
\def\proprtz{K^{{\text{\tiny{\sc{BM}}}^d}}_{r_{1}+r_{2};0}}
\def\propuoz{K^{{\text{\tiny{\sc{BM}}}^d}}_{r_{1}+u_{1};0}}

\def\propuozz{K^{{\text{\tiny{\sc{BM}}}^d}}_{r_{1}+u_{1};z}}

\def\KItzs{{K}^{{\Lambda_{\beta}}}_{t;0,s}}

\def\KItzso{{K}^{{\Lambda_{\beta}}}_{t;0,s_{1}}}
\def\KItzuo{{K}^{{\Lambda_{\beta}}}_{t;0,u_{1}}}
\def\KItzso{K^{{\Lambda_{\beta}}}_{t;0,s_{1}}}

\def\KIszri{K^{{\Lambda_{\beta}}}_{s;0,r_{i}}}
\def\KIszro{K^{{\Lambda_{\beta}}}_{s;0,r_{1}}}
\def\KIszuo{K^{{\Lambda_{\beta}}}_{s;0,u_{1}}}

\def\KIuzro{K^{{\Lambda_{\beta}}}_{u;0,r_{1}}}
\def\KIuzrt{K^{{\Lambda_{\beta}}}_{u;0,r_{2}}}
\def\KIvzrt{K^{{\Lambda_{\beta}}}_{v;0,r_{2}}}

\def\KItzx{K^{\Lambda_{\beta}}_{t;0,x}}
\def\KIhalftzx{K^{\Lambda_{\frac12}}_{t;0,x}}
\def\KIhalftzxsqrt2{K^{\Lambda_{\frac12}}_{t;0,\sqrt2 x}}
\def\KabBmtzxsqrt2{K^{{|B|}}_{t;0,\frac{x}{\sqrt{2}}}}

\def\K{{\mathbb K}}

\def\KBtxy{{\K}^{\text{\tiny{\sc{BTBM}}}^d}_{t;x,y}}
\def\KIBtxy{{\K}^{\text{\tiny{\sc{BM}}}^d,\Lambda_\beta}_{t;x,y}}
\def\KBtx{{\K}^{\text{\tiny{\sc{BTBM}}}^d}_{t;x}}
\def\KIBtx{{\K}^{\text{\tiny{\sc{BM}}}^d,\Lambda_\beta}_{t;x}}

\def\KBtsxy{{\K}^{\text{\tiny{\sc{BTBM}}}^d}_{t-s;x,y}}
\def\KIBtsxy{{\K}^{\text{\tiny{\sc{BM}}}^d,\Lambda_\beta}_{t-s;x,y}}

\def\KIBtsx{{\K}^{\text{\tiny{\sc{BM}}}^d,\Lambda_\beta}_{t-s;x}}

\def\KIBrsx{{\K}^{\text{\tiny{\sc{BM}}}^d,\Lambda_\beta}_{r-s;x}}

\def\KIBsptmrx{{\K}^{\text{\tiny{\sc{BM}}}^d,\Lambda_\beta}_{s+(t-r);x}}

\def\KIBsx{{\K}^{{\text{\tiny{\sc{ISLTBM}}}^d}}_{s;x}}
\def\KIBsx{{\K}^{\text{\tiny{\sc{BM}}}^d,\Lambda_\beta}_{s;x}}

\def\KIBsxpz{{\K}^{\text{\tiny{\sc{BM}}}^d,\Lambda_\beta}_{s;x+z}}

\def\KIBux{{\K}^{\text{\tiny{\sc{BM}}}^d,\Lambda_\beta}_{u;x}}

\def\Ktx{{\K}^{{\text{\tiny{\sc{RW}}}^d_\delta,\Lambda_\beta}}_{t;x}}
\def\Ktxy{{\K}^{\text{\tiny{\sc{RW}}}^d_\delta,\Lambda_\beta}_{t;x,y}}
\def\Ktsx{{\K}^{{\text{\tiny{\sc{RW}}}^d_\delta,\Lambda_\beta}}_{t-s;x}}
\def\Krsx{{\K}^{{\text{\tiny{\sc{RW}}}^d_\delta,\Lambda_\beta}}_{r-s;x}}
\def\Ksx{{\K}^{{\text{\tiny{\sc{RW}}}^d_\delta,\Lambda_\beta}}_{s;x}}
\def\Kux{{\K}^{{\text{\tiny{\sc{RW}}}^d_\delta,\Lambda_\beta}}_{u;x}}
\def\Kvx{{\K}^{{\text{\tiny{\sc{RW}}}^d_\delta,\Lambda_\beta}}_{v;x}}
\def\Ksptmrx{{\K}^{{\text{\tiny{\sc{RW}}}^d_\delta,\Lambda_\beta}}_{s+(t-r);x}}

\def\Ktsxy{{\K}^{{\text{\tiny{\sc{RW}}}^d_\delta,\Lambda_\beta}}_{t-s;x,y}}
\def\Ktsxz{{\K}^{{\text{\tiny{\sc{RW}}}^d_\delta,\Lambda_\beta}}_{t-s;x,z}}
\def\Ktsyz{{\K}^{{\text{\tiny{\sc{RW}}}^d_\delta,\Lambda_\beta}}_{t-s;y,z}}
\def\Krsxz{{\K}^{{\text{\tiny{\sc{RW}}}^d_\delta,\Lambda_\beta}}_{r-s;x,z}}
\def\Kroxz{{\K}^{{\text{\tiny{\sc{RW}}}^d_\delta,\Lambda_\beta}}_{\rho;x,z}}
\def\Ksx{{\K}^{{\text{\tiny{\sc{RW}}}^d_\delta,\Lambda_\beta}}_{s;x}}
\def\Ksxpz{{\K}^{{\text{\tiny{\sc{RW}}}^d_\delta,\Lambda_\beta}}_{s;x+z}}
\def\Kzxn{{\K}^{\text{\tiny{\sc{RW}}}^d_{\h},\Lambda_\beta}_{0;x}}

\def\Ktxn{{\K}^{\text{\tiny{\sc{RW}}}^d_{\h},\Lambda_\beta}_{t,x}}
\def\Ktxzn{{\K}^{\text{\tiny{\sc{RW}}}^d_{\h},\Lambda_\beta}_{t;x,0}}
\def\Ktsxyn{{\K}^{\text{\tiny{\sc{RW}}}^d_{\h},\Lambda_\beta}_{t-s;x,y}}

\def\Ktausxyn{{\K}^{\text{\tiny{\sc{RW}}}^d_{\h},\Lambda_\beta}_{\tau-s;x,y}}

\def\Ktauitxiyn{{\K}^{\text{\tiny{\sc{RW}}}^d_{\h},\Lambda_\beta}_{\tau_i-t;x_i,y}}

\def\Ktaujtxjyn{{\K}^{\text{\tiny{\sc{RW}}}^d_{\h},\Lambda_\beta}_{\tau_j-t;x_j,y}}
\def\Ktxyn{{\K}^{\text{\tiny{\sc{RW}}}^d_{\h},\Lambda_\beta}_{t;x,y}}
\def\Ktxynint{{\K}^{\text{\tiny{\sc{RW}}}^d_{\h},\Lambda_\beta}_{t;[x]_{\h},[y]_{\h}}}

\def\KIIsptmrz{{\K}^{\text{\tiny{\sc{RW}}}^d_{\h},2\Lambda_\beta}_{s+(t-r),s+(t-r);0}}

\def\KIIBsptmrz{{\K}^{\text{\tiny{\sc{BM}}}^d,2\Lambda_\beta}_{s+(t-r),s+(t-r);0}}

\def\KIIsptmrsz{{\K}^{\text{\tiny{\sc{RW}}}^d_{\h},2\Lambda_\beta}_{s+(t-r),s;0}}
\def\KIIBsptmrsz{{\K}^{\text{\tiny{\sc{BM}}}^d,2\Lambda_\beta}_{s+(t-r),s;0}}

\def\KIIsz{{\K}^{\text{\tiny{\sc{RW}}}^d_{\h},2\Lambda_\beta}_{s,s;0}}
\def\KIIBsz{{\K}^{\text{\tiny{\sc{BM}}}^d,2\Lambda_\beta}_{s,s;0}}

\def\KIIuvz{{\K}^{\text{\tiny{\sc{RW}}}^d_{\h},2\Lambda_\beta}_{u,v;0}}
\def\KIIBuvz{{\K}^{\text{\tiny{\sc{BM}}}^d,2\Lambda_\beta}_{u,v;0}}
\def\KIIBuz{{\K}^{\text{\tiny{\sc{BM}}}^d,2\Lambda_\beta}_{u,u;0}}

\def\qsxyn{K^{{\text{\tiny{\sc{RW}}}}_{\h}^d}_{s;x,y}}

\def\qzxn{K^{{\text{\tiny{\sc{RW}}}}_{\h}^d}_{0;x}}
\def\qtxn{K^{{\text{\tiny{\sc{RW}}}}_{\h}^d}_{t;x}}
\def\qsxn{K^{{\text{\tiny{\sc{RW}}}}_{\h}^d}_{s;x}}
\def\qtxyn{K^{{\text{\tiny{\sc{RW}}}}_{\h}^d}_{t;x,y}}
\def\qtxynint{K^{{\text{\tiny{\sc{RW}}}}_{\h}^d}_{t;[x]_{\h},[y]_{\h}}}
\def\qsxynint{K^{{\text{\tiny{\sc{RW}}}}_{\h}^d}_{s;[x]_{\h},[y]_{\h}}}

\def\qrox{K^{{\text{\tiny{\sc{RW}}}^d_{\delta}}}_{r_{1};x}}
\def\qrtx{K^{{\text{\tiny{\sc{RW}}}^d_{\delta}}}_{r_{2};x}}

\def\qroprtz{K^{{\text{\tiny{\sc{RW}}}^d_{\delta}}}_{r_{1}+r_{2};0}}

\def\qrox{K^{{\text{\tiny{\sc{RW}}}^d_{\delta}}}_{r_{1};x}}
\def\qrtx{K^{{\text{\tiny{\sc{RW}}}^d_{\delta}}}_{r_{2};x}}
\def\ut{\tilde{U}_{\beta}}
\def\utn{\tilde{U}_{\beta, n}}
\def\utnlf{\tilde{U}_{\beta,n,l}}
\def\utnk{\tilde{U}_{\beta,n_k}}
\def\utsyn{\tilde{U}_{\beta, n}^{y}(s)}

\def\uttxn{\tilde{U}_{\beta, n}^{x}(t)}
\def\uttxnk{\tilde{U}_{\beta,n_k}^{x}(t)}
\def\uttxnlf{\tilde{U}_{\beta,n,l}^{x}(t)}
\def\uttauxnlf{\tilde{U}_{\beta,n,l}^{x}(\tau)}
\def\xtaunlf{X_{\beta,n,l}^{\tau}}

\def\xttauoxnlf{X_{\beta,n,l}^{\tau_1,x_{1}}(t)}
\def\xttauoynlf{X_{\beta,n,l}^{\tau_1,y}(t)}

\def\xttautxnlf{X_{\beta,n,l}^{\tau_2,x_{2}}(t)}
\def\xttautynlf{X_{\beta,n,l}^{\tau_2,y}(t)}

\def\xttauxnlf{X_{\beta,n,l}^{\tau,x}(t)}
\def\xttauxonlf{X_{\beta,n,l}^{\tau,x_1}(t)}

\def\xttauxinlf{X_{\beta,n,l}^{\tau_{i},x_i}(t)}

\def\mtauxo{M^{\tau_1,x_{1}}(\cdot)}
\def\mtauxt{M^{\tau_2,x_{2}}(\cdot)}
\def\mtautx{M^{\tau,x}(t)}
\def\mtautxi{M^{\tau_i,x_{i}}(t)}
\def\mtautxj{M^{\tau_j,x_{j}}(t)}

\def\mtauuxj{M^{\tau_j,x_{j}}(u)}

\def\ntautxot{N^{\tau_{1,2},x_{1,2}}(t)}
\def\xztauxnlf{{X_{\beta,n,l}^{\tau,x}(0)}}
\def\xtautauxnlf{X_{\beta,n,l}^{\tau,x}(\tau)}
\def\xttxnlf{X_{\beta,n,l}^{t,x}(t)}
\def\utsynlf{\tilde{U}_{\beta,n,l}^{y}(s)}
\def\utrynlf{\tilde{U}_{\beta,n,l}^{y}(r)}
\def\xstauynlf{X_{\beta,n,l}^{\tau,y}(s)}
\def\xttauxrnlf{X_{\beta,n,l}^{\tau,x_r}(t)}
\def\xstauiynlf{X_{\beta,n,l}^{\tau_i,y}(s)}

\def\xstauxinlf{X_{\beta,n,l}^{\tau_i,x_{i}}(s)}
\def\uttx{\tilde{U}_{\beta}^{x}(t)}
\def\uttxD{\tilde{U}^x_{\beta,D}(t)}
\def\uttxR{\tilde{U}^x_{\beta,R}(t)}
\def\uttxnD{\tilde{U}^x_{\beta,n,D}(t)}
\def\uttxinD{\tilde{U}^{x_i}_{\beta,n,D}(t)}

\def\utsxjnD{\tilde{U}^{x_j}_{n,D}(s)}

\def\utuxinD{\tilde{U}^{x_i}_{n,D}(u)}
\def\uttxnR{\tilde{U}^x_{\beta,n,R}(t)}
\def\uttynR{\tilde{U}^y_{n,R}(t)}
\def\utrxnR{\tilde{U}^x_{n,R}(r)}
\def\uttxnlfR{\tilde{U}^x_{\beta,n,l,R}(t)}
\def\utrxnlfR{\tilde{U}^x_{\beta,n,l,R}(r)}
\def\uttynlfR{\tilde{U}^y_{\beta,n,l,R}(t)}

\def\uotn{{\tilde U}_n^{(1)}}
\def\utttn{{\tilde U}_n^{(2)}}
\def\yk{Y_{\beta,k}}
\def\ytxk{Y_{\beta,k}(t,x)}
\def\ytyk{Y_{\beta,k}(t,y)}
\def\yrxk{Y_{\beta,k}(r,x)}
\def\utx{U_{\beta}(t,x)}
\def\uty{U_{\beta}(t,y)}
\def\urx {U_{\beta}(r,x)}

\def\usy{U_{\beta}(s,y)}
\def\ytx{Y_{\beta}(t,x)}
\def\yty{Y_{\beta}(t,y)}
\def\yrx{Y_{\beta}(r,x)}
\def\uity{u_i^y(t)}
\def\ujty{u_j^y(t)}

\def\unx{u_0(x)}
\def\unxo{u_0(x_1)}
\def\unxt{u_0(x_2)}
\def\unxi{u_0(x_i)}

\def\uny{u_0(y)}
\def\un{u_0}
\def\zivtx{U_{\beta}^{(0)}(t,x)}
\def\zivsy{U_{\beta}^{(0)}(s,y)}
\def\nivtx{U_{\beta}^{(n)}(t,x)}
\def\nivsy{U_{\beta}^{(n)}(s,y)}
\def\npoivtx{U_{\beta}^{(n+1)}(t,x)}

\def\nmoivsy{U_{\beta}^{(n-1)}(s,y)}
\def\dnptx{D_{\beta,n,p}(t,x)}

\def\dnmopsy{D_{\beta,n-1,p}(s,y)}
\def\mnpt{D^*_{\beta,n,p}(t)}

\def\mnmops{D^*_{\beta,n-1,p}(s)}
\def\sW{\mathscr W}

\def\sTns{{\mathscr T}_{n,s}}
\def\sAn{{\mathscr A}_{n}}
\def\sAns{{\mathscr A}_{n}^{2}}

\def\wtxn{W_n^{x}(t)}
\def\wtyn{W_n^{y}(t)}
\def\wtxm{W_m^{x}(t)}
\def\wsxn{W_n^{x}(s)}
\def\wsyn{W_n^{y}(s)}
\def\DxBt{{\mathbb D}^x_B(t)}
\def\DxkBet{{\mathbb D}^{x,k}_{B,e}(t)}
\def\DxiBet{{\mathbb D}^{x,\infty}_{B,e}(t)}
\def\DxBsnkchnt{{\mathbb D}^{x}_{B,SC}(t)}
\def\P{{\mathbb P}}
\def\Pt{\tilde{\mathbb P}}
\def\EP{{\mathbb E}_{\P}}
\def\EPt{{\mathbb E}_{\Pt}}
\def\E{{\mathbb E}}
\def\N{{\mathbb N}}
\def\Ns{{\mathbb N}^{*}}
\def\Rd{{\mathbb R}^{d}}
\def\Q{\mathbb Q}
\def\TQ{{\T}_{\Q}}
\def\Zd{{\mathbb Z}^{d}}
\def\R{\mathbb R}
\def\Rs{{\mathbb R}^2}
\def\B{{\mathbb B}}
\def\S{\mathbb S}
\def\T{\mathbb T}
\def\X{\mathbb X}
\def\Rp{{\R}_+}

\def\Rpkop{(0,\infty)^k}
\def\Rpkmoop{(0,\infty)^{k-1}}

\def\Xnd{{\mathbb X}_n^{d}}
\def\Xmd{{\mathbb X}_m^{d}}
\def\Xd{{\mathbb X}^{d}}

\def\Xnldf{{\mathbb X}_{n,l}^{d}}

\def\sF{{\mathscr F}}

\def\sFt{{\mathscr F}_t}
\def\sFs{{\mathscr F}_s}
\def\sFr{{\mathscr F}_r}

\def\OFP{(\Omega,\sF,\P)}
\def\OFFtP{(\Omega,\sF,\{\sFt\},\P)}
\def\OFFtPt{(\tilde{\Omega},\tilde{\sF},\{\tilde{\sFt}\},\tilde{\P})}
\def\OFPt{(\tilde{\Omega},\tilde{\sF},\tilde{\P})}

\def\OtauotFtauotPtauott{(\tilde{\Omega}_{\tau_{1,2}},\tilde{\sF}_{\tau_{1,2}},\tilde{\P}_{\tau_{1,2}})}
\def\Skspace{(\Omega^S,\sF^S,\{\sFt^S\},\P^S)}
\def\locm{{\mathscr{M}}_2^{c,loc}}
\def\C{\mathrm C}
\def\H{\mathrm H}
\def\eqdef{:=}
\def\eqd{\overset{{\mathscr L}}{=}}

\def\sm{\setminus}

\def\lap{\Delta}

\def\df#1#2{\ds{\frac{#1}{#2}}}
\def\tf#1#2{\ts{\frac{#1}{#2}}}
\def\lbl#1{\label{#1}}
\def\intrd{\int_{\Rd}}
\def\intrdzt{\int_{\Rd}\int_0^t}
\def\intzt{\int_0^t}

\def\pa{\partial}
\def\pas{\partial^2}
\def\lqv{\left<}
\def\rqv{\right>}

\def\lab{\left|}
\def\rab{\right|}
\def\lpa{\left(}
\def\rpa{\right)}

\def\lbk{\left[}
\def\rbk{\right]}
\def\lbr{\left\{}
\def\rbr{\right\}}
\def\bdf{\begin{defn}}
\def\edf{\end{defn}}
\def\bcr{\begin{cor}}
\def\ecr{\end{cor}}
\def\bnt{\begin{notn}}
\def\ent{\end{notn}}
\def\brm{\begin{rem}}
\def\erm{\end{rem}}
\def\blm{\begin{lem}}
\def\elm{\end{lem}}
\def\bpf{\begin{pf}}
\def\bpfs{\begin{pf*}}
\def\epf{\end{pf}}
\def\epfs{\end{pf*}}
\def\eqrf{\eqref}
\def\beq{\begin{equation}}
\def\beqs{\begin{equation*}}
\def\eeq{\end{equation}}
\def\eeqs{\end{equation*}}
\def\bsp{\begin{split}}
\def\esp{\end{split}}
\def\bc{\begin{cases}}
\def\ec{\end{cases}}
\def\bt{\begin{tabular}}
\def\et{\end{tabular}}

\def\bthm{\begin{thm}}
\def\ethm{\end{thm}}
\def\bpr{\begin{prop}}
\def\epr{\end{prop}}
\def\bfr{\begin{framed}}
\def\efr{\end{framed}}
\def\bsh{\begin{shaded}}
\def\esh{\end{shaded}}
\def\bcm{\iffalse}
\def\ig{\iffalse}
\def\babs{\begin{abstract}}
\def\eabs{\end{abstract}}
\def\bit{\begin{itemize}}
\def\eit{\end{itemize}}
\def\ben{\begin{enumerate}}

\def\een{\end{enumerate}}
\def\babs{\begin{abstract}}
\def\eabs{\end{abstract}}
\def\Dn{\Delta_n}
\def\Dns{\Delta_n^2}

\def\btsie{e_{_{\mbox{\tiny BTBM}}}^{\mbox{\tiny SIE}}(a,u_0)}
\def\isltsie{e_{_{\mbox{\tiny $\beta$-ISLTBM}}}^{\mbox{\tiny SIE}}(a,u_0)}
\def\isltsieb{e_{_{\mbox{\tiny $\beta$-ISLTBM}}}^{\mbox{\tiny SIE}}(a,b,u_0)}

\def\isltrwsien{e_{_{\mbox{\tiny $\beta$-ISLTRW}}}^{\mbox{\tiny SIE}}(a,u_0,n)}
\def\isltrwsienl{e_{_{\mbox{\tiny $\beta$-ISLTRW}}}^{\mbox{\tiny t-SIE}}(a,u_0,n,l)}
\def\isltrwsienlaux{e_{_{\mbox{\tiny $\beta$-ISLTRW}}}^{\mbox{\tiny aux-SIE}}(a,u_0,n,l,\tau)}

\def\sAKtauot{{\mathscr{A}}_{\Upsilon}^{\tau_{1,2}}}
\def\sAKtauotim{{\mathscr{A}}_{\Upsilon,i_m}^{\tau_{1,2}}}

\def\utsyn{{\tilde U}_n^{y}(s)}

\def\p{\partial}
\def\paopd{\partial^{1+d}}
\def\Lb{L_{\beta}}
\def\ILb{\Lambda_{\beta}}

\def\nuset{\lbr2^{k};k\in\N\rbr}
\begin{document}
\title[$\beta$-Inverse-stable-L\'evy-time Brownian motion SIE{\scriptsize s} on $\Rp\times\Rd$]{Time-fractional and {memoryful} $\Delta^{2^{k}}$ SIE{\scriptsize s} on $\Rp\times\Rd$: 
how far can we push white noise?}
\author{Hassan Allouba}
\address{Department of Mathematical Sciences, Kent State University, Kent, Ohio 44242
\\email: allouba@math.kent.edu}
\thanks{}
\begin{abstract}
High order and fractional PDEs have become prominent in theory and in modeling many phenomena.  Here, we focus on the  regularizing effect of a large class of memoryful high-order or time-fractional PDEs---through their fundamental solutions---on stochastic integral equations (SIEs) driven by space-time white noise.  Surprisingly, we show that maximum spatial regularity is achieved in the fourth-order-bi-Laplacian case; and any further increase in the \emph{spatial}-Laplacian order is \emph{entirely} translated into additional \emph{temporal }regularization of the SIE.  We started this program in \cite{Abtbmsie,Abtpspde}, where we introduced two different stochastic versions of the fourth order memoryful  PDE associated with the Brownian-time Brownian motion (BTBM): (1) the BTBM SIE and (2) the BTBM SPDE, both driven by space-time white noise.  Under wide conditions, we showed the existence of random field locally-H\"older solutions to the BTBM SIE with striking and unprecedented  time-space H\"older exponents, in spatial dimensions $d=1,2,3$.   In particular, we proved that  the spatial regularity of such solutions is nearly locally Lipschitz in $d=1,2$.  This gave, for the first time,  an example of a  space-time white noise driven equation whose solutions are smoother than the corresponding Brownian sheet in either time or space.  

In this paper, we introduce the $2\beta^{-1}$-order  $\beta$-inverse-stable-L\'evy-time Brownian motion ($\beta$-ISLTBM) SIEs, $\beta\in\lbr1/2^{k};k\in\N\rbr$, driven by space-time white noise.  Based on the dramatic regularizing effect of the BTBM density ($\beta=1/2$), and since the kernels in these $\beta$-ISLTBM SIEs are fundamental solutions to higher order Laplacian PDEs; one may suspect that we get even more dramatic spatial regularity than the BTBM SIE case.   We show, however, that the BTBM SIE spatial regularity and its random field third spatial dimension limit are maximal among all $\beta$-ISLTBM SIEs---no matter how high we take the order $1/\beta$  of the Laplacian.  This gives a limit as to how far we can push the SIEs spatial regularity when driven by the rough white noise.  Furthermore, we show that increasing the order of the Laplacian $\beta^{-1}$ beyond the BTBM bi-Laplacian manifests entirely as increased temporal regularity of our random field solutions that asymptotically approaches that of the Brownian sheet as $\beta\searrow0$.  Our solutions are both direct and lattice limit solutions.  We treat many stochastic fractional PDEs and their corresponding higher order SPDEs, including BTBM and $\beta$-inverse-stable-L\'evy-time Brownian motion SPDEs, in separate articles.  
\end{abstract}
\maketitle
\newpage
\tableofcontents
\section{Introduction, motivation, and statement of results}\lbl{intro}
Lately, many phenomena in mathematical physics, fluids dynamics and turbulence models, mathematical finance, and the modern theory of stochastic processes have been related to and described through deterministic fractional and higher order evolution equations (e.g., see \cite{AN,Abtbs}, \cite{Aks}--\cite{AX}, \cite{BanDeB}, \cite{BMN}--\cite{Carr}, \cite{DeB}, \cite{Fun}--\cite{HO}, \cite{keyantuo-lizama,KuLag}, \cite{MeerSik}--\cite{BOap09}, and \cite{T}); and it is only natural to investigate these important equations under the influence of a driving random noise.

In the two recent articles \cite{Abtbmsie,Abtpspde} we introduced two new stochastic versions of fourth order memory-preserving (which we coin \emph{memoryful}) deterministic PDEs related to Brownian-time processes (BTPs)\footnote{A
BTP, in its simplest form,  is a process $X^x\lpa\lab B_t\rab\rpa$ in which $X^x$ is a Markov process starting at $x\in\Rd$ and
$B$ is an independent one dimensional BM starting at $0$.
A Brownian-time Brownian motion (BTBM) is a BTP in which
$X^x$ is also a Brownian motion.
BTPs include many new and quite interesting processes (see \cite{Abtp1,Abtp2,DeB,Nanesd}), which we are currently investigating in several directions (e.g., \cite{AN,AD,AL,AX}).
With the exception of the Markov snake of Le Gall (\cite{LeG}), BTPs fall outside the classical theory of Markov, Gaussian, or semimartingale processes.  We label BTP PDEs as memoryful since the initial data is part of the PDE itself (see \eqref{btppdedet})}---introduced in \cite{Abtp1,Abtp2}---driven by space-time white noise: 
\ben
\item  the space-time-white-noise-driven Brownian-time Brownian motion (BTBM) SPDE
\beq\lbl{btpspde}
\bc
\pa_{t}U= \df{\Delta\un }{\sqrt{8\pi t}}+\df18\Delta^2 U+
a(U)\paopd_{t,x} W,&(t,x)\in(0,\infty)\times\Rd; \cr
U(0,x) = \unx,&x\in\Rd,
\ec
\eeq
where $\paopd_{t,x} W$ is the space-time white noise on $\Rp\times\Rd$---and on a probability space $\OFP$---that corresponds to the Brownian sheet $W$; and 
\item the stochastic integral equation we called BTBM SIE
\beq\lbl{btpsol}
\bsp
U(t,x)=\intrd\KBtxy \uny dy+ \intrdzt\KBtsxy a(U(s,y))\sW(ds\times dy)
\end{split}
\eeq
where $\KBtxy$ is the density of a $d$-dimensional Brownian-time Brownian motion given by:
\beq\lbl{btpden}\KBtxy=2\int_0^\infty \psxy\ptzs ds
\eeq
with  $\psxy=\tf{e^{-|x-y|^2/2s}}{(2\pi s)^{{d}/{2}}}$ and $\ptzs=\tf{e^{-s^2/2t}}{\sqrt{2\pi t}}$; and  where $\sW$ is the white noise on $\Rp\times\Rd$.  
\een
Unlike the deterministic case $a\equiv0$, \eqref{btpspde} and \eqref{btpsol} behave differently, and each is quite interesting in its own right.
Each of these two equations gives a different stochastic interpretation of the memoryful BTBM PDE in \cite{Abtp1,Abtp2}:
\beq\lbl{btppdedet}
\bc
\pa_{t} u= \ds{\df{\Delta\un }{\sqrt{8\pi t}}+\df18\Delta^2 u};&(t,x)\in(0,\infty)\times\Rd, \cr
u(0,x) = \unx;&x\in\Rd.
\ec
\eeq
and its equivalent integral form\footnote{For a review of the BTPs higher order and fractional PDEs connections and generalizations, as well as connection to the important Kuramoto-Sivashinsky PDE, we refer  the reader to \cite{Abtp1,Abtp2,Aks,MNV09,Nanesd,BOap09} and the references therein.  The connection of BTPs to their fourth order PDEs (including \eqref{btppdedet}) was first given in \cite{Abtp1}.  Also, their connection to time-fractional PDEs was first established implicitly via the half derivative generator in \cite{Abtp1}.  In \cite{MNV09,Nanesd,BOap09} the equivalence between a large class of high order and time-fractional PDEs, including \eqref{btppdedet} and \begin{equation}\label{BTBM-half}
\begin{cases}
\displaystyle\p_{t}^{\frac12} u=\frac{1}{\sqrt{8}}\lap u;& t\in(0,\infty),x\in\Rd,\\
u(0,x)=\unx;&x\in\Rd,
 \end{cases}
\end{equation} was established explicitly, using the Caputo fractional derivative.  
For a discussion of interesting aspects of these PDEs see also the introduction in \cite{Abtbmsie}.  In the new multiparameter-time case the reader is referred to  \cite{Abtbs,AN}.}
\beq\lbl{btpsoldet}
\bsp
u(t,x)=\intrd\KBtxy \uny dy.
\end{split}
\eeq
As proven in \cite{Abtpspde,Abtbmsie}, the SIE \eqref{btpsol}---which we also denote by $\btsie$---has real random field solutions in $d=1,2,3$ with striking H\"older regularity in which the time-space H\"older exponents are $\lpa{{\tf{4-d}{8}}^{-},{\lpa\tf{4-d}{2}\wedge 1\rpa}^{-}}\rpa$, as we recall precisely in \secref{btbmsieresrec} below\footnote{In particular, as was established in \cite{Aksspde}, the BTBM SIE \eqref{btpsol} has nearly locally Lipschitz solutions  in $d=1,2$.  This fact provided for the first time a counterexample to the common folklore {\it non}-wisdom that ``a solution to a space-time-white-noise-driven equation cannot have a solution that is more regular, temporally or spatially, than the Brownian-sheet in the underlying white noise''}, and it is similar in regularity to the following L-Kuramoto-Sivashinsky (L-KS)\footnote{The L in the name refers to the linearized PDE part.  Such L-KS SPDE is treated in \cite{Aksspde}.} SPDE 
  \begin{equation} \label{lksdet}
 \begin{cases} \p_{t} U=\ds
-\frac18\lap^2U-\frac12\lap U-\frac12U+a(U)\pa^{1+d}_{t,x} W, & (t,x)\in(0,\infty )\times\Rd;
\cr U(0,x)=\unx, & x\in\Rd,
\end{cases}
\end{equation}
obtained from the linearized KS PDE in \cite{Aks} by adding a multiplicative space-time white noise term (see \cite{Aksspde}).  In \cite{AMfrac,AMfracwav}, we treat a large class of higher order and fractional---and rougher---SPDEs, including \eqref{btpspde} and its equivalent time-fractional SPDE
\begin{equation}\label{btbm-half-spde}
\begin{cases}
\displaystyle\p_{t}^{\frac12} U=\frac{1}{\sqrt{8}}\lap U+
a(U) \pa^{1+d}_{t,x} W;& t\in(0,\infty),x\in\Rd,\\
U(0,x)=\unx;&x\in\Rd,
 \end{cases}
\end{equation} 
where $\p_{t}^{\frac12}$ is a factional derivative in time (see e.g.~\cite{MeerSik}).  

In this article, we focus on a large class of fascinating stochastic integral equations driven by space-time white noise and generalizing the BTBM SIE \eqref{btpsol}: the $\beta$-inverse-stable-L\'evy-time Brownian motion SIEs ($\beta$-ISLTBM SIEs), which we discuss in more details  in \secref{frstres} below.  These SIEs are obtained from the BTBM SIE in \eqref{btpsol} by replacing the BTBM density with the fundamental solution to the $2\beta^{-1}=2\nu$ order, $\beta^{-1}\in\lbr2^k;k\in\N\rbr$, memoryful PDEs
\beq\lbl{2nupde}\bc
\ds\p_{t}{u_{\beta}(t,x)}=\sum_{\kappa=1}^{\nu-1}\frac{\Delta^{\kappa}\unx}{2^{\kappa}t^{1-\kappa/\nu}}{E_{\beta,\kappa}}+\frac{\Delta^{\nu}u_{\beta}(t,x)}{2\nu};
& (t,x)\in(0,\infty)\times\Rd\\
u_{\beta}(0,x)=u_0(x),&x\in\Rd
\ec
\eeq
and their equivalent time-fractional PDEs
\beq\lbl{fracbetapde}\bc
\ds\p^{\beta}_{t}{u_{\beta}(t,x)}=\frac12\Delta u_{\beta}(t,x)
& (t,x)\in(0,\infty)\times\Rd\\
u_{\beta}(0,x)=u_0(x),&x\in\Rd,
\ec
\eeq
where\footnote{As usual, $\E$ denotes the expectation operator.} $E_{\beta,\kappa}=\frac{\E\lpa\Lambda_{\beta}(1)\rpa^{\kappa}}{\kappa!}$,   the process $\ILb$ is the $\beta$-inverse-stable-L\'evy motion described in \secref{frstres}  below, and $\p_{t}^{\beta}$ is the well known Caputo fractional derivative of order $\beta\in\lbr1/2^k;k\in\N\rbr$ in time (see e.g.~\cite{MeerSik}).    

Based on the dramatic regularizing effect of the BTBM density on the space-time white noise driven BTBM SIE \eqref{btpsol} as just described above (see also \thmref{lip} below), and due to the fact that the kernels in the $\beta$-ISLTBM SIEs of this article are fundamental solutions to the higher order PDEs \eqref{2nupde}; one may suspect that we get even more dramatic spatial regularity than the BTBM SIE case, possibly obtaining random field solutions in arbitrarily high spatial dimensions as $\beta\searrow0$ ($\nu\nearrow\infty$) instead of just $d=3$ as in the BTBM case ($\beta=1/\nu=1/2$).   We show, however, that the BTBM SIE spatial regularity and its random field third spatial dimension limit  are maximal among all $\beta$-ISLTBM SIEs; no matter how small we take $\beta$ (how high we take the order $\beta^{-1}$ of the Laplacian).  Further, we show that increasing the order $\beta^{-1}$ of the \emph{spatial} Laplacian beyond the BTBM order of $2$ translates \emph{entirely} into temporal regularization of our $\beta$-ISLTBM SIEs\footnote{I.e., the extra regularizing ``energy'' of spatial Laplacians of orders higher than that of the bi-Laplacian is converted to extra temporal regularity, when faced with the extremely rough driving space-time white noise.}.  This surprising result is the regularity content of our two main theorems: \thmref{mainthm1} and \thmref{mainthm2} below.
\subsection{Recalling the Brownian-time Brownian motion SIE case}\lbl{btbmsieresrec}
Before stating our first main result, it is instructive to recall the BTBM SIE results in \cite{Abtbmsie}\footnote{Earlier, in \cite{Abtpspde}, the additive noise case $a\equiv1$ for $\btsie$ was considered; and the existence of a pathwise unique continuous BTBM SIE solution $U(t,x)$ for $x\in\Rd$ and $d=1,2,3$, such that
$$\sup_{x\in\Rd}\EP|U(t,x)|^{2p}\le C\lbk1+t^{\tf{(4-d)p}{4}}\rbk;\quad t>0,\ d=1,2,3,\  p\ge1,$$
was proved.}.
Following \cite{Abtbmsie}, we denote by\footnote{Throughout the paper, $\T=[0,T]$ for some fixed but arbitrary $T>0$.  Here and in the sequel $\mathrm{C}_{b}^{p,\gamma}(\Rd;\R)\subset\mathrm{C}_{b}^p(\Rd;\R)$ denotes the space of bounded $p$-times continuously differentiable functions such that all derivatives up to (and including) the $p$-th order are bounded and all $p$-th order derivatives are H\"older continuous, with some H\"older exponent $0<\gamma\le1$.   Also, the boundedness conditions on $\un$ and its derivatives may easily be relaxed as in \cite{AN}.} $\H^{\gamma_t^{-},\gamma_s^{-}}(\T\times\Rd;\R)$ the space of real-valued locally H\"older functions on $\T\times\Rd$ whose time and space H\"older exponents are in $(0,\gamma_t)$ and $(0,\gamma_s)$, respectively.  The first main result in \cite{Abtbmsie} is now restated.
\bfr
\bthm[Allouba \cite{Abtbmsie}]\lbl{lip}  
Fix $0<\gamma\le1$.  Assume the following Lipschitz and growth conditions 
\beq\lbl{lc}
\bc
(a) \hspace{-2.5mm}&\lab a(u)-a(v)\rab\le C\lab u-v\rab \mbox \quad u,v\in\R;\\
(b) \hspace{-2.5mm}&a^2(u)\le C (1+u^2);\quad u\in\R,\\
(c) \hspace{-2.5mm} &u_0\in \mathrm{C}_{b}^{2,\gamma}(\Rd;\R)\mbox{ and nonrandom } ,\ \forall\ d=1,2,3.
\ec
\tag{\mbox{Lip}}
\eeq
hold.  Then there exists a pathwise-unique strong solution $(U,\sW)$ to $\btsie$
on $\Rp\times\Rd$, for $d=1,2,3$, which is $L^p(\Omega)$-bounded 
on $\T\times\Rd$ for all $p\ge2$.  Furthermore, $U\in\H^{{\tf{4-d}{8}}^{-},{\lpa\tf{4-d}{2}\wedge 1\rpa}^{-}}\lpa\T\times\Rd;\R\rpa$ for every $d=1,2,3$.  
\ethm
\efr
\thmref{lip} states that the stochastic kernel integral equation \eqref{btpsol} has ultra regular strong\footnote{Here strong is in the stochastic sense of the noise $\sW$ and its probability space $\OFFtP$ being fixed a priori.  Throughout this article, whenever needed, we will assume that our filtrations satisfy the usual conditions without explicitly stating so.} solutions on $\Rp\times\Rd$, namely $U\in\H^{{\tf{3}{8}}^{-},1^{-}}(\T\times\R;\R)$, $U\in\H^{{\tf{1}{4}}^{-},1^{-}}(\T\times\R^{2};\R)$, and $U\in\H^{{\tf{1}{8}}^{-},{\tf12}^{-}}(\T\times\R^{3};\R)$.   I.e., in space, the BTBM paths have a rather remarkable---and initially-surprising---nearly local Lipschitz regularity for $d=1,2$; and nearly local H\"older $1/2$ regularity in $d=3$.  This is remarkable because the BTBM kernel is able, in $d=1,2$, to spatially regularize such solutions beyond the traditional H\"older-$1/2^{-}$ spatial regularity of the underlying Brownian sheet corresponding to the driving space-time white noise\footnote{As noted in \cite{Abtbmsie}, it is important to note here that the common ``folklore wisdom'' of solutions of space-time-white-noise driven equations not being smoother than the associated Brownian sheet---in either space or time---originated from  the predominant case of SPDEs, in which either the underlying kernel is that of a Brownian motion or the spatial operator is a Laplacian.  The kernel $\KBtxy$, however, is much more regularizing to the space-time-white-noise driven $\btsie$ than the density of BM is to its corresponding equation. This becomes evidently clear in \lemref{3rdQinequality}, \lemref{4thQinequality}, and \lemref{2ndQinequality} (compare to the more traditional BM and random walk case in \cite{Asdde2}).}.  This degree of smoothness is unprecedented for space-time white noise driven kernel equations or their corresponding SPDEs; and the BTBM SIE is thus the first such example.  In time, our solutions are locally $\gamma$-H\"older continuous with dimension-dependent exponent $\gamma\in\lpa0,\tf{4-d}{8}\rpa$ for $d=1,2,3$.  This is in sharp contrast to traditional second order reaction-diffusion (RD) and other heat-operator-based SPDEs driven by space-time white noise, whose fundamental kernel is the Brownian motion density and whose real-valued mild solutions are confined to the case $d=1$.  In this regard, the dichotomy between the rougher paths of BTBMs as compared to standard Brownian motions on the one hand (quartic vs.~quadratic variations) and the stronger regularizing properties of the BTBM density vs.~the BM one on the other hand is certainly another interesting point to make\footnote{We observe in passing that---roughly speaking---the paths of $\btsie$ in $d=1$ are effectively $3/2$ times as smooth as the RD SPDE paths in $d=1$, in $d=2$ the BTBM SIE is as smooth as an RD SPDE in $d=1$, and in $d=3$ our BTBM SIE is half as smooth as an RD SPDE in $d=1$.  Also, for $d=2,3$, the spatial regularity is roughly four times the temporal one, and in $d=1$ the spatial regularity is maximized at a near Lipschitz vs near H\"older $3/8$ in time.}.   

\subsection{The $\beta$-inverse-stable-L\'evy-time Brownian motion SIE: the first main theorem}\lbl{frstres}
In the first main result of this article, we generalize the first BTBM SIE result in \cite{Abtbmsie} \thmref{lip} to the interesting case of the inverse-stable-L\'evy-time Brownian motion SIE  with index $\beta=1/\nu$, $\nu\in\lbr2^k;k\in\N\rbr$ ($\beta$-ISLTBM SIE), which we now motivate and introduce\footnote{Throughout this article we assume that $\nu=\beta^{-1}\in\lbr2^k;k\in\N\rbr$, where $\N$ is the set of natural numbers.  The case $\beta^{-1}=2$ is the BTBM SIE case, with a minor scaling of the Brownian motion as discussed in \cite{AN}.}.  This generalization allows us to better appreciate how hard it is to smooth away space-time white noise.
\subsubsection{Recalling $\beta$-ISLTBM}
Inverse stable subordinator---which we also call $\beta$-inverse-stable-L\'evy motion and denote by $\Lambda_{\beta}$---arise in the work of Meerschaert et al.~\cite{limitCTRW,Zsolution} as scaling limits of continuous time random walks.  Let $S(n)=Y_1+\cdots+Y_n$ a sum of independent and identically distributed random variables with $EY_n=0$ and $EY_n^2<\infty$.  The scaling limit $c^{-1/2}S([ct])\Rightarrow B(t)$ as $c\to\infty$ is a Brownian motion $B$ at time $t$, which is normal with mean zero and variance proportional to $t$.  Consider $Y_n$ to be the random jumps of a particle.  If we impose a random waiting time $T_n$ before the $n$th jump $Y_n$, then the position of the particle at time $T_n=J_1+\cdots+J_n$ is given by $S(n)$.  The number of jumps by time $t>0$ is $N(t)=\max\{n:T_n\leq t\}$, so the position of the particle at time $t>0$ is $S(N(t))$, a subordinated process.  If $\P(J_n>t)=t^{-\beta}l(t)$ for some $0<\beta<1$, where $l(t)$ is slowly varying, then the scaling limit $c^{-1/\beta}T_{[ct]}\Rightarrow \Lb(t)$ is a strictly increasing stable L\'evy motion $L_{\beta}$ at time $t$ and with index $\beta$, sometimes called a stable subordinator.  The jump times $T_n$ and the number of jumps $N(t)$ are inverses $\{N(t)\geq x\}=\{T(\ceil{x})\leq t\}$ where $\ceil{x}$ is the smallest integer greater than or equal to $x$.  It follows that the scaling limits are also inverses $c^{-\beta}N(ct)\Rightarrow \ILb(t)$ where $\ILb(t)=\inf\{x:L(x)> t\}$, so that $\{\ILb(t)\leq x\}=\{\Lb(x)\geq t\}$.  We call the process $\ILb$ a $\beta$-inverse-stable-L\'evy motion.  Since $N({ct})\approx c^{\beta}\Lambda(t)$, the particle location may, for large $c$, be approximated by  $c^{-\beta/2}S(N({[ct]}))\approx (c^\beta)^{-1/2}S(c^\beta \ILb(t))\approx B(\ILb(t))$, a Brownian motion subordinated to the inverse or hitting time (or first passage time) process of the stable subordinator $\Lb$.  The random variable $\Lb(t)$ has a smooth density.  For properly scaled waiting times, the density of $\Lb(t)$ has Laplace transform $e^{-ts^\beta}$ for any $t>0$, and the random variables $\Lb(t)$ and $t^{1/\beta}\Lb(1)$ are identically distributed.  Writing $g_\beta(u)$ for the density of $\Lb(1)$, it follows that $\Lb(t)$ has density $t^{-1/\beta}g_\beta(t^{-1/\beta}u)$ for any $t>0$.  Using the inverse relation $\P(\ILb(t)\leq x)=\P(\Lb(x)\geq t)$ and taking derivatives, it follows that $\ILb(t)$ has density
\begin{equation}\label{Etdens}
\KItzx=t\beta^{-1}x^{-1-1/\beta}g_\beta(tx^{-1/\beta}),
\end{equation}

As noted above, we assume throughout this article that $\nu=\beta^{-1}\in\lbr2^k;k\in\N\rbr$.  In this case, there is a simple connections between $k$-iterated Brownian-time Brownian motion and $\beta$-ISLTBM.   We denote by $$\mathbb{B}^{x}_{\underset{i=1}{\overset{k}{\bigcirc}} B_{i}}(t):=B^{x}\lpa\lab B_{k}\lpa\cdots B_{2}\lpa\lab B_{1}(t)\rab\rpa\cdots\rpa\rab\rpa$$ a $k$-iterated Brownian-time Brownian motion at time $t$; where $\lbr B_{i}\rbr_{i=1}^{k}$ are independent copies of a one dimensional scaled Brownian motion starting at zero, with density $\frac{1}{\sqrt{4\pi t}}\exp\left(-\frac{z^2}{4t}\right)=\lpa1/\sqrt2\rpa K^{\text{\tiny{\sc{BM}}}}_{t;0,z/\sqrt2}$, and independent from the standard $d$-dimensional Brownian motion $B^{x}$, which starts at $x\in\Rd$.  By $\mathbb{B}^{x}_{\Lambda_{1/2^{k}}}(t)=B^{x}\lpa\Lambda_{1/2^{k}}(t)\rpa$ we mean a $d$-dimensional $\beta$-ISLTBM---with $\beta=1/2^{k}$---starting at $x\in\Rd$ and evaluated at time $t$; in which the outer BM $B^{x}$ and the inner $\Lambda_{1/2^{k}}$ are independent.  

\blm[The $\beta$-ISLTBM density]\lbl{conkitbtbmisltbm} The probability distributions of $\mathbb{B}^{x}_{\underset{i=1}{\overset{k}{\bigcirc}} B_{i}}(t)$ and $\mathbb{B}^{x}_{\Lambda_{1/2^{k}}}(t)$ are the same for every $k=1,2,\ldots$ and every $t\ge0$.  In particular, when $\beta=1/2^{k}$, $k\in\N$, the $\Lambda_{\beta}$ and the  $\beta$-ISLTBM transition densities are given by
\beq\lbl{denofkitbtbm}
\bsp
\KItzso&=2^{\frac{k}{2}}\int_{\Rpkmoop} \ptzsk\prod_{i=0}^{k-2}\ptzskiimo ds_{2}\cdots ds_{k}
\\ \KIBtx&=2^{\frac{k}{2}}\int_{\Rpkop} \psoxz \ptzsk\prod_{i=0}^{k-2}\ptzskiimo ds_{1}\cdots ds_{k},
\end{split}
\eeq
respectively\footnote{We are using the convention $\prod_{i=0}^{-1}c_{i}=1$ for any $c_{i}$ and the convention $\int_{\Rp^{0}}f(s)ds=f(s)$, for every $f$.  Also, we use the convention that the case $k=0$ ($\beta=1$) in the $\beta$-ISLTBM is the standard $d$-dimensional Brownian motion case.}.
\elm
\bpf
Let  $\beta=1/2^{k}$, $k\in\N$.  By Corollary 3.1 in \cite{Nanesd} we get that the distributeions are the same. Now, equation (0.14) in \cite{AN} gives us that
\begin{equation}\begin{split}\label{Ydens}
\KIhalftzx=\frac{2}{\sqrt{4\pi t}}\exp\left(-\frac{x^2}{4t}\right)
=\frac{2}{\sqrt{2}}K^{\text{\tiny{\sc{BM}}}}_{t;0,\frac{x}{\sqrt2}}.
\end{split}\end{equation}
This, together with Lemma 3.1 and Lemma 3.2 in \cite{Nanesd} and a simple conditioning argument using the independence of all the Brownian motions, we  immediately obtain \eqref{denofkitbtbm} as asserted.
\epf
We now define our $\beta$-ISLTBM SIE as the stochastic integral equation:
\beq\lbl{isltbmsie}
\bsp
U_{\beta}(t,x)=\intrd\KIBtxy \uny dy+ \intrdzt\KIBtsxy a(\usy)\sW(ds\times dy)
\end{split}
\eeq
where $\KIBtxy$ is the transition density of a $d$-dimensional $\beta$-ISLTBM, starting from $x\in\Rd$, ${\mathbb B}^{x}_{\ILb}:=\lbr B^{x}(\ILb(t)),t\ge0\rbr$ given by\footnote{Compare with  the  expression of $\KIBtxy$ \lemref{conkitbtbmisltbm} in terms of scaled BM transition densities.}:
\beq\lbl{btpden}\KIBtxy=\int_0^\infty \psxy\KItzs ds.
\eeq
We also denote the $\beta$-ISLTBM SIE \eqrf{isltbmsie} by $\isltsie$.  Just as in the BTBM SIE case, $\isltsie$ is one of two {\em different} stochastic versions\footnote{The other stochastic version is the $2\nu$ or the time-fractional $\beta$ order SPDE obtained from \eqref{2nupde} or from \eqref{fracbetapde} by adding the white noise term as in \cite{AMfrac}.} of the higher order ($2\nu=2\beta^{-1}$) memoryful PDEs \eqref{2nupde} and their equivalent time fractional PDEs \eqref{fracbetapde}.  

Of course, in the deterministic case, both \eqref{2nupde} and \eqref{fracbetapde} are equivalent to their integral form 
\beq\lbl{betainteg}
u_{\beta}(t,x)=\int_{\Rd}\KIBtxy dy.
\eeq
\subsubsection{First theorem: $2\beta^{-1}$ order SIEs regularity and third dimension maximality}
Our first main theorem is now stated.
\bfr
\bthm[Spatio-temporal regularity and third dimension maximality: direct solution]\lbl{mainthm1}  Fix $\beta=1/\nu$, $\nu\in\lbr2^k;k\in\N\rbr$.
Assume the following Lipschitz, growth, and initial smoothness  conditions 
\beq\lbl{lcnd}
\bc
(a) \hspace{-2.5mm}&\lab a(u)-a(v)\rab\le C\lab u-v\rab \mbox \quad u,v\in\R;\\
(b) \hspace{-2.5mm}&a^2(u)\le C (1+u^2);\quad u\in\R,\\
(c) \hspace{-2.5mm} &u_0\in \mathrm{C}_{b}^{2\nu-2,\gamma}(\Rd;\R)\mbox{ and nonrandom } ,\ \forall\ d=1,2,3.
\ec
\tag{\mbox{Lip}}
\eeq
hold.  Then there exists a pathwise-unique strong solution $(U_{\beta},\sW)$ to $\isltsie$
on $\Rp\times\Rd$, for {$d=1,2,3$}, which is $L^p(\Omega)$-bounded 
on $\T\times\Rd$ for all $p\ge2$.  Furthermore, $U_{\beta}\in\H^{\lpa{\tf{2\nu-d}{4{\nu}}}\rpa^{-},{\lpa\tf{4-d}{2}\wedge 1\rpa}^{-}}(\T\times\Rd;\R)$ for every $d=1,2,3$.
\ethm
\efr

\thmref{mainthm1} states that, for $\beta=1/\nu$ and $\nu\in\lbr2^k;k\in\N\rbr$,  these $2\beta^{-1}$ order $\beta$-ISLTBM SIEs have quite interesting locally-H\"older solutions with temporal and spatial H\"older exponents given by $\lpa{\tf{2\nu-d}{4\nu}}\rpa^{-}$ and ${\lpa\tf{4-d}{2}\wedge 1\rpa}^{-}$, respectively, for $d=1,2,3$.   Comparing this regularity with the corresponding result for the fourth order BTBM SIE in \thmref{lip}, we see that the spatial regularity---spatial H\"older exponent and the maximum spatial dimension of 3---is identical.  Since, the fundamental density (fundamental solution) estimates leading to the regularity conclusions of \thmref{mainthm1}---\lemref{2ndQinequality} to \lemref{3rdQinequality}---are sharp\footnote{We will have more to say about the regularity of these $\beta$-ISLTBM SIEs in \cite{AX}.  We also briefly note that by third dimension maximality, we mean maximality among \emph{integer} dimensions.}, this means that there is a limit as to how far we can push against the powerful roughening effect of the driving space-time white noise.  Despite the fact that these SIEs are co-driven by fundamental solutions of arbitrarily high order ($2\beta^{-1}$) PDEs involving the spatial $\beta^{-1}$-Laplacian operators, we can obtain locally H\"older real random field solutions only up to three spatial dimensions and with spatial H\"older exponents up to the maximal BTBM bi-Laplacian case ($\beta^{-1}=2$), for all $\nu=\beta^{-1}\in\lbr2^k;k\in\N\rbr$, no matter how large  $\beta^{-1}$ is.  

To appreciate the richness of the  regularizing effect of these $\beta$-ISLTBM SIEs, however,  we need to look beyond just the spatial dimensionality and regularity aspects.   So, we will now examine the conclusion of \thmref{mainthm1} regarding the maximum temporal (effective) H\"older exponent\footnote{The effective H\"older exponent is the minimum of the spatial and temporal H\"older exponents, which of course determine how smooth the random field solutions are as functions of both time and space together.}, as $\beta\searrow0$.  As observed above, the strong roughening influence of the space-time white noise prevents further spatial smoothing of our $\beta$-ISLTBM SIEs beyond the BTBM bi-Laplacian case, no matter how large $\beta^{-1}$ gets.  However, all of the extra smoothing ``energy'' resulting from increasing the \emph{spatial} Laplacian order $\beta^{-1}$ cannot simply be ``destroyed'' by the white noise; and it is converted instead into \emph{temporal} regularization of these $\beta$-ISLTBM SIEs (as $\beta\searrow0$).  \thmref{mainthm1} describes precisely this temporal effect in terms of H\"older exponents.  In particular, the maximum effective regularity of the $\beta$-ISLTBM SIEs increases asymptotically to the well-known  H\"older $(1/2)^{-}$ regularity of the Brownian sheet; i.e., the maximum effective H\"older exponent  $\lpa{\tf{2\beta^{-1}-d}{4\beta^{-1}}}\rpa^{-}\nearrow\frac12^{-}$ as $\beta\searrow0$ for every $d=1,2,3$.  The following table summarizes our regularity findings and compares them to the more standard and classical case of reaction-diffusion SPDEs driven by space-time white noise.

 \begin{table}[h]
\bt[t]{|c|c|c|c|c|}\hline
$d$ & \multicolumn{2}{c|}{Random Field Solutions}
&\multicolumn{2}{c|}{H\"older Exponent (time, space)}
\\
 \cline{2-5}
 & RD SPDE & $\beta$-ISLTBM SIE& RD SPDE & $\beta$-ISLTBM SIE\\
 \hline
$1$ & Yes & Yes & $\lpa\tf14^{-}, \tf12^{-}\rpa$ &$\lpa\lpa{\tf{2\nu-1}{4\nu}}\rpa^{-},{1^{-}}\rpa$\\
$2$ & No & Yes & N/A &  $\lpa\lpa{\tf{2\nu-2}{4{\nu}}}\rpa^{-},{1}^{-}\rpa$ \\
$3$ & No & Yes & N/A & $\lpa\lpa{\tf{2\nu-3}{4{\nu}}}\rpa^{-},{\lpa\tf{1}{2}\rpa}^{-}\rpa$  \\
\hline
\et
\smallskip\\
\caption{$\beta$-ISLTBM SIEs ($\nu=\beta^{-1}\in\lbr2^k;k\in\N\rbr$) vs RD SPDEs ($\beta=1$).}
\lbl{tbl}
\end{table}

To prepare for the statement of our results under the less-than-Lipschitz conditions in \eqref{cnd} (\thmref{mainthm2} below), we now introduce the spatial lattice version of $\isltsie$ as well as introduce the new associated process we call $\beta$-inverse-stable-L\'evy-time random walk and define the lattice limit solutions involved in the statement of \thmref{mainthm2}.   The main machinery we use in the proof in this case is our K-martingale approach, which we introduced and used in the BTBM SIE case in \cite{Abtbmsie}.  We recall this approach, adapting it to our setting\footnote{All we need to adapt it here is to replace the BTRW kernel of \cite{Abtbmsie} with the $\beta$-ISLTRW one in \eqref{btrwdsty} below. }, in Section \ref{Kmartsec}. 
\subsection{The spatial lattice version and the second main result}  As in \cite{Abtbmsie}, we now spatially discretize $\isltsie$.   This accomplishes at least two things: (1) it gives a multiscale view of $\isltsie$ and (2) it allows us to prove our existence and regularity results without the Lipschitz condition on $a$. 
\subsubsection{$\beta$-inverse-stable-L\'evy-time random walk on the lattice}\lbl{btc}
In \cite{Asdde1,Asdde2}, standard continuous-time random walks
on a sequence of refining spatial lattices $$\lbr\Xnd:=\prod^{d }_{i=1 }\lbr\ldots, -2\h,-\h,0,\h, 2\h, \ldots\rbr=\h\Zd\rbr_{n\ge1}$$
(with the step size $\h\searrow0$ as $n\nearrow\infty$) played a crucial role---through their densities---in
obtaining our results for second order
RD SPDEs.  In \cite{Abtbmsie}, in the fourth order Brownian-time setting, that role is played by
Brownian-time random walks on $\Xnd$:
\beq\lbl{btrwdef}
\S^x_{B,\h}(t)\eqdef S^x_{\h}\lpa\lab B_t\rab\rpa;\quad 0\le t<\infty, x\in\Xnd
\eeq
where $S_{\h}^x(t)$ is a standard $d$-dimensional continuous-time symmetric RW
starting from $x\in\Xnd$ and $B$ is an independent one-dimensional
BM starting at $0$.  The subscript $\delta_n$ in \eqref{btrwdef}
is to remind us that the lattice step size is $\delta_n$ in each of the $d$ directions.

In this article, we replace Brownian-time random walk with $\beta$-inverse-stable-L\'evy-time random walk ($\beta$-ISLTRW):
\beq\lbl{isltrwdef}
\S^x_{\ILb,\h}(t)\eqdef S^x_{\h}\lpa\ILb(t)\rpa;\quad 0<\beta<1,\ 0\le t<\infty,\ x\in\Xnd
\eeq
 It is then clear that the transition probability (density) $\Ktxyn$ of the $\beta$-ISLTRW $\S_{\ILb,\h}^x(t)$ on $\Xnd$ is given by\footnote{Throughout this article, $\Ktxn\eqdef\Ktxzn$ (with a similar convention for all transition densities).}
\beq\lbl{btrwdsty}
\Ktxyn=2\int_0^\infty\qsxyn\KItzs ds; \quad 0<\beta<1,\ 0< t<\infty, \ x,y\in\Xnd
\eeq
where $\qtxyn$ is the continuous-time random walk transition density
starting at $x\in\Xnd$ and going to $y\in\Xnd$ in time $t$,
in which the times between transitions are exponentially distributed
with mean $\h^{2d}$.    I.e., $\qtxn$ is the fundamental solution to the deterministic
heat equation on the lattice $\Xnd:$
\beq\lbl{latheat}
\displaystyle \frac{du_{n}^x(t)}{dt}
 = \frac{1}{2}\Delta_nu_{n}^x(t);\hspace{2mm}(t,x)\in(0,\infty)\times\Xnd
\eeq
where $\mathscr{A}_n\eqdef \Delta_n/2$ is the generator of the RW $S^x_{\h}(t)$ on $\Xnd$.

By mimicking our proof of Theorem 0.3 in \cite{AN}, we easily get a $2\nu$  order differential-difference equation connection to $\beta$-ISLTRW:
\blm[$\beta$-ISLTRW's DDE]\lbl{fodde}   Fix $\beta=1/\nu$, $\nu\in\lbr2^k;k\in\N\rbr$.
Let $u^x_{\beta,n}(t)=\E \lbk u_0\lpa \S^x_{B,\h}(t)\rpa\rbk$ with $u_0$ as in \eqref{cnd}.
Then $u_{\beta,n}$ solves the following $2\nu$ order differential-difference equation (DDE) on
$\Rp\times\Xnd:$
\beq\lbl{btrwdde}\bc
\df{d u^x_{\beta,n}(t)}{d t}=\sum_{\kappa=1}^{\nu-1}\frac{\Delta_{n}^{\kappa}\unx}{2^{\kappa}t^{1-\kappa/\nu}}{E_{\beta,\kappa}}+\frac{1}{2\nu}\Delta^{\nu}_{n}u_{\beta,n}^{x}(t)
& (t,x)\in(0,\infty)\times\Xnd\\
u_{\beta,n}^x(0)=u_0(x),&x\in\Xnd
\ec
\eeq
where $E_{\beta,\kappa}=\frac{\EP\lpa\Lambda_{\beta}(1)\rpa^{\kappa}}{\kappa!}$.  Moreover, $\Ktxn$ solves \eqref{btrwdde} on $[0,\infty)\times\Xnd$, with
\beq
u_0(x)=\Kzxn=\qzxn=
\bc
1,&x=0\cr
0,&x\neq0.
\ec
\eeq
\elm
\ig 
Just as with BTPs in \cite{Abtp1,Abtp2}, and the BTRW in \cite{Abtbmsie}, we can define a large class of processes containing our $\beta$-ISLTRW and $\beta$-ISLTBM.  We list elements of this class in the discrete setting; their continuous cousins are obtained by replacing the outside $d$-dimensional chain (discrete-valued process) with a $d$-dimensional continuous process.    We now introduce this new class of discrete-valued processes.  Suppressing the $n$ in the lattices $\Xnd$, let
$\Gamma_{\beta}$ be a one-dimensional  at $0$ and let $D^x$ be
an independent $d$-dimensional $\Xd$-valued continuous-time Markov chain starting at
$x$, both defined on a probability space $\OFFtP$.  We call the process
$\DxBt\eqdef D^x(|B_t|)$  a Brownian-time chain (BTC).  A BTRW
is a special case of BTCs in which $D^x$ is a continuous-time random walk.
Excursions-based Brownian-time chains (EBTCs) are obtained from BTCs by breaking up the
path of $|B_t|$ into excursion intervals---maximal intervals $(r, s)$ of time on which
$|B_t| > 0$---and, on each such interval, we pick an independent copy of the Markov
chain $D^x$ from a finite or an infinite collection.  BTCs and EBTCs may be regarded as
canonical constructions to some quite interesting new processes:
\begin{enumerate}
\item  Markov snake chain: when $|B_t|$ increases we pick a new chain $D^x$, we denote this process
by $\DxBsnkchnt$.
\item  $k$-EBTCs: let $D^{x,1},\ldots,D^{x,k}$ be independent copies of $D^x$ starting from point $x\in\Xd$.
On each $|B_t|$ excursion interval, use one of the copies chosen at random.   we denote such a process by
$\DxkBet$.  
When $k=1$ we obtain a BTC.
\item  $\infty$-EBTCs: we use an independent copy of $D^x$ on each $|B_t|$ excursion interval.
This is the $k\to\infty$ of (2),  It is intermediate between (1) and (2).
Here, we go forward to a new independent chain only after $|B_t|$ reaches $0$.
This process is denoted by $\DxiBet$.
\end{enumerate}
\fi
\subsubsection{Lattice $\beta$-ISLTRW SIEs and limits solutions to $\beta$-ISLTBM SIEs}\lbl{limitsol}
The crucial role of the $\beta$-ISLTRW density in our approach to the $\beta$-ISLTBM SIEs\eqref{btpsol}
becomes even clearer from the following definition of our approximating spatially-discretized
equations:
\begin{defn}[Lattice $\beta$-ISLTRW SIEs]\label{btrwsiedef}
By the $\beta$-ISLTRW SIEs associated with the BTBM SIE $\btsie$ we mean
the system $\lbr\isltrwsien\rbr_{n=1}^\infty$ of spatially-discretized stochastic integral equations
on $\Rp\times\Xnd$ given by
\beq\lbl{btrwsie}
\bsp
\uttxn=\sum_{y\in\Xnd}\Ktxyn\uny+\sum_{y\in\Xnd}\int_0^t\Ktsxyn a(\utsyn)\df{d\wsyn}{\sdh},
\end{split}
\eeq
where the $\beta$-ISLTRW density is given by \eqref{btrwdsty}.  For each $n\in\N$, we think of $\{\wtxn; t\ge0\}$ as
a sequence of independent standard
Brownian motions indexed by the set $\Xnd$ (independence within the same lattice).  We also
assume that if $m\neq n$ and $x\in\Xmd\cap\Xnd$ then $\wtxm=\wtxn$, and if $n>m$ and
$x\in\Xnd\sm\Xmd$ then $\wtxm=0$.
\end{defn}
\bnt\lbl{notdetsto}
We will denote the deterministic and the random parts of \eqref{btrwsie} by $\uttxnD$ and $\uttxnR$
(or $\uttxD$ and $\uttxR$ when we suppress the dependence on $n$), respectively,
whenever convenient.
\ent
We define two types of solutions to $\beta$-ISLTRW SIEs: direct solutions and limit solutions.
\bdf[Direct $\beta$-ISLTRW SIE Solutions]\lbl{btrwdsol}
A direct solution to the $\beta$-ISLTRW SIE system $\lbr\isltrwsien\rbr_{n=1}^\infty$ on $\Rp\times\Xnd$
with respect to the Brownian $($in $t$\/$)$ system $\lbr\wtxn;t\ge0\rbr_{(n,x)\in\N\times\Xnd}$ on
the filtered probability space $\OFFtP$ is a sequence of real-valued processes
$\lbr\tilde{U}_{n}\rbr_{n=1}^\infty$  with continuous sample paths in $t$ for each fixed $x\in\Xnd$ and $n\in\N$
such that, for every $(n,x)\in\N\times\Xnd$, $\uttxn$ is $\sFt$-adapted, and equation \eqref{btrwsie} holds
$\P$-a.s.  A solution is said to be strong if $\{\wtxn;t\ge0\}_{(n,x)\in\N\times\Xnd}$ and $\OFFtP$ are fixed a priori$;$
and with
\begin{equation}
\sFt=\sigma\left\{\sigma\left(\wsxn;0\le s\le t,x\in\Xnd,n\in\N\right)\cup
\mathscr N\right\};\quad t\in\Rp,
\label{filt2}
\end{equation}
where $\mathscr N$ is the collection of null sets
$$\left\{O:\exists\,G\in\mathscr{G},O\subseteq G\ \mbox{and}\ \P(G)=0\right\}$$
and where 
$$\mathscr{G}=\sigma\left(\bigcup_{t\ge0}\sigma\left(\wsxn;0\le s\le t,x\in\Xnd,n\in\N\right)\right).$$
A solution is termed weak if we are free to
choose $\OFFtP$ and the Brownian system on it and without requiring $\sFt$ to
satisfy $\eqnref{filt2}$.  Replacing $\Rp$ with $\T:=[0,T]$---for some $T>0$ in the
above, we get the definition of  a solution to the $\beta$-ISLTRW SIE system $\lbr\isltrwsien\rbr_{n=1}^\infty$ on
$\T\times\Rd$.
\edf
The next type of $\beta$-ISLTRW SIE solutions we define is the first step in our K-martingale approach  of \cite{Abtbmsie}, which we recall in \secref{Kmartsec}.   By first reducing $\isltrwsien$ to the simpler finite dimensional noise setting, it takes full advantage of the notion of $\beta$-ISLTRW SIEs limit solutions to $\beta$-ISLTBM SIEs.
\bdf[Limit $\beta$-ISLTRW SIE Solutions]\lbl{btrwlsol}
Let $l\in\N$.  By the $l$-truncated $\beta$-ISLTRW SIE on $\Rp\times\Xnd$ we mean the $\beta$-ISLTRW SIE obtained from \eqref{btrwsie}
by restricting the sum in the stochastic term to the finite $d$-dimensional lattice $\Xnldf:=\Xnd\cap\lbr[-l,l]^d;l\in\N\rbr$
and leaving unchanged the deterministic term $\uttxnD$:
\beq\lbl{trnctdsie}
\uttxnlf=\bc
\ds\uttxnD
&+\ \ds\sum_{y\in\Xnldf}\int_0^t\kappa^{x,y}_{\h,s,t}\lpa\utsynlf\rpa{d\wsyn}; x\in\Xnldf,\\
\ds\uttxnD;&x\in\Xnd\setminus\Xnldf
\ec
\eeq
where $$\kappa^{x,y}_{\beta,\h,s,t}\lpa\utrynlf\rpa:=\df{\Ktsxyn}{\sdh} a(\utrynlf),\quad \forall r,s<t.$$  We denote \eqref{trnctdsie} by
$\isltrwsienl$.  Fix $n\in\N$, a solution to the system of truncated $\beta$-ISLTRW SIEs $\lbr\isltrwsienl\rbr_{l=1}^\infty$
on $\Rp\times\Xnd$ with respect to the Brownian $($in $t$\/$)$ system $\lbr\wtxn;t\ge0\rbr_{x\in\Xnd}$ on
the filtered probability space $\OFFtP$ is a sequence of real-valued processes $\lbr\tilde{U}_{\beta,n,l}\rbr_{l\in\N}$
with continuous sample paths in $t$ for each fixed $x\in\Xnd$ and $l\in\N$,
such that, for every $(l,x)\in\N\times\Xnd$, $\uttxnlf$ is $\sFt$-adapted, and equation \eqref{trnctdsie} holds
$\P$-a.s.  We call $\utn$ a limit solution to the $\beta$-ISLTRW SIE \eqref{btrwsie} if $\utn$ is a limit of the truncated
solutions $\utnlf$ (as $l\to\infty$).    If desired, we may indicate the limit type (a.s., in $L^p$, weak, \ldots, etc).
\edf
\brm\lbl{smuthdet}
In both \eqref{trnctdsie} and \eqref{btrwsie}, $\uttxnD=\E \lbk u_0\lpa \S^x_{B,\h}(t)\rpa\rbk$.  So, by \lemref{fodde},
$\uttxnD$ is differentiable in time $t$ and satisfies \eqref{btrwdde}.
Also, using linear interpolation, we can extend the definition of an already continuous-in-time process $\uttxn$
on $\Rp\times\Xnd$, so as to obtain a continuous process on $\Rp\times\Rd$, for each
$n\in\N$, which we will also denote by $\uttxn$.    Henceforth,
any such sequence $\{\utn\}$ of interpolated $\utn$'s will be called a continuous
or an interpolated solution to the system $\lbr\isltrwsien\rbr_{n=1}^\infty$.  Similar comments
apply to solutions of the truncated $\isltrwsienl$.
\erm
We now define solutions to $\isltsie$ based entirely on their approximating $\lbr\isltrwsien\rbr$, through their limit.
Since we defined direct and limit solutions to $\isltrwsien$, for each fixed $n$, we get two types of $\beta$-ISLTRW SIEs limit solutions to $\isltsie$:
direct $\beta$-ISLTRW SIEs limit solutions and $\beta$-ISLTRW SIE double limit solutions.  The ``double'' in the second type of solutions
reminds us that we are taking two limits, one from truncated to nontruncated fixed lattice (as $l\to\infty$)
and the other limit is taken as the lattice mesh size shrinks to zero (as $\h\searrow0$ or equivalently as $n\nearrow\infty$).

\begin{defn}[$\beta$-ISLTRW SIEs limits solutions to $\isltsie$]\label{limitsolns}
We say that the random field $U$ is a $\beta$-ISLTRW SIE limit solution
to $\isltsie$ on $\Rp\times\Rd$ iff there is a solution $\{\uttxn\}_{n\in\N}$
to the lattice SIE system $\lbr\isltrwsien\rbr_{n\in\N}$ on a probability space
$\OFFtP$ and with respect to a Brownian system $\{\wtxn;t\ge0\}_{(n,x)\in\N\times\Xnd}$
such that $U$ is the limit or a modification of the limit of $\lbr\utn\rbr_{n\in\N}$ $($or a subsequence thereof).   A $\beta$-ISLTRW SIE limit solution $U$ is called a direct $\beta$-ISLTRW SIEs limit solution or a $\beta$-ISLTRW SIEs double limit solution
according as $\{\uttxn\}_{n\in\N}$ is a sequence of direct or limit solutions to $\lbr\isltrwsien\rbr_{n\in\N}$.   The limits may be taken
in the a.s., probability, $L^p$, or weak sense\footnote{When desired, the types of the solution and the limit are explicitly stated
$($e.g., we say strong $($weak\/$)$ $\beta$-ISLTRW SIEs weak,
in probability, $L^p(\Omega)$, or a.s.~limit solution to indicate that the solution to the approximating
SIEs system is strong $($weak\/$)$
and that the limit of the SIEs is in the weak, in the probability, in the $L^p(\Omega)$, or in the  a.s.~sense, respectively\/$)$.  Of course, we may also take limits in any other suitable sense.}.  We say that uniqueness in law holds if whenever $U^{(1)}$ and $U^{(2)}$ are $\beta$-ISLTRW SIEs limit solutions they have the same law.   We say that pathwise uniqueness holds for $\beta$-ISLTRW SIEs limit solutions if whenever $\lbr\uotn\rbr$ and $\lbr\utttn\rbr$ are lattice SIEs solutions on the same probability space
and with respect to the same Brownian system, their limits $U^{(1)}$ and $U^{(2)}$ are indistinguishable.
\end{defn}
\bcm
\begin{rem}\label{defnrem}
\end{rem}
\fi
\subsubsection{Second main theorem: the lattice-limits solutions case} \lbl{mres}
We can now state our second main result of the paper.  The following theorem gives our lattice-limits solutions result for $\isltsie$ under the non-Lipschitz conditions \eqref{cnd} on $a$.  Our limits solutions result under Lipschitz conditions is stated in \thmref{latlimlip}\footnote{The type of limit solutions in the Lipschitz case  is direct limit solutions as opposed to the double limit solution in \thmref{mainthm2}.}, which is proved in Appendix \ref{appB}.  
\bfr
\bthm[Spatio-temporal regularity and third dimension maximality: lattice-limits solutions] \lbl{mainthm2}  Fix $\beta=1/\nu$, $\nu\in\lbr2^k;k\in\N\rbr$.
 Assume the conditions \beq\lbl{cnd}\tag{\mbox{NLip}}
\bc
(a) &a(u) \mbox { is continuous in $u$;}\quad u\in\R,\\
(b) &a^2(u)\le C (1+u^2);\quad u\in\R,\\
(c) &u_0\in \mathrm{C}_{b}^{2\nu-2,\gamma}(\Rd;\R)\mbox{ and nonrandom } ,\ \forall\ d=1,2,3.
\ec
\eeq
 hold.  Then,
there exists a $\beta$-ISLTRW SIE double weak-limit solution to $\isltsie$, $U$, such that
$U(t,x)$ is $L^p(\Omega,\P)$-bounded on $\T\times\Rd$ for every $p\ge2$ and $U_{\beta}\in\H^{\lpa{\tf{2\nu-d}{4{\nu}}}\rpa^{-},{\lpa\tf{4-d}{2}\wedge 1\rpa}^{-}}(\T\times\Rd;\R)$ for every $d=1,2,3$.   \ethm
\efr
\brm\lbl{scnv}
Of course, we can use change of measure---as we did in our earlier work on Allen-Cahn SPDEs and other second order SPDEs (see e.g.~\cite{Acom,Acom1,Acom2} and all our change of measure references in \cite{Asdde2} for results and conditions)---to transfer existence, uniqueness, and law equivalence results between $\isltsie$ and the $\beta$-ISLTBM SIE with measurable drift $\isltsieb$: 
\beq\lbl{btpsiedrft}
\bsp
U_{\beta}(t,x)=&\intrd\KIBtxy \uny dy+ \intrdzt\KIBtsxy b(\usy)ds dy
\\&+ \intrdzt\KIBtsxy a(\usy)\sW(ds\times dy),
\end{split}
\eeq
under the same conditions on the drift/diffusion ratio.  If it is desired to investigate $\isltsieb$ on a bounded domain in $\Rd$ with a regular boundary, we simply replace the $\beta$-ISLTBM density $\KIBtxy$ in \eqref{btpsiedrft} with its boundary-reflected or boundary-absorbed version (the $\beta$-ISLTBM density in which the outside $d$-dimensional BM is either reflected or absorbed at the boundary).
\erm
The proof of \thmref{mainthm2} under the conditions \eqref{cnd} is neither standard nor straightforward---even after obtaining the new non-trivial spatio-temporal regularity estimates (in  \lemref{4thQinequality} and  \lemref{3rdQinequality} below) on the unconventional kernel $\KIBtxy$.   This is because standard techniques, like the classical martingale problem approach, do not apply directly to kernel equations like the $\beta$-ISLTBM SIE $\isltsie$ or its discretized version $\isltrwsien$ under \eqrf{cnd}.  This leads us to use our K-martingale approach, introduced in \cite{Abtbmsie}.

\section{Key estimates}\lbl{pfs}
\subsection{Density regularity estimates and third dimension maximality}\lbl{destm}
The first set of estimates\footnote{As is customary, all constants may change their value from one line to the next without changing their notation.  Also, to simplify notation, we will often suppress the dependence on $\beta$ without further notice.  We will denote the Euclidean norm on $d$-dimensional spaces by $\lab\cdot\rab$.} we need are bounds on the square of the $\beta$-inverse-stable-L\'evy-time Brownian motion density $\KIBtxy$ and its associated lattice $\beta$-inverse-stable-L\'evy-time random walk density $\Ktxn$ and their temporal and spatial differences.  We obtain these estimates for both kernels simultaneously. The method of proof is to reduce, via an asymptotic argument, these estimates for the $\beta$-ISLTRW to the corresponding ones for the $\beta$-ISLTBM density $\KIBtx$ and perform the computations in the continuous setting of the $\beta$-ISLTBM.  Since all the results in this part hold for all $n\ge N^*$ (equivalently for all $\h\le\delta_{N^*}$) for some positive integer $N^*$, we will suppress the dependence on $n$, except when it is needed or helpful, to simplify the notation.   Also, whenever we need these estimates, we assume that  $n\ge N^*$ without explicitly stating it every time; and when we do, we let\footnote{We adopt these simplifications with lattice computations throughout the paper.} 
\beq\lbl{star}
\Ns:=\lbr n\in\N;n\ge N^*\rbr
\eeq

We start by observing that in the classical setting of Brownian motion and its discretized version continuous-time random walk on $\Xnd=\h\Zd$, we have  the following well known asymptotic result relating their densities (see e.g., \cite{SZ})
\beq\lbl{as1}
\qtxynint\sim\ptxy\hd\mbox{ as } n\to\infty\mbox{ (as $\h\to0$); }\forall t>0,\ x,y\in\Rd, 
\eeq
where for each $x\in\Rd$ we use $[x]_{\h}$ to denote the element of $\Xnd$ obtained by replacing each coordinate $x_{i}$ with $\h$ times the integer part of $\h^{-1}x_{i}$, and $a_n\sim b_n$ as $n\to\infty$ means $a_n/b_n\to1$ as $n\to\infty$.
Now, for every continuous and bounded $\un:\Rd\to\R$, we have 
\beq\lbl{minusx}
\lim_{\h\searrow0}\sum_{y\in\Xnd\setminus\lbr x\rbr} \KIBtxy \uny\hd=\intrd\KIBtxy\uny dy;\  t>0,\ x\in\Rd,\ d\ge1,
\eeq
and by the dominated convergence theorem we obtain
\beq\lbl{BTRWsol2BTBMsol}
\bsp
&\lim_{\h\searrow0}\lab\sum_{y\in\Xnd}\Ktxynint\uny-\sum_{y\in\Xnd\setminus\lbr x\rbr} \KIBtxy \uny\hd\rab\\
&=\lab\int_{0}^{\infty}\lbr\lim_{\h\searrow0}\sum_{y\in\Xnd\setminus\lbr x\rbr}\lbk\qsxynint-\psxy\hd\rbk\uny\rbr\KItzs ds\rab=0
\end{split}
\eeq
for $t>0$, $x\in\Rd$, and $d\ge1$; since, by \eqref{as1}, 
$$\lim_{\h\searrow0}\sum_{y\in\Xnd}\qsxynint\uny=\lim_{\h\searrow0}\sum_{y\in\Xnd}\psxy\uny\hd=\intrd\psxy\uny dy$$
for every $(s,x)\in(0,\infty)\times\Rd$.  We then straightforwardly get the following result. 
\blm\lbl{btrwasbtp}  For every continuous and bounded $\un:\Rd\to\R$ and for every $d\ge1$ 
\beq\lbl{ddegreen2pdegreen}
\lim_{\h\searrow0}\sum_{y\in\Xnd}\Ktxynint\uny=\intrd\KIBtxy\uny dy;\forall (t,x)\in(0,\infty)\times\Rd,
\eeq
and the following asymptotic relation holds between the $\beta$-ISLTBM and $\beta$-ISLTRW densities: 
\beq\lbl{btrwbtp}
\Ktxynint\sim\KIBtxy\hd \mbox{ as } n\to\infty\ (\mbox{as }\h\to0);\ t>0,\ x,y\in\Rd,\ x\neq y.
\eeq
\elm
\brm\lbl{ddetopde}
Equation \eqref{ddegreen2pdegreen} confirms the intuitively clear fact that the kernel form of the $\beta$-ISLTRW DDE \eqref{btrwdde} 
converges pointwise---as $\h\searrow0$---to the kernel form of its continuous version, the $\beta$-ISLTRW PDE in \cite{Nanesd,AN}.  We also remind the reader that the right hand side of \eqref{ddegreen2pdegreen} is in $\mathrm{C}^{1,2\beta^{-1}}$ for all $(t,x)\in(0,\infty)\times\Rd$ under the $\un$ conditions in \eqref{cnd}.  
 \erm

Our first regularity lemma for the densities is now stated.  It implies, among other things, that 
there is a considerable smoothing effect of $\KIBsx$ as $\beta$ gets smaller; however it also implies that our SIEs don't possess  random field solutions beyond the third spatial dimension, no matter how small $\beta$ gets. 
\begin{lem}[Smoothing and third dimension maximality]\label{2ndQinequality}
There are constants $C$ and $\tilde{C}$, depending only on $d$ and $\beta=1/\nu$, $\nu\in\lbr2^k;k\in\N\rbr$, and a $\d^*>0$ such that for all $\d\le\d^*$
\[\intrd\lbk\KIBtx\rbk^2dx={{C}{t^{\tf{-d}{2{\nu}}}}}\mbox{ and }\sum_{x\in\Xd} \lbk\Ktx\rbk^2 \leq {\tilde{C}{\dd}{t^{\tf{-d}{2{\nu}}}}};\]for all $t>0,\ d=1,2,3$.
{Hence, } 
\[\int_0^t\intrd\lbk\KIBsx\rbk^2dx ds={Ct^{\tf{2\nu-d}{2{\nu}}}}\mbox{ and } \int_0^t\sum_{x\in\Xd} \lbk\Ksx\rbk^2  ds \leq {\tilde{C}\dd t^{\tf{2\nu-d}{2{\nu}}}};\]
for all $t>0,\ d=1,2,3$.  In addition, $\intrd\lbk\KIBtx\rbk^2dx=\int_0^t\intrd\lbk\KIBsx\rbk^2dx ds=\infty,$ for all $d\ge4$.
\end{lem}
\bpf
First, fix an arbitrary $\beta^{-1}=\nu\in\lbr 2^{k},k\in\N\rbr$.
Using the definition of $\KIBtx$, \lemref{btrwasbtp} and \lemref{conkitbtbmisltbm} here together with Lemma 3.1 and Lemma 3.2 in \cite{Nanesd} we obtain\footnote{Recall that we are using the convention $\int_{\Rp^{0}}f(s)ds=f(s)$, for every $f$.}
\beq\lbl{allcomp}
\bsp
&\lim_{\d\searrow0}\sum_{x\in\Xd} \df{\lbk\Ktx\rbk^2}{\dd}=\intrd\lbk\KIBtx\rbk^{2} dx\\
&=\int_0^\infty\int_0^\infty\lbk\int_{\Rd} \psoxz \puoxz dx\rbk \KItzso \KItzuo ds_{1}du_{1}\\
&=\int_0^\infty\int_0^\infty\lbk\df{1}{\lbk2\pi (s_{1}+u_{1})\rbk^{d/2}}\rbk \KItzso \KItzuo ds_{1}du_{1}\\
&=\lbr\int_0^\infty\int_0^\infty\lbk\df{2^{k}}{\lbk2\pi (s_{1}+u_{1})\rbk^{d/2}}\rbk\right.
\\&\times \lpa\int_{(0,\infty)^{k-1}} \ptzsk\prod_{i=0}^{k-2}\ptzskiimo ds_{2}\cdots ds_{k}\rpa\\
&\left.\times \lpa\int_{(0,\infty)^{k-1}} \ptzuk\prod_{i=0}^{k-2}\ptzukiimo du_{2}\cdots du_{k}\rpa ds_{1}du_{1}\rbr
\end{split}
\eeq 
Gathering the two inside integrals and transforming to polar coordinates $(s_{i},u_{i})\mapsto(\rho_{i},\theta_{i})$, letting  $\underline{\rho}=(\rho_{1},\ldots,\rho_{k})$ and $\underline{\theta}=(\theta_{1},\ldots,\theta_{k})$, and noticing that all $\rho_{i}$ for $i=2,3,\ldots,\rho_{k}$ cancel when $k\ge2$;  equation \eqref{allcomp} becomes\footnote{Equation \eqref{sumofsqofisltbmden} is the reason for the third spatial dimension maximality.}
\beq\lbl{sumofsqofisltbmden}
\bsp
& C\int\limits_{(0,{\pi/2})^{k}}\int\limits_{(0,\infty)^{k}}\df{e^{\frac{-\rho_{k}^2}{4t}}\ds\prod_{i=0}^{k-2}e^{-\lbk\frac{\rho_{k-i-1}^2\cos^{2}(\theta_{k-i-1})}{4\rho_{k-i}\cos(\theta_{k-i})}+\frac{\rho_{k-i-1}^2\sin^{2}(\theta_{k-i-1})}{4\rho_{k-i}\sin(\theta_{k-i})}\rbk}}{\rho_{1}^{\frac{d}{2}-1}t
\lbk\sin(\theta_{1})+\cos(\theta_{1})\rbk^{\frac{d}{2}}\ds\prod_{i=0}^{k-2}\sqrt{\sin(\theta_{k-i})\cos(\theta_{k-i})}}d\underline{\rho} d\underline{\theta}\\
&=\bc{C}{t^{\tf{-d}{2{\nu}}}};& d=1,2,3,\cr\infty;&d\ge4.\ec
\end{split}
\eeq
Then there is a $\d^{*}>0$ such that, whenever $\d\le\d^{*}$, we obtain 
$$\df1\dd\sum_{x\in\Xd} \lbk\Ktx\rbk^2\le{\tilde{C}{t^{\tf{-d}{2{\nu}}}}};\ d=1,2,3,$$
with a finite constant $\tilde{C}>C$.  The last assertion of the lemma trivially follows upon integration over the time interval $(0,t]$.
\epf
The following lemma is key to our H\"older regularity result in time.  We give a probabilistically-flavored proof using the notion of $2$-$\beta$-inverse-stable-L\'evy-times random walk and $2$-$\beta$-inverse-stable-L\'evy-times Brownian motion given below.  
\begin{lem}[Kernel temporal regularity]\label{4thQinequality}
There is a constant $C$, depending only on $d$ and $\beta=1/\nu$, $\nu\in\lbr2^k;k\in\N\rbr$, and a $\d^*>0$ such that for $\d\le\d^*$ 
\beq\lbl{tmpkernel}
\bc
\ds\int_0^t\intrd{\lbk\KIBtsx - \KIBrsx\rbk}^2 dx ds \leq C(t-r)^{\tf{2\nu-d}{2{\nu}}},\\
\ds\int_0^t\sum_{x\in\Xd} {\lbk\Ktsx - \Krsx\rbk}^2 ds \leq {C}\dd(t-r)^{\tf{2\nu-d}{2{\nu}}},\ec
\eeq
for $0<r<t$ and $d=1,2,3$, with the convention that $\Ktx=0=\KBtx$ if $t<0$.
\end{lem}
\bpf
We will prove that 
\beq\lbl{tmphl}
\int_0^t\sum_{x\in\Xd} {\lbk\Ksptmrx - \Ksx\rbk}^2 ds \leq C\dd(t-r)^{\tf{2\nu-d}{2{\nu}}}; \quad d=1,2,3.
\eeq
for all $\delta\le\delta^{*}$, for some $\d^{*}>0$, simultaneously with its corresponding $\beta$-inverse-stable-L\'evy-time Brownian motion density statement.  The first step is to show the identity 
\beq\lbl{ido}
\bsp
&\sum_{x\in\Xd} {\lbk\Ksptmrx - \Ksx\rbk}^2
\\  & \ =\KIIsptmrz+\KIIsz-2\KIIsptmrsz
\end{split}
\eeq
where 
\beq
\KIIuvz=\int_{0}^{\infty}\int_{0}^{\infty}\qroprtz\KIuzro\KIvzrt dr_{1} dr_{2}
\eeq 
is the density of the $2$-$\beta$-inverse-stable-L\'evy-times random walk
\beq\lbl{twobtrwdef}
\S^0_{\Lambda_{\beta}^{(1)},\Lambda_{\beta}^{(2)},\h}(u,v)\eqdef S^0_{\h}\lpa\Lambda_{\beta}^{(1)}(u)+ \Lambda_{\beta}^{(2)}(v) \rpa;\quad 0\le u,v<\infty,
\eeq
in which the $d$-dimensional random walk $S^0_{\h}$ (on $\Xnd$) and the two identically-distributed one-dimensional processes $\Lambda_{\beta}^{1}$ and $\Lambda_{\beta}^{2}$ are all independent.
But, 
\beq\lbl{idz}
\bsp
&\sum_{x\in\Xd} \Kux\Kvx\\&=\int_{0}^{\infty}\int_{0}^{\infty}\lbk\sum_{x\in\Xd}\qrox\qrtx\rbk\KIuzro\KIvzrt dr_{1} dr_{2}\\
\\&=\int_{0}^{\infty}\int_{0}^{\infty}\qroprtz\KIuzro\KIvzrt dr_{1} dr_{2}=\KIIuvz.
\end{split}
\eeq
The identity \eqref{ido} immediately follows from \eqref{idz}.  Similarly,  we get the corresponding identity for the $\beta$-inverse-stable-L\'evy-time Brownian motion setting
\beq\lbl{idtwo}
\bsp
&\int_{\Rd} {\lbk\KIBsptmrx - \KIBsx\rbk}^2 dx
\\\ &=\KIIBsptmrz+\KIIBsz-2\KIIBsptmrsz
\end{split}
\eeq
where 
\beq
\KIIBuvz=\int_{0}^{\infty}\int_{0}^{\infty}\proprtz\KIuzro\KIvzrt dr_{1} dr_{2}
\eeq 
is the density of the $2$-$\beta$-inverse-stable-L\'evy-times Brownian motion
\beq\lbl{2btbmdef}
\X^0_{\Lambda_{\beta}^{(1)},\Lambda_{\beta}^{(2)}}(u,v)\eqdef X^0\lpa\Lambda_{\beta}^{(1)}(u)+ \Lambda_{\beta}^{(2)}(v) \rpa;\quad 0\le u,v<\infty,
\eeq
in which the $d$-dimensional BM $X^{0}$ and the two identically-distributed one-dimensional processes $\Lambda_{\beta}^{(1)}$ and $\Lambda_{\beta}^{(2)}$ are all independent.
Using the identities \eqref{ido} and \eqref{idtwo}, along with a similar asymptotic argument to the one we used in the proof of 
\lemref{2ndQinequality} together with the dominated convergence theorem, yield
\beq\lbl{hlbd}
\bsp
&\lim_{\d\searrow0}\frac{1}{\dd}\lbk\int_{0}^{t}\KIIsptmrz ds + \int_{0}^{t} \KIIsz ds-2\int_{0}^{t} \KIIsptmrsz ds\rbk
\\&=\lim_{\d\searrow0}\int_0^t\sum_{x\in\Xd} \frac{{\lbk\Ksptmrx - \Ksx\rbk}^2}{\dd} ds\\&=\int_{0}^{t}\int_{\Rd} {\lbk\KIBsptmrx - \KIBsx\rbk}^2 dxds
\\&=\lbk\int_{0}^{t}\KIIBsptmrz ds + \int_{0}^{t} \KIIBsz ds-2\int_{0}^{t} \KIIBsptmrsz ds\rbk
\\&=\lbk\int_{0}^{t}\tilde{\K}_{2s+2(t-r)} ds + \int_{0}^{t} \tilde{\K}_{2s} ds-2\int_{0}^{t} \tilde{\K}_{2s+(t-r)}ds\rbk
\\&=\lbk\int_{0}^{\tf{t-r}{2}}\tilde{\K}_{2s} ds-\int_{\tf{t-r}{2}}^{t-r}\tilde{\K}_{2s} ds-\int_{t}^{t+\tf{t-r}{2}}\tilde{\K}_{2s} ds
+\int_{t+\tf{t-r}{2}}^{2t-r}\tilde{\K}_{2s} ds\rbk
\end{split}
\eeq
for $d=1,2,3$, where $\tilde{\K}_{w}$ is defined in terms of $\KIIBuvz$ by the relation
\beq\lbl{trans}
\bsp
&\tilde{\K}_{w}=\KIIBuvz\iff w=u+v\mbox{ and }(u,v)\mbox{ has one of the forms} \\&(u,v)=(a,a)\mbox{ or }(u,v)=(a+b,a)\mbox{ or }(u,v)=(a,a+b)\mbox{ for some}\ a,b\ge0.
\end{split}
\eeq 

We observe that
\beq\lbl{Ktildebtrw}
\bsp
\tilde{\K}_{2u}=\KIIBuz&=\int_{0}^{\infty}\int_{0}^{\infty}\proprtz\KIuzro\KIuzrt dr_{1} dr_{2}\\
&=\int_{0}^{\infty}\int_{0}^{\infty}\lbk\int_{\Rd}\prox\prtx dx\rbk\KIuzro\KIuzrt dr_{1} dr_{2}\\
&=\int_{\Rd} \lbk\KIBux\rbk^{2}dx={C} {u^{\tf{-d}{2{\nu}}}};\quad d=1,2,3
\end{split}
\eeq
The last assertion follows from the computation in \eqref{allcomp} and \eqref{sumofsqofisltbmden}.  It is clear then that $\tilde{\K}_{2u}$ is decreasing in $u$, for every $\nu=1/\beta\in\lbr2^k;k\in\N\rbr$.  Thus, the sum of the last three terms of the \eqref{hlbd} is $\le0$.   This and \eqref{Ktildebtrw} give us \eqref{tmphl} for all $\delta\le\delta^{*}$, for some $\d^{*}>0$ and for some constant $C>0$, together with its corresponding $\beta$-inverse-stable-L\'evy-time Brownian motion density statement; and \lemref{4thQinequality} follows at once.
\epf
The following spatial difference second moment inequality for the $\beta$-ISLTRW and $\beta$-ISLTBM densities reflects their critical spatial-regularizing effect on our solutions.  The following lemma captures the surprising fact that we cannot improve on the spatial regularity of the BTBM SIE by decreasing  $\beta$ below $1/2$.  This implies the maximality of the BTBM SIEs spatial regularity among the family of $\beta$-ISLTBM SIE family.
\begin{lem}[Kernel spatial regularity]\label{3rdQinequality}
Let $\beta\in\lbr1/2^k;k\in\N\rbr$ and define the intervals 
$$I_{d}=\bc (0,1];&d=1,\\ 
(0,1);&d=2,\\
(0,\frac12);&d=3.
\ec$$ 
For any given positive numbers $\lbr\alpha_{d}\in I_{d}\rbr_{d=1}^{3}$, there exists a constant $C$ depending only on $\beta$, $d$ and $\lbr\alpha_{d}\rbr_{d=1}^{3}$, and  a $\d^*>0$ such that for $\d\le\d^*$  
\beq\lbl{sptlkernel}
\bc
\ds\int_0^t\intrd {\lbk\KIBsx - \KIBsxpz\rbk}^2 dx ds \leq C |z|^{2\alpha_{d}}t^{p(\alpha_{d},\beta)},\\
\ds\int_0^t\sum_{x\in\Xd} {\lbk\Ksx - \Ksxpz\rbk}^2 ds \leq C\delta^{d} |z|^{2\alpha_{d}}t^{p(\alpha_{d},\beta)},
\ec
\eeq
for $t>0$, where $0<C<\infty$ and $0\le p(\alpha_{d},\beta)<1$ for every $\alpha_{d}\in I_{d}$ for $ d=1,2,3$ and for every $\beta\in\lbr1/2^k;k\in\N\rbr$.  
\end{lem}
\brm\lbl{spatlkerrem}  For a given $\beta^{-1}\in\lbr2,3,4,\ldots\rbr$, and on any compact time interval $\T=[0,T]$, the inequality \eqref{sptlkernel} may---for any given value $\alpha_{d}$---be rewritten as 
\beq\lbl{sptlkernel2}
\bc
\ds\int_0^t\sum_{x\in\Xd} {\lbk\Ksx - \Ksxpz\rbk}^2 ds \leq 
\tilde{C}\delta^{d} |z|^{2\alpha_{d}}
;&\\
\ds\int_0^t\intrd {\lbk\KIBsx - \KIBsxpz\rbk}^2 dx ds \leq \tilde{C} |z|^{2\alpha_{d}};
\ec
\eeq
where, for each $d=1,2,3$ $$\tilde{C}=C\sup_{\substack{\alpha_{d}\in I_{d},\\ \beta\in\lbr1/2^k;k\in\N\rbr.}}T^{p(\alpha_{d},\beta)}<\infty$$ also depends on $T$ in \eqref{sptlkernel2}.  
\erm
\bpf  Let $\beta=1/2^{k}$ for $k\in\N$.  
Starting with the $L^{2}$ estimate involving the spatial difference of the $\beta$-ISLTBM density in \eqref{sptlkernel}, letting $u_{1}=r_{2}$, using the polar transformation $(r_{i},u_{i})\mapsto(\rho_{i},\theta_{i})$,  letting  $\underline{\rho}=(\rho_{1},\ldots,\rho_{k})$ and $\underline{\theta}=(\theta_{1},\ldots,\theta_{k})$, and noticing that all $\rho_{i}$ for $i=2,3,\ldots,\rho_{k}$ cancel when $k\ge2$, we have
\beq\lbl{sptlebaktar}
\bsp
&\int_0^t\intrd {\lbk\KIBsx - \KIBsxpz\rbk}^2dx ds
\\&=\int_0^t\lbk\int_0^\infty\int_0^\infty\intrd\prod_{i=1}^{2}{\lpa\prix - \prixpz\rpa} \KIszri dxdr_{1}dr_{2}\rbk ds
\\&=\int_{0}^{t}\int_0^\infty\int_{0}^{\infty}{\lpa2\propuoz - 2\propuozz\rpa} \KIszro\KIszuo  dr_{1}du_{1}ds
\\&=2\int_{0}^{t}\int_0^\infty\int_{0}^{\infty}\frac{1-e^{-\frac{|z|^{2}}{2(r_{1}+u_{1})}}}{\lbk2\pi (r_1+u_1)\rbk^{d/2}}
\\&\times \lpa\int_{(0,\infty)^{k-1}} \pszrk\prod_{i=0}^{k-2}\ptzrkiimo du_{2}\cdots du_{k}\rpa  dr_{1}du_{1}ds
\\&\times \lpa\int_{(0,\infty)^{k-1}} \pszuk\prod_{i=0}^{k-2}\ptzukiimo du_{2}\cdots du_{k}\rpa  dr_{1}du_{1}ds
\\&\le C\int_{0}^{t}\int\limits_{(0,\frac\pi2)^{k}}\int\limits_{(0,\infty)^{k}}\df{\lpa1-e^{-\frac{|z|^{2}}{2\rho_{1}}}\rpa e^{\frac{-\rho_{k}^2}{4s}}\ds\prod_{i=0}^{k-2}e^{-\lbk\frac{\rho_{k-i-1}^2}{4\rho_{k-i}}\rbk}}{\rho_{1}^{\frac d2-1}s\lbk\sin(\theta)+\cos(\theta)\rbk^{\frac d2}\ds\prod_{i=0}^{k-2}\sqrt{\sin(\theta_{k-i})\cos(\theta_{k-i})}} d\underline{\rho} d\underline{\theta} ds
\\&\le C \int_{0}^{t}\int\limits_{(0,\infty)^{k}}\df{\lpa1-e^{-\frac{|z|^{2}}{2\rho_{1}}}\rpa e^{\frac{-\rho_{k}^2}{4s}} \ds\prod_{i=0}^{k-2}e^{-\lbk\frac{\rho_{k-i-1}^2}{4\rho_{k-i}}\rbk}}{\rho_{1}^{\frac d2-1}s}d\underline{\rho}  ds
\\&\le C \int_{0}^{t}\int\limits_{(0,\infty)^{k}}\df{{{|z|^{2\alpha}}} e^{\frac{-\rho_{k}^2}{4s}} \ds\prod_{i=0}^{k-2}e^{-\lbk\frac{\rho_{k-i-1}^2}{4\rho_{k-i}}\rbk}}{\rho_{1}^{\alpha+\frac d2-1}s} d\underline{\rho} ds\\&
\le\bc
 C_{1} |z|^{2\alpha}t^{p_{1}(\alpha,\beta)}; &d=1,\alpha\in(0,1], \\
 C_{2} |z|^{2\alpha}t^{p_{2}(\alpha,\beta)}; &d=2,\alpha\in(0,1), \\
 C_{3}|z|^{2\alpha}t^{p_{3}(\alpha,\beta)}; &d=3,\alpha\in(0,\tf{1}{2}),
\ec
\end{split}
\eeq
for some finite constants $C_{i}$, $i=1,2,3$, where $C_{2}$ and $C_{3}$ depend on $\alpha$\footnote{See \remref{spatlkerrem} in \cite{Abtbmsie} for a detailed discussion in the BTBM case $\beta=1/2$.}, and where we have used the simple facts that $\min_{0\le\theta\le\pi/2}\lbk\sin(\theta)+\cos(\theta)\rbk=1$ and that $1-e^{-u}\le u^{\alpha}$ for $u\ge0$ and $0<\alpha\le1$.  This proves the $L^{2}$ estimate for the $\beta$-ISLTBM density in \eqref{sptlkernel}.  Then, an asymptotic argument similar to the one in the proofs of \lemref{2ndQinequality} and \lemref{4thQinequality} yields 
\beq\lbl{ha}
\lim_{\d\searrow0}\int_0^t\sum_{x\in\Xd}\frac{{\lbk\Ksx - \Ksxpz\rbk}^2}{\dd} ds=\int_0^t\intrd {\lbk\KIBsx - \KIBsxpz\rbk}^2 dx ds,
\eeq
together with the desired $\beta$-ISLTRW density $L^{2}$ estimate in \eqref{sptlkernel} for all $\d\le\d^{*}$, for some $\d^{*}>0$, with possibly different constants.
\epf

\subsection{Spatio-temporal estimates for $\beta$-ISLTRW and $\beta$-ISLTBM SIEs}\lbl{regandtight}
 In this subsection, and assuming only the less-than-Lipschitz conditions \eqref{cnd} on $a$---together with a temporary moment condition---we obtain spatial and temporal differences moments estimates that are crucial in obtaining the regularity
 of the $\beta$-ISLTRW SIE $\isltrwsien$ for each fixed $n\in\Ns$ (see \eqref{star}), the tightness of the $\beta$-ISLTRW SIEs sequence $\lbr\isltrwsien\rbr_{n\in\Ns}$, as well as the H\"older regularity for their limiting $\beta$-ISLTBM SIE.  To make it more convenient for the proof of our first main result in the direct solution case, \thmref{mainthm1}, we include the corresponding spatio-temporal statements for the  $\beta$-ISLTBM SIE in the same lemmas, together with those for their lattice cousins.  
 
 Fix $n\in\Ns$, and assume $\utn$ solves $\isltrwsien$ in \eqref{btrwsie} and $U_{\beta}$ solves the $\beta$-ISLTBM SIE $\isltsie$ in \eqrf{isltbmsie}.   Suppressing the dependence on $n$,
 let $\tilde{M}_{\beta,q}(t) = \sup_x \mathbb{E}|\tilde{U}_{\beta}^x(t)|^{2q}$, and ${M}_{\beta,q}(t) = \sup_x \mathbb{E}|{U}_{\beta}(t,x)|^{2q}$ for $q\ge1$ and $\beta\in\lbr1/2^{k}, k\in\N\rbr$.  Writing $\ut$ and $U_{\beta}$ in terms of
their deterministic and random parts  $\uttx=\uttxD+\uttxR$ and $U_{\beta}(t,x)=U_{\beta,D}(t,x)+U_{\beta,R}(t,x)$, we observe that $\uttxD$ is smooth in time by \lemref{fodde} and $U_{\beta,D}(t,x)$ is smooth in time and space as it is a solution to PDEs of order $2\beta^{-1}$ as in \cite{AN,Nanesd}.  The next two lemmas give us estimates on the random part.

 \begin{lem}[Spatial differences] \lbl{sptdiff}
Assume that \eqref{cnd} holds and that $M_{\beta,q}(t)$ and $\tilde{M}_{\beta,q}(t)$ are bounded on any time interval\footnote{This is the aforementioned temporary  moment condition.  It is assumed here (in \lemref{sptdiff} and \lemref{tmpdiff} below) only to simplify the presentation and to get to the proof of \thmref{lip} as quickly as possible in \secref{pf1stmain}.  In \secref{regandtight2}, this moment condition is shown to automatically hold under \eqref{cnd}.} $\T=[0,T]$.  There exists a constant $C$ depending only on $q\ge1$, $\max_x |u_0(x)|$, $\beta=1/\nu$, $\nu\in\lbr2^k;k\in\N\rbr$, the spatial dimension
$d=1,2, 3$, $\alpha_{d}$, and $T$ such that
\beq\lbl{sptdiffinq}
\bc
\mathbb{E}\left|\uttxR - \tilde{U}_{\beta,R}^y(t)\right|^{2q}
\le C |x-y|^{2q\alpha_{d}}&,\\
\mathbb{E}\left|U_{\beta,R}(t,x) -U_{\beta,R}(t,y)\right|^{2q}
\le C |x-y|^{2q\alpha_{d}}&,
\ec
\eeq
for all $x,y \in \Xd$, $t\in\T$, and $d=1,2,3;$ where $\lbr\alpha_{d}\rbr_{d=1}^{3}$  are as in \lemref{3rdQinequality}.   I.e., in $d=1$, we may take $\alpha_{1}=1;$ in $d=2$ we may take any fixed $\alpha_{2}\in(0,1);$ and in $d=3$, $\alpha_{3}$ may be taken to be any fixed value in $(0,\frac12)$.
\end{lem}
\begin{pf}
We prove the lattice SIE statement in \eqref{sptdiffinq} for $\ut$; the proof of the statement for $U_{\beta}$ follows the exact same steps, with obvious modifications and will  not be repeated.  Using Burkholder inequality, we have for any $(t,x,y)\in\T\times\X^{2d}$
\begin{equation}
\begin{split}
\E\lab \tilde{U}_{\beta,R}^x(t) - \tilde{U}_{\beta,R}^y(t)\rab^{2q}
\leq C\E{\lab\sum_{z\in\Xd}\int_0^t
{\lbk\Ktsxz- \Ktsyz\rbk}^2a^2
(\tilde{U}_{\beta}^{z}(s)) \df{ds}{\dd}\rab}^q
\end{split}
\label{spdiff1}
\end{equation}
For any fixed but
arbitrary point $(t,x,y)\in\T\times\X^{2d}$ let $\mu_t^{x,y}$ be the measure defined on
$[0,t]\times\Xd$ by
\beqs
\bsp
d\mu_t^{x,y}(s,z)&=\lbk\Ktsxz-\Ktsyz\rbk^2\df{ds}{\dd},
\end{split}
\eeqs
and let
$|\mu_t^{x,y}| = \mu_t^{x,y}([0,t] \times \Xd)$.    We see from \eqref{spdiff1}, Jensen's inequality
applied to the probability measure $\mu_t^{x,y}/\lab\mu_t^{x,y}\rab$,
the growth condition on $a$, the definition of $\tilde{M}_{\beta,q}(t)$, and elementary inequalities, that we have
\begin{equation}
\begin{split}
\mathbb{E}\left| \tilde{U}_{\beta.R}^x(t) - \tilde{U}_{\beta,R}^y(t)\right|^{2q}
&\le C\mathbb{E}\Big[\int_{[0,t] \times \Xd}
\left|a(\tilde{U}_{\beta}^{z}(s))\right|^{2q}\frac{d\mu_t^{x,y}(s,z)}{|\mu_t^{x,y}|}\Big]
{|\mu_t^{x,y}|}^q \\
&\le  C\Big[\int_{[0,t] \times \Xd}\left (1+\tilde{M}_{\beta,q}(s)\right)
\frac{d\mu_t^{x,y}(s,z)}{|\mu_t^{x,y}|}\Big] {|\mu_t^{x,y}|}^q
\end{split}
\label{spdiff2}
\end{equation}
Now, using the boundedness assumption on $\tilde{M}_{\beta,q}$ on $\mathbb{T}$ for $d=1,2,3$,
we get
 \begin{equation*}
\begin{split}
\mathbb{E}\left| \tilde{U}_{\beta,R}^x(t) - \tilde{U}_{\beta,R}^y(t)\right|^{2q}&\le C\lab\mu_t^{x,y}\rab^q
\le \lbk C_{d}t^{p_{d}(\alpha_{d})}\rbk^{q} |x-y|^{2q\alpha_{d}}
\\&\le \tilde{C}_{d}|x-y|^{2q\alpha_{d}};\ \alpha_{d}\in I_{d},
\end{split}
\label{spdiff3}
\end{equation*}
where the last inequality follows from \lemref{3rdQinequality} and \eqref{sptlkernel2} in \remref{spatlkerrem}, and where the constant $\tilde{C}<\infty$ is as in \remref{spatlkerrem}.  
\end{pf}

\begin{lem}[Temporal differences]\label{tmpdiff}
Assume that \eqref{cnd} holds and that ${M}_{\beta,q}(t)$ and $\tilde{M}_{\beta,q}(t)$ are bounded on any time interval $\T=[0,T]$.
 There exists a constant
$C$ depending only on $q\ge1$, $\max_x|u_0(x)|$, $\beta=1/\nu$, $\nu\in\lbr2^k;k\in\N\rbr$, the spatial dimension $d=1,2, 3$, and $T$ such that
\beq\lbl{tmpdiffinq}
\bc
\E\left|{U}_{\beta,R}(t,x) - {U}_{\beta,R}(r,x) \right|^{2q}
\leq C\lab t-r\rab^{\tf{(2\nu-d)q}{2{\nu}}}; & x \in \Rd, t,r \in\mathbb{T},\\
\mathbb{E}\left| \tilde{U}_{\beta,R}^x(t) - \tilde{U}_{\beta,R}^x(r) \right|^{2q}
\leq C\lab t-r\rab^{\tf{(2\nu-d)q}{2{\nu}}}; & x \in \Xd, t,r \in\mathbb{T},\\
\ec
\eeq
for $d=1,2, 3$.
\end{lem}

\begin{pf}  We prove the lattice SIE statement in \eqref{tmpdiffinq} for $\ut$; the proof of the statement for $U_{\beta}$ follows the exact same steps, with obvious modifications. 
Assume without loss of generality  that $r<t$.
Using Burkholder inequality, and using the change of variable $\rho=t-s$,
we have for $(r,t,x)\in\T^2\times\X^{d}$
\begin{equation}
\begin{split}
\E\lab  \tilde{U}_{\beta,R}^x(t) - \tilde{U}_{\beta,R}^x(r) \rab^{2q}
&\leq C\E{\lab\sum_{z\in\Xd}\int_0^r{\lbk\Ktsxz- \Krsxz\rbk}^2a^2(\tilde{U}^{z}(s)) \df{ds}{\dd}\rab}^q
\\&+C\E{\lab\sum_{z\in\Xd}\int_0^{t-r}
{\lbk\Kroxz\rbk}^2a^2
(\tilde{U}^{z}(t-\rho)) \df{d\rho}{\dd}\rab}^q
\end{split}
\label{tmpdiff1}
\end{equation}
For a fixed point $(r,t,x)$ and a fixed $\beta$,
let $\mu_{\beta,t,r}^x$ be the measure defined on $[0,r]\times\Xd$ by
\begin{equation*}
\bsp
d\mu_{\beta,t,r}^x(s,z)&=\lbk\Ktsxz-\Krsxz\rbk^2\df{ds}{\dd}
\end{split}
\end{equation*}
and let $|\mu_{\beta,t,r}^x| = \mu_{\beta,t,r}^x([0,r] \times \Xd)$.
Also, for a fixed $x\in\X^d$ and $\beta$, let $\kappa^x$ be the measure defined
on $[0,t-r]\times\Xd$ by
\begin{equation*}
\bsp
d\kappa_{\beta}^x(\rho)&=\lbk\Kroxz\rbk^2\df{d\rho}{\dd}
\end{split}
\end{equation*}
and let $|\kappa_{\beta}^x| = \kappa_{\beta}^x([0,t-r] \times \Xd)$.
Then, arguing as in Lemma \ref{sptdiff}
above we get that
 \begin{equation*}
 \begin{split}
\mathbb{E}\left| \tilde{U}_{\beta,R}^x(t) - \tilde{U}_{\beta,R}^x(r) \right|^{2q}&\le
C\lpa\lab\mu_{\beta,t,r}^x\rab^q+\lab\kappa_{\beta}^x\rab^q\rpa\le C(t-r)^{\tf{(2\nu-d)q}{2{\nu}}},
\end{split}
\end{equation*}
for $d=1,2,3$, where the last inequality follows from \lemref{2ndQinequality} and
\lemref{4thQinequality}, completing the proof.
\end{pf}

\section{Proof of the first main Theorem}\lbl{pf1stmain}
Here, we prove \thmref{lip}.  We start first by recalling a useful elementary Gronwall-type lemma whose proof can be found in Walsh \cite{W}.
\blm\lbl{grnwlrate}
Let $\lbr g_n(t)\rbr_{n=0}^\infty$ be a sequence of positive functions such that $g_0$
is bounded on $\T=[0,T]$ and
$$g_n(t)\le C \intzt g_{n-1}(s) (t-s)^\alpha ds, \quad n=1,2,\ldots$$
for some constants $C>0$ and $\alpha>-1$.  Then, there exists a (possibly different)
constant $C>0$ and an integer $k>1$ such that for each $n\ge1$ and $t\in\T$
$$g_{n+mk}(t)\le C^m\intzt g_n(s)\df{t-s}{(m-1)!} ds;\quad m=1,2,\ldots.$$
\elm
We are now ready for our proof.
\bpfs{Proof of \thmref{mainthm1}}  For the existence proof, we construct a solution iteratively.  So, given a space-time white noise $\sW$, on some $\OFFtP$,
define
\beq\lbl{itdef}
\bc
\zivtx=\ds\intrd\KIBtxy \uny dy&\cr
\ds\npoivtx=\zivtx\ds+\intrdzt\KIBtsxy a(\nivsy)\sW(ds\times dy)
\ec
\eeq
We will show that, for any $p\ge2$ and all $d=1,2,3$, the sequence $\lbr \nivtx\rbr_{n\ge1}$
converges in $L^p(\Omega)$ to a solution.  Let
$$D_{\beta,n,p}(t,x):=\E\lab\npoivtx-\nivtx\rab^p$$

$$D^*_{\beta,n,p}(t):=\sup_{x\in\Rd}\dnptx.$$
Starting with the case  $p>2$, we bound $D_{\beta,n,p}$ using Burkholder inequality,
the Lipschitz condition $(a)$ in \eqref{lcnd}, and then H\"older inequality with $0\le\epsilon\le1$
and $q=p/(p-2)$ to get
\beqs
\bsp
&\dnptx=\E\lab\intrdzt\KIBtsxy\lbk a(\nivsy)-a(\nmoivsy)\rbk\sW(ds\times dy)\rab^{p}
\\&\le C\E\lab\intrdzt\lpa\KIBtsxy\rpa^2\lbk \nivsy-\nmoivsy\rbk^2ds dy\rab^{p/2}
\\&\le C\lpa\intrdzt\lbk\KIBtsxy\rbk^{2\epsilon q}dsdy\rpa^{p/2q}
\\&\times\intrdzt\lpa\KIBtsxy\rpa^{(1-\epsilon)p}\dnmopsy ds dy
\end{split}
\eeqs
Take $\epsilon=(p-2)/p$ in the above ($2\epsilon q=(1-\epsilon)p=2$), take the supremum over the space variables,
and use \lemref{2ndQinequality} to see that, for $d=1,2,3$ the above reduces to
\beq\lbl{mdiff}
\bsp
\mnpt\le C \lpa t^{\tf{2\nu-d}{2\nu}}\rpa^{\tf{p-2}{2}}\intzt\mnmops\lbk t-s\rbk^{\tf{-d}{2\nu}} ds
\end{split}
\eeq
The case $p=2$ is simpler.  We apply Burkholder's inequality to $D_{n,2}$ and then take the space supremum
to get
\beq\lbl{mdiff1}
\bsp
D^*_{\beta,n,2}(t)\le C \intzt D^*_{\beta,n-1,2}(s)\lbk t-s\rbk^{\tf{-d}{2\nu}} ds
\end{split}
\eeq
I.e., on any time interval $\T=[0,T]$, the integral multiplier on the r.h.s.~of \eqref{mdiff} is bounded; and
if $D^*_{\beta,n-1,p}$ is bounded on $\T$ then so is $D^*_{\beta,n,p}$, for every $p\ge2$.  Now,
$$D^*_{\beta,0,p}(t)\le C\sup_{x\in\Rd}\E\lab\intrdzt\lbk\KIBtsxy\rbk^{2}a^2\lpa\zivsy\rpa dsdy\rab^{\tf p 2}$$
Since $\un$ is bounded and deterministic, then so are  $U^{(0)}$ and $a(U^{(0)})$.  The latter assertion follows
from the growth condition on $a$ in \eqref{lcnd}.  Thus, by  \lemref{2ndQinequality} $D^*_{\beta,0,p}$ is bounded on $\T$  for $d=1,2,3$ and so are all the
$D^*_{\beta,n,p}$.  \lemref{grnwlrate} now implies that for each $d=1,2,3$, the series
$\sum_{m=0}^\infty \lbk D^*_{\beta,n+mk,p}(t)\rbk^{1/p}$ converges uniformly on compacts for
each $n$, which in turn implies that $\sum_{n=0}^\infty\lbk D^*_{\beta,n,p}(t)\rbk^{1/p}$ converges uniformly
on compacts.  Thus $U_{\beta}^{(n)}$ converges in $L^p(\Omega)$ for $p\ge2$, uniformly on $\T\times\Rd$ for $d=1,2,3$.  Let
$U_{\beta}(t,x):=\lim_{n\to\infty}\nivtx$.
It is easy to see that $U_{\beta}$  satisfies \eqref{isltbmsie},
and hence solves the $\beta$-ISLTBM SIE $\isltsie$.  This follows from \eqref{itdef} since
the Lipschitz condition in \eqref{lcnd} gives
\beqs
\bsp
\E\lab a(U_{\beta}(t,x))-a(\nivtx)\rab^2\le C\E\lab U_{\beta}(t,x)-\nivtx\rab^2\to0\quad\mbox{as }n\to\infty
\end{split}
\eeqs
uniformly on $\T\times\Rd$.  Therefore, the stochastic integral term in \eqref{itdef} converges to
the same term with $U_{\beta}^{(n)}$ replaced with the limiting $U_{\beta}$---i.e., it converges to
the corresponding term in $\isltsie$---as $n\to\infty$, for
\beqs
\bsp
&\E\lbk\intrdzt\KIBtsxy\lpa a(U_{\beta}(s,y))- a(\nivsy)\rpa\sW(ds\times dy)\rbk^2
\\&\le  C \intrdzt\lbk\KIBtsxy\rbk^2 \E\lbk U_{\beta}(s,y)-\nivsy\rbk^2 ds dy
\longrightarrow 0
\end{split}
\eeqs
as $n\to\infty$.  It follows that $U_{\beta}$ satisfies the $\beta$-ISLTBM SIE $\isltsie$. Also, the solution
is strong since the $U_{\beta}^{(n)}$ are constructed for a given white noise $\sW$, and the
limit $U_{\beta}$ satisfies \eqref{btpsol} with respect to that same $\sW$.
Clearly $U_{\beta}$ is $L^p(\Omega)$ bounded on $\T\times\Rd$, $d=1,2,3$, for any $p\ge2$ and for any $T>0$.

To show uniqueness fix an arbitrary $\beta^{-1}\in\lbr2^k;k\in\N\rbr$---and suppress the dependence of solutions on $\beta$---and let $d=1,2,3$, let $T>0$ be fixed but arbitrary, and let $U_1$ and $U_2$ be two solutions to our $\beta$-ISLTBM SIE \eqref{isltbmsie} that are $L^2(\Omega)$-bounded on $\T\times\Rd$.
Fix an arbitrary $(t,x)\in\Rp\times\Rd$.
Let $D(t,x)=U_2(t,x)-U_1(t,x)$, $L_2(t,x)=\E D^2(t,x)$,
and $L^*_2(t)=\sup_{x\in\Rd}L_2(t,x)$ (which is bounded on $\T$ by hypothesis).
Then, using \eqref{isltbmsie}, the Lipschitz condition in \eqref{lcnd},
and taking the supremum over the space variable and using \lemref{2ndQinequality}
we have
\beq\lbl{un1}
\bsp
L_2(t,x)&=\intrdzt \E\lbk a(U_2(s,y))- a(U_1(s,y))\rbk^2\lbk\KIBtsxy\rbk^2dsdy
\\&\le C\intrdzt L_2(s,y)\lbk\KIBtsxy\rbk^2dsdy
\\&\le C\intzt L^*_2(s)\intrd\lbk\KIBtsxy\rbk^2dyds
\le C\intzt\df{L^*_2(s)}{(t-s)^{\tf{d}{2\nu}}}ds
\end{split}
\eeq
Iterating and interchanging the order of integration we get
 \beq\lbl{un2}
\bsp
L_2(t,x)&\leq  C\lbr \int_0^t{L^*_2(r)}\lpa\int_r^t\frac{ds}{{(t - s)^{\tf{d}{2\nu}}}{(s-r)^{\tf{d}{2\nu}}}}\rpa dr\rbr\\
 &\le C\lpa\ds\int_0^t L^*_2(s)ds \rpa
\end{split}
\eeq
for any $d=1,2,3$.  Hence,
\beq\lbl{un3}
\bsp
L^*_2(t)\le C\lpa\ds\int_0^t L^*_2(s)ds \rpa
\end{split}
\eeq
for every $t\ge0$.  An easy application of Gronwall's lemma gives that $L^*_2\equiv0$.  So for every $(t,x)\in\Rp\times\Rd$ and $d=1,2,3$ we have $U_1(t,x)=U_2(t,x)$ with probability one.   The indistinguishability of $U_{1}$ from $U_{2}$, and hence pathwise uniqueness, follows immediately from their H\"older regularity, which we now turn to.

For any given $\beta^{-1}=\nu\in\lbr2^k;k\in\N\rbr$, we have just shown that, under the Lipschitz conditions \eqref{lcnd}, our $\beta$-ISLTBM SIE in \eqref{isltbmsie} has an $L^{p}(\Omega)$-bounded solution $U_{\beta}(t,x)$ on $\T\times\Rd$ for any $T>0$ and any $p\ge2$.  Equivalently, ${M}_{\beta,q}(t)= \sup_x \mathbb{E}|U_{\beta}(t,x)|^{2q}$, $q\ge1$, is bounded on any time interval $\T$.   Recalling that the deterministic part\footnote{Of course, the deterministic part of $\isltsie$ is, as discussed before, the integral $\intrd\KIBtxy \uny dy$; and the random part is  $\intrdzt\KIBtsxy a(\usy)\sW(ds\times dy)$.} of $U_{\beta}$ is a $\mathrm{C}^{1,2\nu}(\Rp,\Rd)$ function, we can then use  \lemref{sptdiff} and \lemref{tmpdiff} above, on the random part of $U_{\beta}$ for $d=1,2,3$ to straightforwardly get the desired local H\"older regularity for the direct solution of $\isltsie$, $U_{\beta}$, as follows: we let $q_n=n+d$ for $n\in\{0,1,\ldots\}$ and
let $n=m+d$ for $m=\{0,1,\ldots\}$, we then have from \lemref{sptdiff} and \lemref{tmpdiff}
that
\beq\lbl{directbtpsiehldr}
\bc
\E\lab\utx-\uty\rab^{2n+2d}\le C_{d}\lab x-y\rab^{(2n+2d)\alpha_{d}},\\
\E\lab\utx-\urx\rab^{2m+4d}\le C\lab t-r\rab^{\tf{(2\nu-d)(m+2d)}{2{\nu}}}.
\ec
\eeq
for $d=1,2,3$.  Thus as in Theorem 2.8 p.~53 and Problem 2.9 p.~55 in \cite{KS} we get that the spatial H\"older exponent is
$\gamma_s\in\lpa 0,\tf{2(n+d)\alpha_{d}-d}{2n+2d}\rpa$ and the temporal exponent is $\gamma_t\in\lpa0,\tf{m\lpa1-d/2{\nu}\rpa+d(1-d/{\nu})}{2m+4d}\rpa$
$\forall m,n$.  Taking the limits as $m,n\to\infty$, we get
$\gamma_t\in\lpa0,\tf{2\nu-d}{4{\nu}}\rpa$ and $\gamma_s\in\lpa0,\alpha_{d}\rpa$, for $d=1,2,3$.
The proof is complete.      
\epfs
\section{Proof of the second main Theorem}
\subsection{Regularity and tightness without the Lipschitz condition}\lbl{regandtight2}
As we mentioned in \secref{regandtight}, the finiteness assumption of $M_{\beta,q}(t)$ and $\tilde{M}_{\beta,q}(t)$ on $\T$ in \lemref{sptdiff} and \lemref{tmpdiff} is for convenience only.  We now proceed to show how to remove that assumption by showing it automatically holds under the weaker conditions \eqref{cnd}. It is easily seen that if $a$ is bounded then,  for all spatial dimensions $d=1,2,3$, $\tilde{M}_{\beta,q}$ is bounded on any compact time interval $\T=[0,T]$ (see \remref{bddrem} below).  The following Proposition gives an exponential upper bound on the growth of $\tilde{M}_{\beta,q}$ in time in all $d=1,2,3$ under the conditions in \eqref{cnd}.  The same result holds for $M_{\beta,q}$ with only notational and obvious changes to the following proofs.
\begin{prop}[Exponential bound for $\tilde{M}_{\beta,q}$] \label{Expbd}
Assume that $\uttx$ is a solution of the $\beta$-ISLTRW SIE $\isltrwsien$ on $\T\times\Xd$, and assume that the conditions in \eqref{cnd} are satisfied.  There exists a constant $C$ depending only on $q$, $\max_x|u_0(x)|$, the dimension $d$, $\beta$, and $T$ such that
\begin{equation*}
\tilde{M}_{\beta,q}(t) \leq  C\lpa1+ \ds\int_0^t \tilde{M}_{\beta,q}(s)ds \rpa;  0\le t\le T, 
\end{equation*}
for every $q\ge1, \beta\in\lbr\tf1{2^k};k\in\N\rbr\mbox{and }d=1,2,3$.   Hence, $\tilde{M}_{\beta,q}(t) \leq C\exp{\{Ct\}}\mbox{ for } 0\le t\le T, q\ge1, \beta\in\lbr\tf1{2^k};k\in\N\rbr, \mbox{and }d=1,2,3$.
   In particular, $\tilde{M}_{\beta,q}$ is bounded on $\mathbb{T}$ for all $q\ge1$, $\beta\in\lbr\tf1{2^k};k\in\N\rbr$, and $d=1,2,3$.
 \end{prop}

 The proof of Proposition \ref{Expbd} proceeds via the following lemma and its
 corollary.

\begin{lem} \label{1stboundonUtilde}Under the same assumptions as in \propref{Expbd}
there exists a constant $C$ depending only on $q$, $\max_x |u_0(x)|$, the dimension $d$, $\beta$, and $T$ such that
\begin{equation*}
\tilde{M}_{\beta,q}(t) \leq\bc C \lpa1 + \ds{\int_0^t\df{\tilde{M}_{\beta,q}(s)}{(t-s)^{\tf{d}{2\nu}}} ds}\rpa;
   &0< t\le T,\\
  C;&t=0,
  \ec
\end{equation*}
for every $q\ge1$, $\beta\in\lbr\tf1{2^k};k\in\N\rbr$, and $d=1,2,3$.
\end{lem}

\begin{pf}
Fix $q\ge1$, let $\ds\tilde{U}_{\beta,D}^{x}(t)\overset{\triangle}{=} \sum_{y\in\Xd} \Ktxy u_0(y)$
(the deterministic part of $\tilde{U}_{\beta}$).  Then,  for any $(t,x)\in\T\times\Xd$,
we apply Burkholder inequality to the random term  $\tilde{U}_{\beta,R}^x(t)$ to get
\begin{equation}\label{afterBurkholder}
\begin{split}
\E\lab\tilde{U}_{\beta}^x(t)\rab^{2q}&= \E\lab \sum_{y\in {\Xd}}\int_0^t \Ktsxy
\frac{a(\tilde{U}_{\beta}^y(s))}{\sdd} dW^{y}(s)+ \tilde{U}_{\beta,D}^{x}(t) \rab^{2q} \\
\\&\leq C\lpa\E\lab\sum_{y\in {\Xd}} \int_0^t {\lpa\Ktsxy\rpa}^2\frac{a^2(\tilde{U}_{\beta}^y(s))}{\dd}ds\rab^q + \lab\tilde{U}_{\beta,D}^{x}(t)\rab^{2q}\rpa.
\end{split}
\end{equation}
Now, for a fixed point $(t,x)\in\T\times\Xd$ let $\mu_t^x$ be the measure on $[0,t] \times\Xd$ defined by 
$d\mu_t^x(s,y)=\lbk{\lpa\Ktsxy\rpa}^2/\dd\rbk ds$, and let $|\mu_t^x| = \mu_t^x([0,t] \times \Xd)$.  Then, we can rewrite \eqref{afterBurkholder} as
 \begin{equation}\label{afterBurkholder2}
\begin{split}
 \E\lab\tilde{U}_{\beta}^x(t)\rab^{2q}
\leq C\lpa\E\left|\int_{[0,t] \times \Xd}a^2(\tilde{U}_{\beta}^y(s))
\frac{d\mu_t^x(s,y)}{|\mu_t^x|}\right|^q {|\mu_t^x|}^q +|\tilde{U}_{\beta,D}^{x}(t)|^{2q}
\rpa.
\end{split}
\end{equation}
Observing that $\mu_t^x/|\mu_t^x|$ is a probability measure, we apply Jensen's inequality,
the growth condition on $a$ in \eqref{cnd}, and other elementary inequalities
to \eqref{afterBurkholder2} to obtain
\begin{align*}
&\E\lab\tilde{U}_{\beta}^x(t)\rab^{2q}
\leq C\lpa\E\lbk\int_{[0,t] \times \Xd}\lab a(\tilde{U}_{\beta}^y(s))
\rab^{2q}\frac{d\mu_t^x(s,y)}{|\mu_t^x|}\rbk
{|\mu_t^x|}^q + \lab\tilde{U}_{\beta,D}^{x}(t)\rab^{2q}\rpa\\
& \leq C\lbk\int_{[0,t] \times \Xd}
\lpa1 + \E\lab\tilde{U}_{\beta}^y(s)\rab^{2q}\rpa{d\mu_t^x(s,y)}\rbk
{|\mu_t^x|}^{q-1} + C \lab\tilde{U}_{\beta,D}^{x}(t)\rab^{2q} \\
&= C\lpa\lbk\sum_{y\in {\Xd}}\int_0^t\frac{\lpa\Ktsxy\rpa^2}
{\dd}\lpa1 + \E\lab\tilde{U}_{\beta}^y(s)\rab^{2q}\rpa ds\rbk
{|\mu_t^x|}^{q-1}+\lab\tilde{U}_{\beta,D}^{x}(t)\rab^{2q} \rpa
\end{align*}
Using \lemref{2ndQinequality} we see that  $|\mu_t^x|$ is
uniformly bounded for $t\leq T$ and $d=1,2,3$.
So, using the boundedness of $u_0$, and hence of $\tilde{U}_{\beta,D}^{x}(t)$
by the simple fact that $\sum_{y\in\Xd}\Ktxy=1$,
Lemma \ref{2ndQinequality} and the definition of $\tilde{M}_{\beta,q}(s)$, we get
\begin{align*}
\E\lab\tilde{U}_{\beta}^x(t)\rab^{2q}
&\leq C\lpa1+\sum_{y\in {\Xd}}
\int_0^t\frac{\lpa\Ktsxy\rpa^2}{\dd}\tilde{M}_{\beta,q}(s) ds\rpa\\
&\overset{R_1}{\leq} C \lpa1 + \int_0^t\frac{\tilde{M}_{\beta,q}(s)}{(t-s)^{\tf{d}{2\nu}}} ds
\rpa.
\end{align*}
Here, $R_1$ holds for $d=1,2,3$.   This implies that
\begin{equation*}
\tilde{M}_{\beta,q}(t) \leq C \lpa1 + \int_0^t\frac{\tilde{M}_{\beta,q}(s)}{(t-s)^{\tf{d}{2\nu}}} ds\rpa.
\end{equation*}
Of course, $\tilde{M}_{\beta,q}(0)=\sup_{x}\lab\unx\rab^{2q}\le C$, by the boundedness and nonrandomness
assumptions on $\unx$ in \eqref{cnd}. The proof is complete.
\end{pf}
\brm\lbl{bddrem}
It is clear that for a bounded $a$, $\tilde{M}_{\beta,q}$ is locally bounded in time.  This follows immediately from \lemref{2ndQinequality}
along with \eqref{afterBurkholder2} above.
\erm
\begin{cor}
Under the same assumptions as those in \propref{Expbd}
there exists a constant $C$ depending only on $q$, $\max_x |u_0(x)|$, the dimension $d$, $\beta$, and $T$
such that
\begin{equation*}
\tilde{M}_{\beta,q}(t) \leq C\lpa1+\ds\int_0^t \tilde{M}_{\beta,q}(s)ds\rpa, 0\le t\le T, q\ge1, \beta\in\lbr\tf1{2^k};k\in\N\rbr\mbox{and }d=1,2,3;
\end{equation*}
and hence
\begin{equation*}
\tilde{M}_{\beta,q}(t) \leq C\exp{\{Ct\}}; \quad \forall\,0\le t\le T, \ q\ge1, \beta\in\lbr\tf1{2^k};k\in\N\rbr,\mbox{ and }d=1,2,3.
\end{equation*}
\label{befG}
\end{cor}

\begin{pf}
Iterating the bound in Lemma \ref{1stboundonUtilde} once, and changing the
order of integration, we obtain
\begin{equation}
\begin{split}
 &\tilde{M}_{\beta,q}(t)
 \\ &\leq  C\lbr 1 + C\lbk\int_0^t\frac{ds}{(t-s)^{\tf{d}{2\nu}}}+
 \int_0^t{\tilde{M}_{\beta,q}(r)}\lpa\int_r^t\frac{ds}{{(t-s)^{\tf{d}{2\nu}}}{(s-r)^{\tf{d}{2\nu}}}}
 \rpa dr\rbk\rbr\\
 &\le C\lpa1+ \ds\int_0^t \tilde{M}_{\beta,q}(s)ds \rpa
 \end{split} \label{mom1}
\end{equation}
for $d=1,2,3$.  The proof of the last statement is a straightforward application of Gronwall's
lemma to \eqref{mom1}.  This finishes the proof of \coref{befG} and thus of \propref{Expbd}.
\end{pf}

The regularity, tightness, and weak limit conclusions for the $\beta$-ISLTRW SIEs now follow.

\blm[Regularity and tightness]\lbl{regtight}
Assume that the conditions \eqref{cnd} hold, and that
$\lbr\uttxn\rbr_{n\in\Ns}$ is a sequence of spatially-linearly-interpolated solutions to
the $\beta$-ISLTRW SIEs $\lbr\isltrwsien\rbr_{n\in\Ns}$ in \eqref{btrwsie}.   Then
\begin{enumerate}\renewcommand{\labelenumi}{$($\alph{enumi}$)$}
\item For every $n$, $\uttxn$ is continuous on $\Rp\times\Rd$.  Moreover, with probability one,
the continuous map $(t,x)\mapsto\uttxn$ is locally
$\gamma_t$-H\"older continuous in time with $\gamma_t\in\lpa0,\tf{2\nu-d}{4\nu}\rpa$ for $d=1,2,3$.
\item  There is a $\beta$-ISLTRW SIE weak limit solution to $\isltsie$,  call it $U_{\beta}$, such that $U_{\beta}(t,x)$ is $L^p(\Omega,\P)$-bounded on $\T\times\Rd$ for every $p\ge2$ and $U_{\beta}\in\mathrm{H}^{\lpa{\tf{2\nu-d}{4\nu}}\rpa^{-},\alpha^{-}_{d}}(\T\times\Rd;\R)$ for every $d=1,2,3$ and $\alpha_{d}\in I_{d}$, where $\alpha_{d}$ and $I_{d}$ are as in \lemref{3rdQinequality}.
\end{enumerate}
\elm
\brm
Of course in part (a) above, even without linear interpolation in space, $\uttx$ is locally H\"older continuous in time with
H\"older exponent $\gamma_t\in\lpa0,\tf{2\nu-d}{4\nu}\rpa$ for $d=1,2,3$.
\erm
\bpf  For each $n$, let $\uttxn=\uttxnD+\uttxnR$ be the decomposition of $\uttxn$ in \eqref{btrwsie}
into its deterministic and random parts, respectively.
\begin{enumerate}\renewcommand{\labelenumi}{$($\alph{enumi}$)$}
\item By \lemref{fodde}, $\uttxnD$ is clearly smooth in time;
so it is enough to consider the random term $\uttxnR$.
We let $q_m=m+2$ for $m\in\{0,1,\ldots\}$,
we then have from \lemref{tmpdiff} that
\beq\lbl{siedishldrt}
\E\lab\tilde{U}_{\beta,R}^x(t) - \tilde{U}_{\beta,R}^x(r) \rab^{4+2m}\le C\lab t-r\rab^{\tf{(2\nu-d)(m+2)}{2\nu}}.
\eeq
for $d=1,2,3$.  Thus as in Theorem 2.8 p.~53 \cite{KS} we get that
$\gamma_t\in\lpa0,\tf{m\lpa1-d/2\nu\rpa+2-d-d/{\nu}}{2m+4}\rpa$ for every $m$.  Taking the limit as $m\to\infty$, we get
$\gamma_t\in\lpa0,\tf{2\nu-d}{4\nu}\rpa$ for $d=1,2,3$.
\item   By \lemref{btrwasbtp} it follows that $\uttxnD$ converges pointwise
to the deterministic part of $\isltsie$ in \eqref{btpsol}; i.e.,
\beq\lbl{det2det}
\lim_{n\to\infty}\uttxnD=\intrd\KIBtxy \uny dy.
\eeq
We also conclude from \lemref{sptdiff} and \lemref{tmpdiff} that the sequence $\lbr\uttxnR\rbr_{n\in\Ns}$
is tight on $\mathrm{C}(\T\times\Rd)$  for $d=1,2,3$.  Thus there exists a weakly convergent
subsequence $\lbr\utnk\rbr_{k\in\N}$ and hence a $\beta$-ISLTRW SIE weak limit solution $U$ to $\isltsie$.
Then, following Skorokhod, we construct processes\footnote{As usual, $\eqd$ denotes equal in law or distribution.} $\yk\eqd\utnk$ on some filtered probability space $\Skspace$ such that with probability $1$, as $k\to\infty$, $\ytxk$ converges to a random field $\ytx$ uniformly on compact subsets of
$\T\times\Rd$ for $d=1,2,3$.  Now, for the $\beta$-ISLTRW SIEs limit regularity assertions, clearly the deterministic term on the right hand side of
\eqref{det2det} is $\mathrm{C}^{1,2\nu}$ and bounded as in \cite{AN}, so we use \propref{Expbd}, \lemref{sptdiff}, and \lemref{tmpdiff}
to obtain the regularity results for the random part.  We provide the steps here for completeness.
First, $\yk\eqd\utnk$ and so \propref{Expbd} gives us, for each $p\ge2$:
\beq\lbl{reg111}
\E\lab\ytxk\rab^{p}=\E\lab\uttxnk\rab^{p}\le C<\infty;\forall(t,x,k)\in\T\times\Rd\times\N,
d=1,2,3,
\eeq
for some constant $C$ that is independent of $k, t, x$ but that depends on the dimension $d$.
It follows that, for each $(t, x)\in\T\times\Rd$ the sequence $\lbr|Y_k(t, x)|^p\rbr_k$ is uniformly
integrable for each $p\ge2$ and each $d=1,2,3$. Thus,
\beq\lbl{spdemmntsbd}
\E\lab U_{\beta}(t,x)\rab^{p}=\E\lab\ytx\rab^{p}=\lim_{k\to\infty}\E\lab\ytxk\rab^{p}\le C<\infty;\forall(t,x)\in\T\times\Rd,
\eeq
for all $d=1,2,3$ and $p\ge2$.   Equation \eqref{spdemmntsbd} establishes the $L^p$ boundedness assertion.  In addition, for $q\ge1$ and $d=1,2,3$ we have by \propref{Expbd}
\beq\lbl{ui2}
\bsp
&\E \lab Y_{\beta,k}(t, x) - Y_{\beta,k}(t, y)\rab^{2q} + \E\lab Y_{\beta,k}(t, x) - Y_{\beta,k}(r, x)\rab^{2q}
\\&\le C \lbk\E\lab Y_{\beta,k}(t, x)\rab^{2q} + \E\lab Y_{\beta,k}(t, y)\rab^{2q} + \E\lab Y_{\beta,k}(r, x)\rab^{2q}\rbk
\\&\le C; \ \forall(k, r, t, x, y)\in\N\times\T^2\times \R^2.
\end{split}
\eeq
So, for each $(r, t, x, y)\in\T^2\times\R^2$, the sequences $\lbr\lab Y_{\beta,k}(t, x) - Y_{\beta,k}(t, y)\rab^{2q}\rbr_k$ and
$\lbr\lab Y_{\beta,k}(t, x) - Y_{\beta,k}(r, x)\rab^{2q}\rbr_k$
are uniformly integrable, for each $q\ge1$. Therefore, using
\lemref{sptdiff} and \lemref{tmpdiff}, we obtain
\beq\lbl{spdehldr1}
\bc
\E\lab\utx-\uty\rab^{2q}=\E\lab\ytx-\yty\rab^{2q}&\cr
=\ds\lim_{k\to\infty}\E\lab\ytxk-\ytyk\rab^{2q}\le C_{d} |x-y|^{2q\alpha_{d}};\ \alpha_{d}\in I_{d},\\
\E\lab\utx-\urx\rab^{2q}=\E\lab\ytx-\yrx\rab^{2q}&\cr
=\ds\lim_{k\to\infty}\E\lab\ytxk-\yrxk\rab^{2q}\le C\lab t-r\rab^{\tf{(2\nu-d)q}{2\nu}},
\ec
\eeq
for $d=1,2,\ldots,3$.  The local H\"older regularity is then obtained using exactly the same steps as in \eqref{directbtpsiehldr} and the following conclusions.
\bcm
, we let $q_n=n+d$ for $n\in\{0,1,\ldots\}$ and
let $n=m+d$ for $m=\{0,1,\ldots\}$, we then have from \eqref{spdehldr1}
that
\beq\lbl{spdehldr2}
\bc
\E\lab\utx-\uty\rab^{2n+2d}\le C_{d}\lab x-y\rab^{(2n+2d)\alpha_{d}},\\
\E\lab\utx-\urx\rab^{2m+4d}\le C\lab t-r\rab^{\tf{(4-d)(m+2d)}{4}}.
\ec
\eeq
for $d=1,2,3$.  Thus as in Theorem 2.8 p.~53 and Problem 2.9 p.~55 in \cite{KS} we get that the spatial H\"older exponent is
$\gamma_s\in\lpa 0,\tf{2(n+d)\alpha_{d}-d}{2n+2d}\rpa$ and the temporal exponent is $\gamma_t\in\lpa0,\tf{m\lpa1-d/4\rpa+d(2-d/2)-1}{2m+4d}\rpa$
$\forall m,n$.  Taking the limits as $m,n\to\infty$, we get
$\gamma_t\in\lpa0,\tf{4-d}{8}\rpa$ and $\gamma_s\in\lpa0,\alpha_{d}\rpa$, for $d=1,2,3$.
\fi\end{enumerate}
The proof is complete 
\epf
\subsection{Recalling the K-martingale approach}\lbl{Kmartsec}
For the article to be self-contained, we now recall and briefly discuss the K-martingale approach from \cite{Abtbmsie}---adapting it to this paper's setting\footnote{All we need to adapt it here is a notational change, replacing the BTRW transition density in \cite{Abtbmsie} with the $\beta$-ISLTRW one. }.  This approach is tailor-made for kernel SIEs like $\isltsie$ and other mild formulations for many SPDEs on the lattice.  The first step is to truncate to a finite lattice model as in \eqref{trnctdsie}.  Of course, even after we truncate the lattice, a remaining hurdle to applying
a martingale problem approach is that the finite sum of stochastic integrals in \eqref{trnctdsie} is not a local martingale.
So, we introduce a key ingredient in this K-martingale method: the auxiliary problem associated with the truncated $\beta$-ISLTRW SIE in \eqref{trnctdsie}, which
we now give.  Fix $(l,n)\in\N^2$ and $\tau\in\Rp$.  We define the $\tau$-auxiliary $\beta$-ISLTRW SIE associated with \eqref{trnctdsie} on $[0,\tau]\times\Xnd$ by
\beq\lbl{aux}\tag{Aux}
\xttauxnlf=\bc
\ds\uttxnD+\ds\sum_{y\in\Xnldf}\int_0^t\kappa^{x,y}_{\beta,\h,s,\tau}\lpa\xstauynlf\rpa{d\wsyn}; &x\in\Xnldf,\\
\ds\uttxnD;&x\in\Xnd\setminus\Xnldf
\ec
\eeq
where the independent BMs sequence $\lbr W_n^y\rbr_{y\in\Xnldf}$ in \eqref{aux} is the same for all $\tau>0$, as well as $x\in\Xnldf$.
We denote \eqref{aux} by $\isltrwsienlaux$.  We say that the pair of families $\lpa\lbr\xtaunlf\rbr_{\tau\ge0},\lbr W_n^y\rbr_{y\in\Xnldf}\rpa$ solves $\lbr\isltrwsienlaux\rbr_{\tau\ge0}$ on a filtered probability space $\OFFtP$ if 
there is one family of independent BMs (up to indistinguishability) $\lbr W_n^y(t);0\le t<\infty\rbr_{y\in\Xnldf}$ on $\OFFtP$ such that, for every fixed $\tau\in\Rp$
\begin{enumerate}\renewcommand{\labelenumi}{(\alph{enumi})}
\item
the process $\lbr\xttauxnlf,\sFt;0\le t\le\tau, x\in\Xnd \rbr$
has continuous sample paths in $t$ for each fixed $x\in\Xnd$ and $\xttauxnlf\in\sFt$ for all $x\in\Xnd$ for every $0\le t\le\tau$; and 
\item equation \eqref{aux} holds on $[0,\tau]\times \Xnd$, $\P$-almost surely.
\end{enumerate}
Naturally, implicit in our definition above the assumption that, for each fixed $\tau\in\Rp$, we have
$$\P\lbk\int_0^t\lpa\kappa^{x,y}_{\beta,\h,s,\tau}\lpa\xstauynlf\rpa\rpa^{2}{ds}<\infty\rbk=1;\ \forall x,y\in\Xnldf,0\le t\le\tau.$$
For simplicity, we will sometimes say that $\xtaunlf=\lbr\xttauxnlf,\sFt;0\le t\le\tau, x\in\Xnd \rbr$ is a solution to \eqref{aux} to mean the above.  Clearly, if $\xttauxnlf$ satisfies \eqref{aux} then $\uttauxnlf:=\xtautauxnlf$
satisfies \eqref{trnctdsie} at $t=\tau$ for all $x\in\Xnd$.  Also, for each $n$ and each $d=1,2,3$ 
$$\lab\kappa^{x,y}_{\beta,\h,s,\tau}\lpa\xstauynlf\rpa\rab=\lab\df{\Ktausxyn}{\sdh} a(\xstauynlf)\rab\le\df{\lab a\lpa\xstauynlf\rpa\rab}{\sdh}.$$
In addition, for each fixed $\tau\in\Rp$ and each fixed $x,y\in\Xnldf$ we have for a solution $\xtaunlf$ to \eqref{aux} that
$$\kappa^{x,y}_{\beta,\h,s,\tau}\lpa\xstauynlf\rpa\in\sFs;\quad \forall s\le\tau,$$
since, of course the deterministic ${\Ktausxyn}/{\sdh}\in\sFs$ and $a(\xstauynlf)\in\sFs$.
Thus, if $\xtaunlf$ solves \eqref{aux}; then, for each fixed $\tau>0$ and $x,y\in\Xnldf$, each stochastic integral in \eqref{aux} 
$$I_{\beta,n,l}^{\tau,x,y}=\lbr I_{\beta,n,l}^{\tau,x,y}(t):=\int_0^t\kappa^{x,y}_{\beta,\h,s,\tau}\lpa\xstauynlf\rpa{d\wsyn},\sFt;\ 0\le t\le\tau\rbr$$
is a continuous local martingale in $t$ on $[0,\tau]$.  This is clear since by a standard localization argument we may assume the boundedness of $a$ ($|a(u)|\le C$); in this case we have for each fixed $x,y\in\Xnldf$ and $\tau\in\Rp$ that
\beqs
\bsp
\E\lbk I_{\beta,n,l}^{\tau,x,y}(t)\bigg|\sFr\rbk=\int_0^r\kappa^{x,y}_{\beta,\h,s,\tau}\lpa\xstauynlf\rpa{d\wsyn}=I_{\beta,n,l}^{\tau,x,y}(r),\ r\le t\le\tau.
\end{split}
\eeqs 
So, the finite sum over $\Xnldf$ in \eqref{aux} is also  a continuous local martingale in $t$ on $[0,\tau]$.  I.e., for each $\tau>0$ and $x\in\Xnldf$
$$M_{\beta,n,l}^{\tau,x}=\lbr M_{\beta,n,l}^{\tau,x}(t):=\sum_{y\in\Xnldf}\int_0^t\kappa^{x,y}_{\beta,\h,s,\tau}\lpa\xstauynlf\rpa{d\wsyn},\sFt;\ 0\le t\le\tau\rbr \in\locm$$
with quadratic variation
\beq\lbl{qv}
\lqv M_{\beta,n,l}^{\tau,x}(\cdot)\rqv_t=\sum_{y\in\Xnldf}\int_0^t\lbk\kappa^{x,y}_{\beta,\h,s,\tau}\lpa\xstauynlf\rpa\rbk^2ds
\eeq
where we have used the independence of the BMs $\lbr W_n^y\rbr_{y\in\Xnldf}$ within the lattice $\Xnldf$.
For each $\tau>0$, we call $M_{\beta,n,l}^{x,\tau}$ a kernel local martingale (or K-local martingale).

There is another complicating factor in formulating our K-martingale problem approach that is not present in the standard SDEs setting.
To easily extract solutions to the truncated $\beta$-ISLTRW SIEs in \eqref{trnctdsie} from the family of auxiliary problems
$\lbr\isltrwsienlaux\rbr_{\tau>0}$ in \eqref{aux}, we want the independent BMs sequence
$\lbr W_n^y\rbr_{y\in\Xnldf}$ to not depend on the choices of $\tau$ and $x$.   I.e., we want all the K-local martingales in \eqref{aux}
to be stochastic integrals with respect to the same sequence $\lbr W_n^y\rbr_{y\in\Xnldf}$, regardless of $\tau$ and $x$.
With this in mind, we now formulate the K-martingale problem associated with the auxiliary $\beta$-ISLTRW SIEs in \eqref{aux}.  Let
\beq\lbl{Zspace}
\C_{n,l}:=\lbr u:\Rp\times\lpa\Xnldf\rpa^2\to\Rs; t\mapsto u^{x_1,x_2}(t) \mbox{ is continuous }\forall x_1,x_2\rbr.
\eeq
For $u\in\C_{n,l}$, let $u^{x_1,x_2}(t)=\lpa u_1^{x_1}(t),u_2^{x_2}(t)\rpa$ with $u^x(t)=u^{x,x}(t)$;
and for any $\tau_1,\tau_2>0$ and any $x_1,x_2,y\in\Xnldf$ let
\beq\lbl{ups}
\Upsilon^{x_{i,j},y}_{\h,t,\tau_{i,j}}\lpa u^y(t)\rpa:=\df{\Ktauitxiyn}{\sdh} a(\uity) \df{\Ktaujtxjyn}{\sdh}a(\ujty);\quad 1\le i,j\le2,
\eeq
(we are allowing the cases $\tau_1=\tau_2$ and/or $x_1=x_2$)
where for typesetting convenience we denoted the points $(\tau_i,\tau_j)$ and $(x_i,x_j)$ by $\tau_{i,j}$ and $x_{i,j}$, respectively.
We denote by $\pa_i$ and $\pa^2_{ij}$ the first order partial derivative with respect to the $i$-th argument and the second
order partials with respect to the $i$ and $j$ arguments, respectively.
Let $\C^2=\C^2(\Rs;\R)$ be the class of twice continuously differentiable real-valued functions on $\Rs$  and let
\beq\lbl{twicwcontbdd}
\C_b^2=\lbr f\in\C^2; \mbox{$f$ and its derivatives up to second order are bounded}\rbr.
\eeq
Now, for $\tau_1,\tau_2>0$, for $f\in\C^2_b$, and for $(t,x_1,x_2,u)\in[0,\tau_1\wedge\tau_2]\times\lpa\Xnldf\rpa^2\times \C_{n,l}$ let
\beq\lbl{siegen}
\bsp
\lpa\sAKtauot f\rpa(t,x_1,x_2,u)&:=\sum_{1\le i\le2}\pa_if\lpa u^{x_1,x_2}(t)\rpa \frac{\pa}{\pa t}\uttxinD
\\&+ \frac12\sum_{1\le i,j\le2} \pa^2_{ij}f\lpa u^{x_i,x_j}(t)\rpa\sum_{y\in\Xnldf}\Upsilon^{x_{i,j},y}_{\h,t,\tau_{i,j}}\lpa u^y(t)\rpa
\end{split}
\eeq
Let $X_{\beta,n,l}^\tau=\lbr\xttauxnlf;0\le t\le\tau,x\in\Xnd\rbr$ be a continuous in $t$ adapted real-valued process on a filtered probability space $\OFFtP$.
For every $\tau_1,\tau_2>0$ define the two-dimensional stochastic process $Z_{\beta,n,l}^{\tau_{1,2}}$:
\beq\lbl{Z}
\lbr Z_{\beta,n,l}^{x_{1,2},\tau_{1,2}}(t)=\lpa\xttauoxnlf,\xttautxnlf\rpa; (t,x_1,x_2)\in[0,\tau_1\wedge\tau_2]\times\lpa\Xnldf\rpa^2\rbr
\eeq
with $Z_{\beta,n,l}^{y,\tau_{1,2}}(t)=\lpa\xttauoynlf,\xttautynlf\rpa$ and let $U_0^{x_1,x_2}=(\unxo,\unxt)$.  We say that the family $\lbr X_{\beta,n,l}^\tau\rbr_{\tau\ge0}$ satisfies the K-martingale problem associated with the auxiliary
$\beta$-ISLTRW SIEs in \eqref{aux} on $\Rp\times\Xnd$ if for every $ f\in\C_b^2$, $0<\tau_1,\tau_2<\infty$, $\tau=\tau_1\wedge\tau_2$, $t\in[0,\tau]$, $x_1,x_2\in\Xnldf$, and
$x\in\Xnd\setminus\Xnldf$ we have
\beq\lbl{kmart}\tag{KM}
\bc
f(Z_{\beta,n,l}^{x_{1,2},\tau_{1,2}}(t))-f(U_0^{x_1,x_2})-\ds\int_0^t \lpa\sAKtauot f\rpa(s,x_1,x_2,Z_{\beta,n,l}^{\tau_{1,2}})ds\in\locm;\\
\ds\xttauxnlf=\uttxnD.
\ec
\eeq
We are now ready to state the equivalence of the K-martingale
problem in \eqref{kmart} to the auxiliary SIEs in \eqref{aux} and its implication for the $\beta$-ISLTRW SIE in \eqref{trnctdsie}.  This result is of independent interest and is stated as the following theorem\footnote{This is because it is easily adaptable to many mild formulations of SPDEs, of different orders, not just for the $\beta$-ISLTBM SIEs.  Since we don't prove uniqueness under less than Lipschitz conditions for our $\beta$-ISLTBM SIEs, we have not explicitly mentioned the uniqueness implications of our K-martingale approach.  More on that in future articles.}.
\bthm\lbl{kmarteqaux}
The existence of a solution pair $\lpa\lbr\xtaunlf\rbr_{\tau\ge0},\lbr W_n^y\rbr_{y\in\Xnldf}\rpa$
to $\lbr\isltrwsienlaux\rbr_{\tau\ge0}$ in \eqref{aux} on a filtered probability space $\OFFtP$ 
is equivalent to the existence of a family of processes $\lbr\xtaunlf\rbr_{\tau\ge0}$ satisfying \eqref{kmart}.  Furthermore, if there is $\lbr\xtaunlf\rbr_{\tau\ge 0}$ satisfying
\eqref{kmart} then there is a solution to \eqref{trnctdsie} on $\Rp\times\Xnd$.
\ethm
The proof follows the exact same steps as the proof of Theorem 1.3 in \cite{Abtbmsie} and will not be repeated.
\ig
\subsection{From K-martingale problems to truncated $\beta$-ISLTRW SIEs}\lbl{kmart2trnctd}
We now establish \thmref{kmarteqaux}.
\bpfs{Proof of \thmref{kmarteqaux}}
Assume that $\lpa\lbr\xtaunlf \rbr_{\tau\ge0},\lbr\wtyn\rbr_{y\in\Xnldf}\rpa$ is a solution to
$\lbr\isltrwsienlaux\rbr_{\tau\ge0}$ in \eqref{aux} on a filtered probability space $\OFFtP$; then clearly $\xztauxnlf=\unx$ for $x\in\Xnd$
and $\xttauxnlf=\uttxnD$ for $x\in\Xnd\sm\Xnldf$, for every $\tau\ge0$ and $t\in[0,\tau]$.  Fix any arbitrary $\tau_1,\tau_2>0$  and $x_1,x_2\in\Xnldf$.
By It\^o's formula we have
\beq\lbl{itokmart}
\bsp
&f(Z_{\beta,n,l}^{x_{1,2},\tau_{1,2}}(t))-f(u_0^{x_1,x_2})-\ds\int_0^t \lpa\sAKtauot f\rpa(s,x_1,x_2,Z_{\beta,n,l}^{\tau_{1,2}})ds
\\&=\sum_{i=1}^2\sum_{y\in\Xnldf}\int_0^t\partial_if(Z_{\beta,n,l}^{x_{1,2},\tau_{1,2}}(s))\kappa^{x_i,y}_{\h,s,\tau_i}\lpa\xstauiynlf\rpa{d\wsyn}\in\locm
\end{split}
\eeq
for $t\in[0,\tau_1\wedge\tau_2]$ and \eqref{kmart} is satisfied.

Conversely, if a family of adapted $$\lbr\xtaunlf\rbr_{\tau\ge0}=\lbr\lbr\xttauxnlf;0\le t\le\tau, x\in\Xnd\rbr\rbr_{\tau\ge0}$$ defined on a
probability space $\OFFtP$ satisfies \eqref{kmart}; then fixing any two $\tau_{1},\tau_{2}>0,$ letting $\B_{0,R}=\lbr u=(u_1,u_2)\in\Rs; \lab u\rab\le R\rbr$
and choosing $f_1,f_2,f\in\C_b^2(\R;\R)$ such that $f_i(u)=u_i$ for $i=1,2$ and $f(u)=u_1u_2$ whenever $u\in\B_{0,R}$ we see that
\beq\lbl{2mart}
\bc
\mtautxi&:=\xttauxinlf-\uttxinD\in\locm; i=1,2,\\
\ntautxot&:=\ds\prod_{i=1}^2\xttauxinlf-\prod_{i=1}^2\unxi
\\&\ds-\sum_{\substack{1\le i,j\le2\\i\neq j}}\int_0^t\xstauxinlf d\utsxjnD
\\&\ds-\int_0^t\sum_{y\in\Xnldf}\Upsilon^{x_{1,2},y}_{\h,t,\tau_{1,2}}\lpa Z_{\beta,n,l}^{y,\tau_{1,2}}(s)\rpa ds\in\locm
\ec
\eeq
for all $t\in[0,\tau_{1}\wedge\tau_{2}]$ and $x_1,x_2\in\Xnldf$, where we have used the notation $Z_{\beta,n,l}^{\tau_{1,2}}$ for the two-dimensional process defined in \eqref{Z}.  We then have that
\beq\lbl{qmart}
\bsp
&\prod_{i=1}^2\mtautxi-\sum_{y\in\Xnldf}\int_0^t\Upsilon^{x_{1,2},y}_{\h,t,\tau_{1,2}}\lpa Z_{\beta,n,l}^{y,\tau_{1,2}}(s)\rpa ds
\\&=\ T_1^{x_{1,2},\tau_{1,2}}(t)+T_2^{x_{1,2},\tau_{1,2}}(t)
\end{split}
 \eeq
where
\beq\lbl{m1}
T_1^{x_{1,2},\tau_{1,2}}(t):=\ntautxot-\sum_{\substack{1\le i,j\le2\\i\neq j}}\unxi\mtautxj\in\locm
\eeq
and
\beq\lbl{m2}
\bsp
T_2^{x_{1,2},\tau_{1,2}}(t):&=\sum_{\substack{1\le i,j\le2\\i\neq j}}\int_0^t\lbk\xstauxinlf-\xttauxinlf\rbk d\utsxjnD
+\prod_{i=1}^2\lbk\uttxinD-\unxi\rbk
\\&=\sum_{\substack{1\le i,j\le2\\i\neq j}}\int_0^t\lbk\unxi-\utuxinD\rbk d\mtauuxj\in\locm.
\end{split}
\eeq
Thus,
\beq\lbl{qvkmart}
\lqv\mtauxo,\mtauxt\rqv_t= \sum_{y\in\Xnldf}\int_0^t\Upsilon^{x_{1,2},y}_{\h,t,\tau_{1,2}}\lpa Z_{\beta,n,l}^{y,\tau_{1,2}}(s)\rpa  ds.
\eeq
Equations \eqref{2mart} and \eqref{qvkmart} imply that there exists a set of $r=\#\lbr y;y\in\Xnldf\rbr$ independent
Brownian motions $\lbr\wtyn;t\in\Rp\rbr_{y\in\Xnldf}$ on an extension $\OFFtPt$ such that  
\beq\lbl{repres}
\mtautx=\sum_{y\in\Xnldf}\int_0^t\kappa^{x,y}_{\h,s,\tau}\lpa\xstauynlf\rpa{d\wsyn},\ \forall(x,\tau,t)\in\Xnldf\times\Rp\times[0,\tau].
\eeq
In fact, fixing any $\tau>0$ and labeling the $\lbr x;x\in\Xnldf\rbr$ as $\lbr x_1,\ldots,x_r\rbr$,
the restriction of the desired family of BMs $\lbr\wtyn;t\in\Rp\rbr_{y\in\Xnldf}$ to the time interval $[0,\tau]$ is obtained from the matrix equation (written in differential form)
\beq\lbl{MatrixBMs}
\bsp
\begin{bmatrix}d{W}_n^{x_1}(t)\\\vdots\\d{W}_n^{x_r}(t)\end{bmatrix}&=\begin{bmatrix}\kappa^{x_1,x_1}_{\h,t,\tau}\lpa\xttauxonlf\rpa&\cdots&\kappa^{x_1,x_r}_{\h,t,\tau}\lpa\xttauxrnlf\rpa\\
\vdots&\vdots&\vdots\\
\kappa^{x_r,x_1}_{\h,t,\tau}\lpa\xttauxonlf\rpa&\cdots&\kappa^{x_r,x_r}_{\h,t,\tau}\lpa\xttauxrnlf\rpa\end{bmatrix}^{-1}
\begin{bmatrix}dM^{\tau,x_{1}}(t)\\\vdots\\dM^{\tau,x_{r}}(t)\end{bmatrix}
\\&:=\begin{bmatrix}d\tilde{W}_n^{\tau,x_1}(t)\\\vdots\\d\tilde{W}_n^{\tau,x_r}(t)\end{bmatrix}
\end{split}
\eeq
whenever the middle inverse kernel-diffusion coefficient $r\times r$ matrix, denoted by $A^{-1}$, exists (the determinant  ${\mathrm{det}}(A)\neq0$) almost surely. If this fails we can proceed similar to the standard finite dimensional SDE case cf. Ikeda and Watanabe \cite{IW} or Doob \cite{Do}.  It is now straightforward to verify that, for any $\tau_{1},\tau_{2}>0$ and any $t\in[0,\tau_{1}\wedge\tau_{2}]$, the two families of BMs $\lbr \lbr\tilde{W}_n^{\tau_{k},y}\rbr_{y\in\Xnldf};k=1,2\rbr$ satisfy
\beq\lbl{covbms}
\langle \tilde{W}_{n}^{\tau_{1},x_{i}}(\cdot), \tilde{W}_{n}^{\tau_{2},x_{j}}(\cdot)\rangle_{t}
=\bc 
t,&1\le i=j\le r\\
0,&1\le i\neq j\le r
\ec 
\eeq 
almost surely, whether or not $\tau_{1}=\tau_{2}$.  I.e., we get one family of independent BMs $\lbr\wtyn\rbr_{y\in\Xnldf}$, such that $\lpa\lbr\xtaunlf\rbr_{\tau\ge0},\lbr W_n^y\rbr_{y\in\Xnldf}\rpa$ is a solution pair to  $\lbr\isltrwsienlaux\rbr_{\tau\ge0}$ in \eqref{aux} on $\OFFtPt$.
Hence,  on the probability space $\OFFtPt$, the pair $\lpa\lbr\xttxnlf\rbr_{t\ge0,x\in\Xnd},\lbr\wtyn\rbr_{y\in\Xnldf}\rpa$ solves the $l$-truncated $\beta$-ISLTRW SIE $\isltrwsienl$ in \eqref{trnctdsie} on $\Rp\times\Xnd$.
\epfs
\fi
\subsection{Completing the proof of the second main result}\lbl{lipvnlip}
We now complete the proof of \thmref{mainthm2}.  In \secref{regandtight} and \secref{regandtight2} we assumed the existence of a $\beta$-ISLTRW SIE solution and
we obtained regularity and tightness for the sequence of lattice SIEs $\lbr\isltrwsien\rbr_{n\in\Ns}$.  This, in turn,
implied the existence and regularity for a $\beta$-ISLTRW SIE limit solution to our $\isltsie$ in \eqref{isltbmsie}.   To complete the existence of the desired double limit solution\footnote{The type of our lattice limit  solution to $\isltsie$ in \eqref{isltbmsie} depends on the conditions: under the Lipschitz conditions \eqref{lcnd} we get a direct solution to the lattice SIE $\isltrwsien$ for every $n$ and a direct $\beta$-ISLTRW SIE limit solution to $\isltsie$ (see \thmref{latlimlip}); whereas under the non-Lipschitz conditions in \eqref{cnd} we obtain a limit $\beta$-ISLTRW SIE solution, thanks to our K-martingale approach, and a $\beta$-ISLTRW SIEs double limit solution to $\isltsie$.} for $\isltsie$ it suffices then to prove the existence of a solution to $\isltrwsien$ for each fixed $n\in\Ns$, under the condition \eqref{cnd}, that is uniformly $L^{p}(\Omega,\P)$ bounded on $[0,T]\times \X^d$ for every $T>0$ and every $p\ge2$.  We establish this existence via the K-martingale approach just recalled and adapted from \cite{Abtbmsie}, using \thmref{kmarteqaux}.

First, the following proposition summarizes the results in this case for the $\beta$-ISLTRW SIEs spatial lattice scale\footnote{We remind the reader that we will, without further notice, suppress the dependence on $\beta$ whenever it is more convenient notationally to do so.}.

\bpr[Existence for $\beta$-ISLTRW SIEs with non-Lipschitz $a$]\lbl{ewlip}
Assume the conditions \eqref{cnd} hold.  Then,
\begin{enumerate}\renewcommand{\labelenumi}{$($\alph{enumi}$)$}
\item For every $(n,l)\in\Ns\times\N$, every $\beta=1/\nu\in\lbr1/2^{k},k\in\N\rbr$, and for every $p\ge2$, there exists an $L^p$-bounded solution $\uttxnlf$ to the truncated $\beta$-ISLTRW SIE
\eqref{trnctdsie} on $\T\times\Xnd$.  Moreover, if we linearly interpolate $\uttxnlf$ in space; then, with probability one,
the continuous map $(t,x)\mapsto\uttxnlf$ is locally $\gamma_t$-H\"older continuous in time with
$\gamma_t\in\lpa0,\tf{2\nu-d}{4\nu}\rpa$ for $\nu=\beta^{-1}\in\nuset$ and $d=1,2,3$.
\item  For any fixed $n\in\Ns$, the sequence $\lbr\uttxnlf\rbr_{l\in\N}$ of linearly-interpolated solutions in $(a)$
has a subsequential weak limit $\utn$ in $\C(\T\times\Rd;\R)$.  We thus have a limit solution $\utn$
to $\isltrwsien$, and $\utn$ is locally $\gamma_t$-H\"older continuous in time with
$\gamma_t\in\lpa0,\tf{2\nu-d}{4\nu}\rpa$ for $\nu=\beta^{-1}\in\nuset$ and $d=1,2,3$.
\end{enumerate}
\epr
\bpf
\begin{enumerate}\renewcommand{\labelenumi}{$($\alph{enumi}$)$}
\item  First, recall that the deterministic term $\uttxD$ in \eqref{trnctdsie} is completely determined by $\un$.
Moreover, under the conditions in \eqref{cnd} on $\un$, $\uttxD$ is clearly bounded and it is smooth in time as
in \remref{smuthdet}.  Fix an arbitrary $T>0$, and let $\T=[0,T]$.  We now prove the existence of a family of adapted processes $\lbr\tilde{X}^\tau_{\beta,n,l}\rbr_{\tau\in\T}$ satisfying our K-martingale problem
\eqref{kmart}, which by \thmref{kmarteqaux} implies the existence of a solution to the $l$-truncated $\beta$-ISLTRW SIE
\eqref{trnctdsie} on $\T\times\Xnd$.  On a probability space $\OFFtP$ we prepare  a family of $r$-independent BMs $\lbr\wtyn\rbr_{y\in\Xnldf}$.  For each $\tau\in\T$ and each $i=1,2,\ldots$ define a continuous process $X_{\beta,n,l,i}^\tau$ on $[0,\tau]\times\Xnd$
inductively for $k/2^i\le t\le((k+1)/2^i)\wedge\tau$
$(k=0,1,2,\ldots)$ as follows: $X_{\beta,n,l,i}^{x,\tau}(0)=\unx$ $(x\in\Xnd)$ and if $X_{\beta,n,l,i}^{x,\tau}(t)$ is defined for $t\le k/2^i$, then we define
$X_{\beta,n,l,i}^{x,\tau}(t)$ for $k/2^i\le t\le((k+1)/2^i)\wedge\tau$, by
\beq\lbl{indproc}\bsp
&X_{\beta,n,l,i}^{x,\tau}(t)\\&=\bc
\ds X_{\beta,n,l,i}^{x,\tau}\lpa\tf{k}{2^{i}}\rpa+\ds\sum_{y\in\Xnldf}\kappa^{x,y}_{\h,\frac{k}{2^i},\tau}\lpa X_{\beta,n,l,i}^{y,\tau}(\tfrac{k}{2^i})\rpa\lpa\Delta_{t,\tf{k}{2^{i}}}W_{n}^{y}\rpa\\ +\lbk\uttxnD-\tilde{U}_{n,D}^{x}\lpa\frac{k}{2^{i}}\rpa \rbk; &x\in\Xnldf, \\
\ds\uttxnD;&x\in\Xnd\setminus\Xnldf,
\ec
\end{split}
\eeq
where $\Delta_{t,\tf{k}{2^{i}}}W_{n}^{y}=\wtyn-W_n^y(\tfrac{k}{2^i})$.   Clearly, $X_{\beta,n,l,i}^\tau$ is the solution to the equation
\beq\lbl{indsie}\bsp
&X_{\beta}^{x,\tau}(t)\\&=
\bc
\ds\sum_{y\in\Xnldf}\int_0^t\kappa^{x,y}_{\h,\phi_i(s),\tau}\lpa X^{y,\tau}(\phi_i(s))\rpa{d\wsyn}
+\uttxnD;&x\in\Xnldf,\\
\ds\uttxnD;&x\in\Xnd\setminus\Xnldf
\ec
\end{split}
\eeq
with $X_{\beta}^{x,\tau}(0)=\unx$, where $\phi_i(t)=k/2^i$ for $k/2^i\le t<(k+1)/2^i\wedge\tau$ $(k=0,1,2,\ldots)$.

Now, for $q\ge1$, let $M_{\beta,q,l,i}^\tau(t)=\sup_{x\in\Xnd}\E\lab X_{\beta,n,l,i}^{x,\tau}(t)\rab^{2q}$.  By the boundedness of
$\uttxnD$ over the whole infinite lattice $\Xnd$, we have
\beq\lbl{trnctdlpbdd}
\bsp
M_{\beta,q,l,i}^\tau(t)\le C+\sup_{x\in\Xnldf}\E\lab X_{\beta,n,l,i}^{x,\tau}(t)\rab^{2q}
\end{split}
\eeq
Then, replacing $\Xnd$ by $\Xnldf$ and following the same steps as in the proof of \propref{Expbd}, we get that
\beq\lbl{trnctdulpbdd}
\sup_{\tau\in\T}\sup_{t\in[0,\tau]}M_{\beta,q,l,i}^\tau(t)\le C,\quad d=1,2,3,
\eeq
where, here and in the remainder of the proof, the constant $C$ depends only on $q$, $\beta$, $\max_x|u_0(x)|$, the spatial dimension $d=1,2,3$, and $T$ but may change its value from one line to the next.  Remembering that $\h\searrow0$ as $n\nearrow\infty$ and $n\in\Ns$,
the independence in $l$ is trivially seen since \lemref{2ndQinequality} implies
$$\sum_{y\in\Xnldf} \lbk\Ktxyn\rbk^2\le\sum_{y\in\Xnd} \lbk\Ktxyn\rbk^2\le \df{C}{t^{d/2\nu}};\forall d=1,2,3, l\in\N$$
Similarly, letting $X_{n,l,i,R}^{x,\tau}$ denote the random part of $X_{\beta,n,l,i}^{x,\tau}$ on the truncated lattice $\Xnldf$, using
\eqref{trnctdulpbdd}, and repeating the arguments in \lemref{sptdiff} and \lemref{tmpdiff}---replacing $\Xnd$ by $\Xnldf$ and
noting that \lemref{4thQinequality} and \lemref{3rdQinequality} hold on $\Xnldf$---we obtain
\beq\lbl{trnctdspttmpdiff}
\bsp
&\E\lab X_{\beta,n,l,i,R}^{x,\tau_{1}}(t) - X_{\beta,n,l,i,R}^{y,\tau_{1}}(t)\rab^{2q}+\E\lab X_{\beta,n,l,i,R}^{x,\tau_{2}}(t) - X_{\beta,n,l,i,R}^{y,\tau_{2}}(t)\rab^{2q}\\&\le C_{d} |x-y|^{2q\alpha_{d}};\alpha_{d}\in I_{d},\\
&\E\lab X_{\beta,n,l,i,R}^{x,\tau_{1}}(t) - X_{\beta,n,l,i,R}^{x,\tau_{1}}(r)\rab^{2q}+\E\lab X_{\beta,n,l,i,R}^{x,\tau_{2}}(t) - X_{\beta,n,l,i,R}^{x,\tau_{2}}(r)\rab^{2q}\\&\le C\lab t-r\rab^{\tf{(2\nu-d)q}{2\nu}},
\end{split}
\eeq
for all $x,y \in \Xnldf$,  $r,t\in[0,\tau_{1}\wedge\tau_{2}]$, $\tau_{1},\tau_{2}\in\T$, and $d=1,2,3$.   It follows that, for every point $\tau_{1,2}=\lpa\tau_{1},\tau_{2}\rpa\in\T^{2}$, there is a subsequence $\lbr\lpa\tilde{X}_{\beta,n,l,i_{m}}^{\tau_{1}},\tilde{X}_{\beta,n,l,i_{m}}^{\tau_{2}}\rpa\rbr_{m=1}^\infty$ on a probability space $\OtauotFtauotPtauott$ such that $\lpa\tilde{X}_{\beta,n,l,i_{m}}^{\tau_{1}},\tilde{X}_{\beta,n,l,i_{m}}^{\tau_{2}}\rpa\overset{\mathscr L}{=} \lpa{X}_{\beta,n,l,i_{m}}^{\tau_{1}},{X}_{\beta,n,l,i_{m}}^{\tau_{2}}\rpa $ and $$\lpa\tilde{X}_{\beta,n,l,i_{m}}^{x,\tau_{1}}(t),\tilde{X}_{\beta,n,l,i_{m}}^{x,\tau_{2}}(t)\rpa\longrightarrow \lpa\tilde{X}_{\beta,n,l}^{x,\tau_{1}}(t),\tilde{X}_{\beta,n,l}^{x,\tau_{2}}(t)\rpa$$  uniformly on compact subsets
of $[0,\tau_{1}\wedge\tau_{2}]\times\Xnd$, as $m\to\infty$ a.s.  Let $\TQ=\T\cap\Q$, where $\Q$ is the set of rationals, and define the product probability space 
$$\OFPt:=\lpa{\bigotimes_{\tau_{1,2}\in\TQ^{2}}}\tilde{\Omega}_{\tau_{1,2}},\bigotimes_{\tau_{1,2}\in\TQ^{2}}\tilde{\sF}_{\tau_{1,2}},\bigotimes_{\tau_{1,2}\in\TQ^{2}}\tilde{\P}_{\tau_{1,2}}\rpa.$$ 
If $s<t$, then for every $ f\in\C_b^2(\Rs;\R)$, $\tau_1,\tau_2\in\TQ\setminus\lbr0\rbr$,
$t\in[0,\tau_1\wedge\tau_2]$, $x_1,x_2\in\Xnldf$, and for every bounded continuous $F:\C\lpa\Rp;\Rs\rpa\to\R$
that is ${\mathscr{B}}_s\lpa\C\lpa\Rp;\Rs\rpa\rpa:=\sigma\lpa z(r); 0\le r\le s\rpa$-measurable function, we have 
\beq\lbl{mmart}
\bsp
&\EPt\lbk\lbr f(\tilde{Z}_{\beta,n,l}^{x_{1,2},\tau_{1,2}}(t))-f(\tilde{Z}_{\beta,n,l}^{x_{1,2},\tau_{1,2}}(s))\right.\right.
\\&-\left.\left.\int_s^t \lpa\sAKtauot f\rpa(r,x_1,x_2,\tilde{Z}_{\beta,n,l}^{\tau_{1,2}})dr\rbr F\lpa\tilde{Z}_{\beta,n,l}^{x_{1,2},\tau_{1,2}}(\cdot)\rpa\rbk \\
&=\lim_{m\to\infty}\EPt\lbk\lbr f(\tilde{Z}_{\beta,n,l,i_m}^{x_{1,2},\tau_{1,2}}(t))-f(\tilde{Z}_{\beta,n,l,i_m}^{x_{1,2},\tau_{1,2}}(s))\right.\right.
\\&-\left.\left.\int_s^t \lpa\sAKtauotim f\rpa(r,x_1,x_2,\tilde{Z}_{\beta,n,l,i_m}^{\tau_{1,2}})dr\rbr F\lpa\tilde{Z}_{\beta,n,l,i_m}^{x_{1,2},\tau_{1,2}}(\cdot)\rpa\rbk=0,
\end{split}
\eeq
where, by a standard localization argument, we have assumed that $a$ is also bounded; and where $\tilde{Z}_{\beta,n,l}^{\tau_{1,2}}$ and $\tilde{Z}_{\beta,n,l,i_m}^{\tau_{1,2}}$ are obtained from the definition of
${Z}_{\beta,n,l}^{\tau_{1,2}}$ in \eqref{Z} by replacing $X_{\beta,n,l}^{\tau_j}$ by $\tilde{X}_{\beta,n,l}^{\tau_j}$ and $\tilde{X}_{\beta,n,l,i_m}^{\tau_j}$, $j=1,2$,
respectively.  The operator $\sAKtauotim$ is obtained from $\sAKtauot$ by replacing $\Upsilon^{x_{i,j},y}_{\h,t,\tau_{i,j}}\lpa u^y(t)\rpa$ in \eqref{siegen} by $\Upsilon^{x_{i,j},y}_{\h,\phi_{i_m}(t),\tau_{i,j}}\lpa u^y(\phi_{i_m}(t))\rpa$.  Also, obviously, for any $\tau\in\TQ$ and $t\in[0,\tau]$
\beq\lbl{mmartd}
\tilde{X}_{\beta,n,l}^{x,\tau}(t)=\lim_{m\to\infty}\tilde{X}_{\beta,n,l,i_m}^{x,\tau}(t)=\uttxnD; \qquad x\in\Xnd\setminus\Xnldf, \mbox{ a.s.~}\Pt.
\eeq
It follows from \eqref{mmart} and \eqref{mmartd} that  $\lbr\tilde{X}_{\beta,n,l}^{\tau}\rbr_{\tau\in\TQ}$ satisfies the K-martingale problem \eqref{kmart} with respect to the filtration $\{\tilde{\sFt}\}$, with
$$\tilde{\sFt}=\bigcap_{\epsilon>0}\sigma\lbr \tilde{X}_{\beta,n,l}^{x,\tau}(u);u\le (t+\epsilon)\wedge\tau, \tau\in\TQ\cap(t,T]\rbr.$$
Thus, by \thmref{kmarteqaux}, with $\tau\in\Rp$ replaced by $\tau\in\TQ$, there is a solution $\uttxnlf$ to the $l$-truncated $\beta$-ISLTRW SIE \eqref{trnctdsie} on $\TQ\times\Xnd$.  Use continuous extension in time of $\uttxnlf$ to extend its definition to $\T\times\Xnd$, and denote the extension also by $\uttxnlf$.  Clearly $\uttxnlf$ solves the $l$-truncated $\beta$-ISLTRW SIE \eqref{trnctdsie} on $\T\times\Xnd$.

Now, for $q\ge1$, let $M_{\beta,q,l}(t)=\sup_{x\in\Xnd}\E\lab\uttxnlf\rab^{2q}$.  As above, the boundedness of $\uttxnD$, implies
\beq\lbl{trnctdlpbdd2}
\bsp
M_{\beta,q,l}(t)\le C+\sup_{x\in\Xnldf}\E\lab\uttxnlf\rab^{2q}.
\end{split}
\eeq
Then, replacing $\Xnd$ by $\Xnldf$ and following the same steps as in the proof of \propref{Expbd}, we get that
\beq\lbl{trnctdulpbdd2}
M_{\beta,q,l}(t)\le C,\quad\forall t\in\T,\beta\in\lbr1/2^{k};k\in\N\rbr\mbox{ and }d=1,2,3.
\eeq

Similarly, letting $\uttxnlfR$ denote the random part of $\uttxnlf$ on the truncated lattice $\Xnldf$, using
\eqref{trnctdulpbdd2}, and repeating the arguments in \lemref{sptdiff} and \lemref{tmpdiff}---replacing $\Xnd$ by $\Xnldf$ and
noting that the inequalities in \lemref{4thQinequality} and \lemref{3rdQinequality} trivially hold if we replace $\Xnd$ by $\Xnldf$---we
obtain
\beq\lbl{trnctdspttmpdiff2}
\bsp
&\E\lab\uttxnlfR - \uttynlfR\rab^{2q}\le C_{d} |x-y|^{2q\alpha_{d}};\ \alpha_{d}\in I_{d},\\
&\E\lab\uttxnlfR - \utrxnlfR\rab^{2q}\le C\lab t-r\rab^{\tf{(2\nu-d)q}{2\nu}},
\end{split}
\eeq
for all $x,y \in \Xnldf$, $r,t\in\T$, and $d=1,2,3$.   By \remref{smuthdet},
$\uttxnD$ is differentiable in $t$.  So, linearly interpolating $\uttxnlf$ in space and using  \eqref{trnctdspttmpdiff2} and
arguing as in the proof of part (a) of \lemref{regtight}, we get that the continuous map $(t,x)\mapsto\uttxnlf$
is locally $\gamma_t$-H\"older continuous in time with $\gamma_t\in\lpa0,\tf{2\nu-d}{4\nu}\rpa$ for $\nu=\beta^{-1}\in\nuset$ and $d=1,2,3$.
\item  Clearly, $\uttxnD$ in \eqref{trnctdsie} is the same for every $l$, so it is enough to show convergence of the random part $\uttxnlfR$.
Using \eqref{trnctdspttmpdiff2} we get tightness for $\lbr\uttxnlfR\rbr_l$ and consequently a subsequential weak limit $\utn$,
which  is our limit solution for $\isltrwsien$.  For the regularity assertion, $\uttxnD$ is smooth and bounded as noted above.  So,
using \eqref{trnctdulpbdd2} and \eqref{trnctdspttmpdiff2}, and imitating the argument in the proof of part (b) of \lemref{regtight}
(remembering that here we are taking the limit as $l\to\infty$); we
get the desired $L^p$ boundedness for $\utn$ as in \propref{Expbd} and the spatial and temporal moments bounds in
\lemref{sptdiff} and \lemref{tmpdiff}\beq\lbl{est}
\bc
\E\lab\uttxn\rab^{2q}\le C&\\
\E\lab\uttxnR - \uttynR\rab^{2q}\le C_{d} |x-y|^{2q\alpha_{d}};\ \alpha_{d}\in I_{d},&\\
\E\lab\uttxnR - \utrxnR\rab^{2q}\le C\lab t-r\rab^{\tf{(2\nu-d)q}{2\nu}},
\ec
\eeq
for $(t,x,n)\in\T\times\Xnd\times\Ns$ and for $\nu=\beta^{-1}\in\nuset$, $d=1,2,3$, and $q\ge1$ and the desired H\"older regularity follows.
\end{enumerate}
The proof is complete.
\epf
We now get \thmref{mainthm2} for $\btsie$ as the following corollary.
\bcr \thmref{mainthm2} holds.
\ecr
\bpf
The desired conclusion follows upon using the argument in the proof of part (b) of \lemref{regtight} along with \defnref{limitsolns} and the $L^p $-boundedness and the spatial and temporal moments bounds for $\lbr \utn\rbr_n$ that we got in \eqref{est} above.
\epf
\ig
\section{Conclusions }\lbl{conc}
We considered the multiplicative noise case for our recently introduced (see \cite{Abtpspde})
 fourth order BTBM stochastic integral equation $\btsie$ in \eqref{btpsol} on $\Rp\times\Rd$,
 under both Lipschitz and less than Lipschitz conditions on the diffusion coefficient $a$.
 We showed striking spatio-temporal dimension-dependent H\"older regularity for such SIEs that are not only real-valued up to spatial dimension $d=3$, but---even more impressively---are spatially nearly locally Lipschitz, and hence roughly twice as smooth in space as the Brownian sheet of the driving noise in $d=1,2$.    This gives, for the first time, an example of a kernel---that is also the density of an interesting stochastic process---that is able to regularize a space-time white noise driven equation so that its solutions are pushed beyond the traditional nearly H\"older $1/2$ spatial regularity.   Of course this contrasts sharply with  second order reaction-diffusion (RD) SPDEs driven by space-time white noise whose counterpart real-valued solutions are confined to the case $d=1$.  
 
 We analyzed our BTBM SIE using both a direct and a numerically-flavored lattice approach similar to  our discretized method for the simpler second order RD SPDEs driven by space-time white noise (see \cite{Asdde1,Asdde2}).  We used the second approach in the case of non-Lipschitz conditions.  In it, we discretize space and formulate solutions to the resulting spatially-discrete stochastic integral equations
 in terms of the density of Brownian-time random walk or BTRW---which we introduce in this article along with the
 general class of Brownian-time chains (BTCs), of which $\beta$-ISLTRW is a special case.
 As with their continuous counterpart BTPs, which we introduced in \cite{Abtp1,Abtp2},
 BTCs are interesting new processes outside of the current well established
 theory; and we believe they merit further study.  In the course of proving our results, we prove several interesting facts about the BTRW.   These include a connection to fourth order differential-difference equation that is proved in Appendix \ref{appA} and different estimates which lead to the definition of $2$-Brownian-times Brownian motion and $2$-Brownian-times random walk.    We define two notions of solutions to the lattice model: direct solutions and limit solutions (from a finite truncation of the lattice to the whole lattice).  These solutions (both direct and limit) are then used to define two types of $\beta$-ISLTRW SIEs limit solutions to $\btsie$ (direct limit solutions and double limit solutions), as the size of the lattice mesh shrinks to zero.  
 
To deal with existence under the non-Lipschitz condition on $a$, we introduce our K-martingale approach, which is tailor-made for kernel SIEs as $\btsie$ as well as for other mild formulations of many other different SPDEs.  It is a delicate variant of the well known, and by now classic, martingale problem approach of Stroock and Varadhan for SDEs.  This K-martingale approach starts by constructing an auxiliary problem to a truncated lattice version of \eqref{btpsol}, for which the existence of solutions implies solutions existence for the truncated lattice model.  We then formulate a martingale problem equivalent to the auxiliary problem (the K-martingale problem).    
  A key advantage of the K-martingale approach is that it is a unified framework in which the existence and uniqueness\footnote{As we noted earlier, we don't prove uniqueness under less-than-Lipschitz conditions, and we therefore don't discuss further the implications of our K-martingale approach to uniqueness in this article.} of many kernel stochastic integral equations, which are the mild formulation for many SPDEs, may be treated using only variants of the kernel SIE.  This includes SPDEs of different orders (including second and fourth), so long as the corresponding spatially-discretized kernel (or density) satisfies Kolmogorov-type bounds on its temporal and spatial differences.  In essence, what the K-martingale approach implies is that if the kernel in the lattice model is nice enough for the lattice model to converge as the lattice mesh shrinks to zero (under appropriate assumptions on $a$), then it is nice enough to guarantee a solution for the lattice model.  We use it here to prove the existence of $\beta$-ISLTRW SIEs double limit solutions to \eqref{btpsol} under the conditions \eqref{cnd}, but just as with the Stroock-Varadhan method, it can handle uniqueness as well.  
The densities of $\beta$-ISLTRW and BTBM have a considerable regularizing effect on stochastic kernel equations driven by space-time white noise as compared to the standard Brownian motion or continuous time random walk densities (the usual green functions for second order RD equations and their spatially-discretized versions).  The unconventional memory-preserving fourth order PDEs associated with the BTP density are
highly regular: their solutions are $\mathrm{C}^{1,4}(\Rp\times\Rd,\R)$ for all times and all $d\ge1$, despite the fact that
a positive bi-Laplacian term is part of the equation \cite{Abtp1,Abtp2}.  However, this bi-Laplacian is coupled in 
a very specific way---dictated by the BTBM probability density function---with a Laplacian acting on the smooth initial data and whose coefficient grows arbitrarily large as time approaches the initial time (zero) at the rate of $1/\sqrt{8\pi t}$.  One way to understand this smoothing effect is to note that the BTBM density in $\btsie$ is intimately connected---and shares regularity properties---with the kernel associated with our recently introduced imaginary-Brownian-time-Brownian-angle process, which we used to give a solution to a Kuramoto-Sivashinsky-type PDE in \cite{Aks}.   In the stochastic setting of this article, this BTBM density smoothing effect on $\btsie$ is evidenced in a regularity of solutions that is much higher than typical second order space-time white noise driven RD SPDEs.  This regularizing effect is such that we are able to obtain $\gamma$-H\"older continuous solutions to our BTBM SIE $\btsie$ for spatial dimensions $d=1,2,3$.   We show that the H\"older exponent is dimension-dependent with $\gamma\in(0,\tf{4-d}{8})$, $d=1,2,3$.  In addition, we are able to show ultra spatial regularity by showing a nearly local Lipschitz regularity for $d=1,2$, and nearly local H\"older $1/2$ regularity in $d=3$.   To get the smoothing effect of the BTBM density, we prove BTBM and $\beta$-ISLTRW estimates that enable us to extract a $\beta$-ISLTRW SIEs weak limit (direct limit in the Lipschitz case and double limits in the non-Lipschitz one) H\"older continuous solution to $\btsie$ in spatial dimensions $d=1,2,3$, in spite of the presence of the driving space-time white noise.  Again, the effective H\"older exponent depends on the space dimension through the expression $(4-d)/8$.  This ultra regularity in multispatial dimensions naturally motivates the study of the variations (temporal and spatial) of BTBM SIEs, which we undertake---among other aspects of BTBM SIEs---in \cite{AX} and a followup article.

Encoding this smoothing effect from our BTBM SIE $\btsie$ into a fourth order SPDE
involving the bi-Laplacian coupled with a Laplacian term requires extra parameters.  We
give what we call the fourth order parametrized SPDE corresponding to $\btsie$, linking
the spatially-discretized $\beta$-ISLTRW SIE to the diagonals of a parametrized stochastic differential-difference equation
($\beta$-ISLTRW PSDDE) on the lattice.

  As we recently started doing for PDEs \cite{Aks}, we adapt the methods presented here and in \cite{Aks} to give an
  entirely new approach---in terms of our Linearized Kuramoto-Sivashinsky process (or imaginary-Brownian-time-Brownian-angle process) and related processes---to study the multi-spatial dimensions SPDEs version of famous fourth order applied mathematics PDEs like the Kuramoto-Sivashinsky (several different versions), the Cahn-Hilliard, and the Swift-Hohenberg PDEs.   We  illustrate this in upcoming papers (\cite{Aksspde,AL,AD}) and planned followup papers.  Traditional semigroup analytical methods alone are not adequate for this since the existence of the KS semigroup in $d>1$ is not settled analytically.   In another direction, we have discovered interesting connections between BTBM and related processes and stochastic fractional PDEs, we address these connections and their consequences in upcoming articles as well.

  Also, SIEs corresponding to other BTP processes we introduced in \cite{Abtp2} may also be studied by adapting
  and generalizing our approach here.    We believe BTPs, their PDEs, their SIEs, and their discretized cousins (the BTCs and their equations) can play a useful role by adding new, currently unavailable, insights and models to the ever growing mathematical finance theory (see \cite{Carr} for an example).    We also hope to explore these aspects in future papers.  
\fi
\appendix
\ig
\section{Proof of BTRW-DDE Connection}\lbl{appA}
In this appendix, we give the proof of \lemref{fodde} linking the density of $\beta$-ISLTRW to fourth order differential-difference equations.
\bpfs{Proof of \lemref{fodde}}
Let $u^x_n(t)=\E \lbk u_0\lpa \S^x_{B,\h}(t)\rpa\rbk$ with $u_0$ as in \eqref{cnd}.
Observe that
\beq\lbl{sg}
\E \lbk u_0\lpa \S^x_{B,\h}(t)\rpa\rbk=2\int_0^\infty\ptzs\E\lbk u_0\lpa S^x_{\h}(s)\rpa\rbk ds
\eeq
where $\ptzs$ the transition density of the one dimensional BM $B(t)$. Differentiating \eqref{sg} with
respect to $t$ and putting the derivative under the integral, which is easily justified
by the dominated convergence theorem, then using the fact that
$\partial\ptzs/\partial t=\tf12\partial^2 \ptzs/\partial s^2$ we have
\beq\lbl{tder}
\bsp
\df{d}{d t}\E \lbk u_0\lpa \S^x_{B,\h}(t)\rpa\rbk&=
2\int_0^\infty\df{\pa}{\pa t}\ptzs\E\lbk u_0\lpa S^x_{\h}(s)\rpa\rbk ds
\\&=\int_0^\infty\df{\pas}{\pa s^2}\ptzs\E\lbk u_0\lpa S^x_{\h}(s)\rpa\rbk ds.
\end{split}
\eeq
Letting $\sTns \unx=\E \un(S^x_{\h}(s))$ be the action of the semigroup $\sTns$
associated with the standard continuous-time symmetric random walk $S^x_{\h}$ on the lattice $\Xnd$.
Then, the generator of $S^x_{\h}$ on $\Xnd$ is given by $\mathscr{A}_n=\Dn/2$.
Alternatively, noting that $\qtxn$ is
the fundamental solution to the deterministic heat equation \eqref{latheat} on the lattice $\Xnd$, we
get
\beqs
\bsp
\df{d}{d s}\sTns \unx&=\df{d}{d s}\E \lbk u_0\lpa S^x_{\h}(s)\rpa\rbk=
\sum_{y\in\Xnd}\uny \df{d}{d s}\qsxyn
\\&=\df12\sum_{y\in\Xnd}\uny\Dn\qsxyn
=\df12\Dn\sum_{y\in\Xnd}\uny\qsxyn
\\&=\df12\Dn\E \lbk u_0\lpa S^x_{\h}(s)\rpa\rbk=\mathscr{A}_n\sTns \unx.
\end{split}
\eeqs
So, we integrate \eqref{tder} by parts twice, and we observe that the boundary terms
always vanish at $\infty$ (as $s\nearrow\infty$) and that we have $(\pa/\pa s)\ptzs=0$ at
$s=0$ but $\ptzz>0$.  This gives us
\beq\lbl{fstp}
\bsp
\df{d}{d t}u^x_n(t)&=\df{d}{d t}\E \lbk u_0\lpa \S^x_{B,\h}(t)\rpa\rbk=
-\int_0^\infty\df{\pa}{\pa s}\ptzs\df{d}{d s}\sTns \unx ds
\\&=\ptzz\sAn\unx+
\int_0^\infty\ptzs\sAns\sTns\unx ds
\\&=\ptzz\sAn\unx+\sAns\int_0^\infty\ptzs\sTns\unx ds
\\&=\df{\Dn\unx}{\sqrt{8\pi t}}+\df18\Dns u^x_n(t).
\end{split}
\eeq
Obviously, $u^x_n(0)=\unx$, and we have proven \eqref{btrwdde}.

Of course, if $u^x_n(t)=\Ktxn$, then $\unx=\Kzxn$ as given in \lemref{fodde},
and we have by \eqref{btrwdsty} and by the steps above
\beqs
\bsp
\df{d}{d t}u^x_n(t)&=\df{d}{d t}\Ktxn=2\int_0^\infty\df{\pa}{\pa t}\ptzs\qsxn ds
=\int_0^\infty\df{\pas}{\pa s^2}\ptzs\qsxn ds
\\&=\df{\Dn\qzxn}{\sqrt{8\pi t}}+\df18\Dns u^x_n(t)
=\df{\Dn\Kzxn}{\sqrt{8\pi t}}+\df18\Dns u^x_n(t).
\end{split}
\eeqs
The proof is complete.
\epfs
\fi
\section{Limit solutions in the  Lipschitz case}\lbl{appB}
We now state prove our lattice-limit solution  existence, uniqueness, and regularity  for our BTBM SIE on $\Rp\times\Rd$ under Lipschitz conditions.
\bthm[Lattice-limits solutions: the Lipschitz case]\lbl{latlimlip}
Under the Lipschitz conditions 
there exists a unique-in-law direct $\beta$-ISLTRW SIE weak-limit solution to $\btsie$, $U$, such that
$U(t,x)$ is $L^p(\Omega,\P)$-bounded on $\T\times\Rd$ for every $p\ge2$ and $U\in\H^{\lpa{\tf{2\nu-d}{4\nu}}\rpa^{-},{\lpa\tf{4-d}{2}\wedge 1\rpa}^{-}}(\T\times\Rd;\R)$ for every $d=1,2,3$. 
\end{thm}
\thmref{latlimlip} follows as a corollary to the results of \secref{regandtight} combined with the following proposition.  
\bpr\lbl{lipprop}
Under the Lipschitz conditions $\eqref{lcnd}$ there exists a unique direct solution to $\isltrwsien$, $\utn$, on some filtered probability space $\OFFtP$ that is $L^p(\Omega,\P)$-bounded on $[0,T]\times \X_{n}^d$ for every $T>0$, $p\ge2$, $n\in\Ns, $and $d=1,2,3$.  
\epr
The proof of \propref{lipprop} follows the same steps as the non-discretization Picard-type direct proof of the corresponding part in the continuous case in \secref{pf1stmain}, with obvious changes, and we leave the details to the interested reader.
\ig
\bpf
For the existence proof, fix a lattice $\Xd$ on which we construct a solution iteratively.
So, given a collection of independent BMs $\lbr W^y(t)\rbr_{y\in\Xd}$,  define for each $(t,x)\in\Rp\times\Xd$
\beq\lbl{itdefdis}
\bc
\zivtx=\ds\sum_{y\in\Xd} \Ktxy u_0(y)&\cr
\ds\npoivtx=\zivtx
\ds+\sum_{y\in {\Xd}}\int_0^t \Ktsxy a(\nivsy)\df{dW^y(s)}{\sdd}&\cr
\ec
\eeq
Just as in the proof of \thmref{lip} (the BTBM SIE case) we show that, for any $p\ge2$ and all $d=1,2,3$, the sequence $\lbr \nivtx\rbr_{n\ge1}$
converges in $L^p(\Omega)$ to a solution.  Let
$$D_{n,p}(t,x):=\E\lab\npoivtx-\nivtx\rab^p$$
and
$$M_{n,p}(t):=\sup_{x}\dnptx.$$
Following the same steps as in the continuous-space BTBM SIE case we get
\beq\lbl{mdiffdis}
\bsp
\mnpt\le C \lpa t^{\tf{4-d}{4}}\rpa^{\tf{p-2}{2}}\intzt\mnmops\lbk t-s\rbk^{-\tf{d}{4}} ds
\end{split}
\eeq
where we have used \lemref{2ndQinequality}.  Again, exactly as in the proof of \thmref{lip},
we conclude that
there exists an  $L^p(\Omega)$-limit $\uttx$ that satisfies the $\beta$-ISLTRW SIE $\isltrwsien$ in \eqref{btrwsie}
with respect to the given BMs $\lbr W^y(t)\rbr_{y\in\Xd}$, and such that $\uttx$ is $L^p(\Omega)$
bounded on $\T\times\Xd$, for any $p\ge2$ and for any $T>0$.  Therefore, by \lemref{regtight} above, that
solution (and its spatially-linearly-interpolated version) is almost surely H\"older continuous in time with exponent $\gamma_t\in\lpa0,\tf{4-d}{8}\rpa$ for $d=1,2,3$.

Again, the uniqueness proof follows exactly the same steps as its counterpart in \thmref{lip} with the
space being $\Xd$ instead of $\Rd$ and with $\KBtxy$ replaced with $\Ktxy$.   So, let $d=1,2,3$
and let $\tilde{U}_1$ and $\tilde{U}_2$ be two solutions to $\isltrwsien$ for a fixed $n$ that are $L^2(\Omega)$-bounded
on $\T\times\Xd$,  for every $T>0$.  Fix an arbitrary $(t,x)\in\Rp\times\Xd$.
Let $D(t,x)=\tilde{U}_2^x(t)-\tilde{U}_1^x(t)$, $L_2(t,x)=\E D^2(t,x)$,
and $L^*_2(t)=\sup_{x\in\Xd}L_2(t,x)$ (which is bounded on $\T$ for every $T>0$ by hypothesis).
Then, using \eqref{btrwsie}, the Lipschitz condition in \eqref{lcnd},
and taking the supremum over the space variable and using \lemref{2ndQinequality}
we have
\beq\lbl{undis1}
\bsp
L_2(t,x)&=\sum_{y\in\Xd}\int_0^t \E\lbk a(\tilde{U}_2^y(s))- a(\tilde{U}_1^y(s))\rbk^2\lbk\Ktsxy\rbk^2ds
\\&\le C\intzt L^*_2(s)\sum_{y\in\Xd}\lbk\Ktsxy\rbk^2ds
\\&\le C\intzt\df{L^*_2(s)}{(t-s)^{\tf d4}}ds
\end{split}
\eeq
Iterating and interchanging the order of integration we get
\beq\lbl{undis3}
\bsp
L^*_2(t)\le C\lpa\ds\int_0^t L^*_2(s)ds \rpa
\end{split}
\eeq
for every $t\ge0$.  An easy application of Gronwall's lemma gives that $L^*_2\equiv0$.  So for every $(t,x)\in\Rp\times\Xd$ and $d=1,2,3$ we have
$\tilde{U}_1^x(t)=\tilde{U}_2^x(t)$ with probability one.  So that, using \lemref{regtight} (a), the spatially-linearly-interpolated
versions of $\tilde{U}_1^x(t)$ and $\tilde{U}_2^x(t)$ are indistinguishable---and hence pathwise uniqueness holds
for such interpolated solutions to $\isltrwsien$---on $\Rp\times\Rd$.
\epf
\fi
\bcr  \thmref{latlimlip} holds.
\ecr
\bpf  The conclusion follows from \propref{lipprop}, \lemref{sptdiff}, \lemref{tmpdiff}, and \lemref{regtight} (b).
\epf

\brm
With extra work, it is possible to prove the existence of a strong limit solution under Lipschitz conditions.  We plan to address that in a future article. 
\erm
\bcm
\bthm[Existence, uniqueness, and regularity of direct solutions for dimensions $d=1,2,3$]\lbl{lip}
Fix an arbitrary $T>0$, and let $\T=[0,T]$.  Assume that \eqref{lcnd} holds.
Then there exists a pathwise- unique strong solution $(U,\sW)$ to $\btsie$
on $\Rp\times\Rd$, for $d=1,2,3$, which is $L^p(\Omega)$-bounded 
on $\T\times\Rd$ for all $p\ge2$.  Furthermore, $U\in\H^{{\tf{4-d}{8}}^{-},\alpha_{d}^{-}}(\T\times\Rd;\R)$ for every $d=1,2,3$ and every $\alpha_{d}\in I_{d}$, where the intervals $I_{d}$ are given in \thmref{mainthm2}.
\ethm
\fi
\section{Glossary of frequently used acronyms and notations}\lbl{glossary}
\begin{enumerate}\renewcommand{\labelenumi}{\Roman{enumi}.}
\item {\textbf{Acronyms}}\vspace{2mm}
\bit
\item BM: Brownian motion
\item BTBM: Brownian-time Brownian motion. 
\item BTBM SIE: Brownian-time Brownian motion stochastic integral equation.
\item BTP:  Brownian-time processe.
\item BTP SIE: Brownian-time process stochastic integral equation.
\item BTC: Brownian-time chain.
\item BTRW: Brownian-time random walk.
\item $\beta$-ISLTRW DDE: Brownian-time random walk differential-difference equation.
\item $\beta$-ISLTRW SIE: Brownian-time random walk stochastic integral equation.
\item DDE: Differential difference equation.
\item KS: Kuramoto-Sivashinsky.
\item RW: Random walk.
\item SIE: Stochastic integral equation.
\eit
\vspace{2.5mm}
\item {\textbf{Notations}}\vspace{2mm}
\bit
\item $\N$: The usual set of natural numbers $\lbr1,2,3,\ldots\rbr$.
\item $\qtxyn$:  The $d$-dimensional continuous-time random walk transition density.
starting at $x\in\Xnd$ and going to $y\in\Xnd$ in time $t$.
\item $\psxy$: The density of a $d$-dimensional BM.
\item $\ptzs$: The  density of a $1$-dimensional BM, starting at $0$. 
\item $\KBtxy$: The kernel or density of a $d$-dimensional Brownian-time Brownian motion.\vspace{0.5 mm}
\item $\Ktxyn$: The kernel or density of a $d$-dimensional Brownian-time random walk on a spatial lattice with step size $\h$ in each of the $d$-dimensions.
\item $\btsie$: The BTBM SIE with diffusion coefficient $a$ and initial function $\un$.
\item $\isltrwsien$: The $\beta$-ISLTRW SIE on the lattice $\Xnd=\h\Zd$ with diffusion coefficient $a$ and initial function $\un$.
\eit
\end{enumerate}

\end{document}